\documentclass[11pt]{article}
\usepackage{ebezier}
\usepackage{amssymb,amsmath,mathrsfs,color} 
 \allowdisplaybreaks

\catcode`!=11
\let\!int\int \def\int{\displaystyle\!int}
\let\!lim\lim \def\lim{\displaystyle\!lim}
\let\!sum\sum \def\sum{\displaystyle\!sum}
\let\!sup\sup \def\sup{\displaystyle\!sup}
\let\!inf\inf \def\inf{\displaystyle\!inf}
\let\!cap\cap \def\cap{\displaystyle\!cap}
\let\!max\max \def\max{\displaystyle\!max}
\let\!min\min \def\min{\displaystyle\!min}
\let\!frac\frac \def\frac{\displaystyle\!frac}
\catcode`!=12

\newtheorem{theorem}{\bf Theorem}[section]
\newtheorem{lemma}{\bf Lemma}[section]
\newtheorem{definition}{\bf Definition}[section]
\newtheorem{proposition}{\bf Proposition}[section]
\newtheorem{remark}{\bf Remark}[section]

\def\Proof{\it{Proof.}\rm\quad}

\def\pd#1#2{\frac{\partial#1}{\partial#2}}

\catcode`@=11
\@addtoreset{equation}{section}
\catcode`@=12

\begin{document}

\title{Global Smooth Supersonic Flows in Infinite Expanding Nozzles}
\author{Chunpeng Wang\footnote{The research is supported by a grant from the National Natural Science Foundation of China (No. 11222106).}
\\
\small School of Mathematics, Jilin University, Changchun
130012, P. R. China
\\
\small The Institute of Mathematical Sciences, The Chinese
University of Hong Kong, Shatin, NT, Hong Kong
\\
\small (email: wangcp@jlu.edu.cn)
\\[2mm]
Zhouping Xin\footnote{The research is supported by Zheng Ge Ru Foundation, Hong Kong RGC Earmarked Research
Grants CUHK-4041/11P and CUHK-4048/13P, a grant from the Croucher Foundation, and a Focus Area Grant from the
Chinese University of Hong Kong.}
\\
\small The Institute of Mathematical Sciences and Department of
Mathematics,
\\
\small The Chinese University of Hong Kong, Shatin, NT, Hong Kong
\\
\small (email: zpxin@ims.cuhk.edu.hk)
}
\date{}
\maketitle

\begin{abstract}
This paper concerns smooth supersonic flows with Lipschitz continuous speed
in two-dimensional infinite expanding nozzles,
which are governed by a quasilinear hyperbolic equation being singular at the sonic and vacuum state.
The flow satisfies the slip condition on the walls and the flow velocity
is prescribed at the inlet.
First, it is proved that if the incoming flow is away from the sonic and vacuum state
and its streamlines are rarefactive at the inlet,
then a flow in a straight nozzle never approaches the sonic and vacuum state in any bounded region.
Furthermore, a sufficient and necessary condition
of the incoming flow at the inlet is derived
for the existence of a global smooth supersonic flow in a straight nozzle.
Then, it is shown that for each incoming flow satisfying this condition,
there exists uniquely a global smooth supersonic flow in a symmetric nozzle with convex upper wall.
It is noted that such a flow may contain a vacuum.
If there is a vacuum for a global smooth transonic flow in a symmetric nozzle with convex upper wall,
it is proved that for the symmetric upper part of the flow,
the first vacuum point along the symmetric axis
must be located at the upper wall
and the set of vacuum points is the closed domain
bounded by the tangent half-line of the upper wall at this point to downstream
and the upper wall after this point.
Moreover, the flow speed is globally Lipschitz continuous in the nozzle,
and on the boundary between the gas and the vacuum,
the flow velocity is along this boundary
and the normal derivatives of the flow speed and the square of the sound speed both are zero.
As an immediate consequence,
the local smooth transonic flow obtained in \cite{WX3}
can be extended into a global smooth transonic flow
in a symmetric nozzle whose upper wall after the local flow is convex.
\\[6pt]
\\
{\sl 2000 MR Subject Classification:} 76J20 76N10 35L20 35L81
\end{abstract}

\newpage

\tableofcontents

\newpage

\section{Introduction}

This paper concerns the global existence of smooth supersonic potential flows
in two-dimensional infinite expanding nozzles,
which arises naturally in physical experiments and engineering designs.
Here, a smooth supersonic flow means that
its velocity is continuous and its speed is Lipschitz continuous.
Assume that there is a smooth incoming supersonic flow at the inlet of an infinite expanding nozzle.
It is interesting to study whether the smooth supersonic flow exists globally or not.
For such a global smooth supersonic flow, what is the asymptotic behavior of the flow at the downstream,
is there a vacuum in the nozzle, what is the vacuum set
and what is the behavior of the flow
on the boundary between the gas and the vacuum?
If a smooth supersonic flow does not exist globally, what is the singularity?

One of the main motivations of this paper lies in
the global existence of smooth transonic flows
in two-dimensional infinite de Laval nozzles.
In \cite{WX3}, we studied smooth transonic potential flows of Meyer type
in a class of two-dimensional finite de Laval nozzles.
It was shown that a flow with nonexceptional points is unstable
for a $C^1$ small perturbation in the shape of the nozzle
and thus we sought smooth transonic flows of Meyer type whose sonic points are exceptional.
It was proved that there exists uniquely such a smooth transonic flow near the throat of the nozzle,
whose acceleration is Lipschitz continuous,
if the wall of the nozzle is sufficiently flat.
A natural question is whether this local smooth transonic flow
can be extended smoothly and globally.
For a supersonic flow, there is also a discontinuous extension pattern called a transonic shock,
which is supersonic upstream and turn to subsonic across the transonic shock.
Courant and Friedrichs had described such a pattern as follows (\cite{CF}).
Given the appropriately large receiver pressure, if the
upstream flow is still supersonic behind the throat of the nozzle, then at a certain place
in the divergent part of the nozzle a shock front intervenes and the gas is compressed
and slowed down to subsonic speed, which was proved rigidly
in recent works \cite{LXY,LXY1,LXY2} for two-dimensional full Euler system.
See also \cite{CKL,CFM,CFM1,CFM3,XY} for related works.
The reason of the intervention of a transonic shock in the pattern
lies in that there is a large receiver pressure, which can be achieved only when the flow is subsonic,
at the outlet of the nozzle.
We are interested in the smooth extension of a supersonic flow in an infinite expanding nozzle
when the flow is not required to be subsonic downstream.
That is to say, it is concerned how to exclude transonic shocks, supersonic shocks and sonic state
in the extension of a supersonic flow in an infinite expanding nozzle.

\vskip5mm
\hskip40mm
\setlength{\unitlength}{0.6mm}
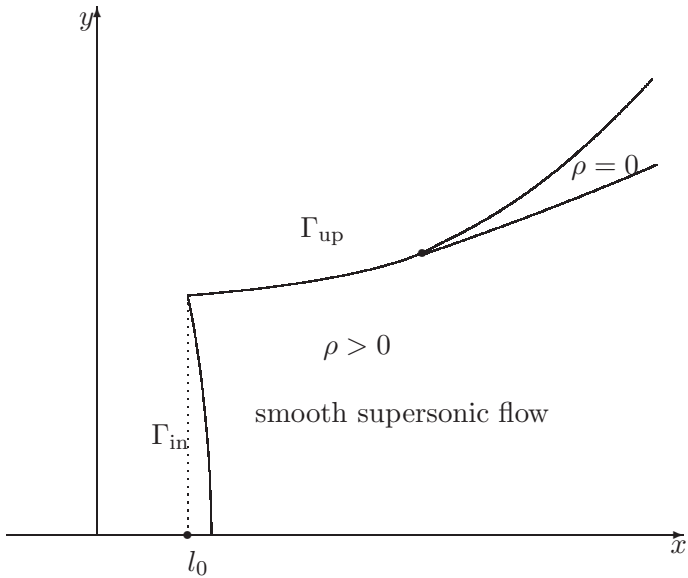
\begin{picture}(250,120)
\put(-20,0){\vector(1,0){150}}
\put(0,0){\vector(0,1){117}}
\put(127,-4){$x$} \put(-4,113){$y$}

\put(60,48){\cbezier(-40,5)(-20,6.5)(5,10)(15,16)}
\put(60,48){\cbezier(15,16)(25,21)(40,28)(63,53)}

\put(60,0){\cbezier(-40,53)(-37.5,42)(-34.8,20)(-34.7,0)}

\put(45,66){$\Gamma_{\text{\rm up}}$}
\put(12,20){$\Gamma_{\text{\rm in}}$}

\put(20,0){\qbezier[30](0,0)(0,26.5)(0,53)}
\put(20,0){\circle*{1.5}}
\put(20,-8){$l_0$}
\put(35,25){smooth supersonic flow}
\put(50,40){$\rho>0$}

\put(105,80){$\rho=0$}
\put(72,62.5){\circle*{1.5}}
\put(75,55){\qbezier(-3,7)(26,17)(49,27)}

\end{picture}
\vskip10mm

As usual, it is assumed that the nozzle is symmetric with respect to the $x$-axis.
For convenience, we consider only the upper part of the nozzle.
The upper wall and the inlet are given by $\Gamma_{\text{\rm up}}:y=f(x)\,(l_0\le x<l_1)$
and $\Gamma_{\text{\rm in}}:x=\Upsilon(y)\,(0\le y\le f(l_0))$, respectively, where
$0<l_0<l_1\le+\infty$, $f\in C^{2}([l_0,l_1))$,
$\Upsilon\in C^{2}([0,f(l_0)])$ and
\begin{align*}
f(l_0)>0,\quad
\lim_{x\to l_1^-}(x+f(x))=+\infty,\quad
\inf_{(l_0,l_1)}f'\ge0,\quad
\Upsilon(f(l_0))=l_0,\quad
\Upsilon'(0)=0,\quad\Upsilon'(f(l_0))=-f'(l_0).
\end{align*}
Consider the smooth supersonic flow problem in the infinite nozzle
$\Omega$ bounded by the upper wall $\Gamma_{\text{\rm up}}$, the lower wall $x$-axis
and the inlet $\Gamma_{\text{\rm in}}$.
The flow velocity is prescribed at the inlet.
In particular, only the flow speed needs to be prescribed
if it is assumed that the incoming flow velocity is along
the normal direction at the inlet.
The flow satisfies the slip condition on the walls.
Therefore, the problem can be formulated as follows
\begin{align}
\label{int-phy1}
&\mbox{div}(\rho(|\nabla\varphi|^2)\nabla\varphi)=0,\quad&&(x,y)\in\Omega,
\\
\label{int-phy2}
&\varphi(\Upsilon(y),y)=0,\quad&&0<y<f(l_0),
\\
\label{int-phy3}
&|\nabla\varphi(\Upsilon(y),y)|=q_0(y),\quad&&0<y<f(l_0),
\\
\label{int-phy4}
&\rho(|\nabla\varphi(x,0)|^2)\pd{\varphi}y(x,0)=0,\quad&&\Upsilon(0)<x<l_1,
\\
\label{int-phy5}
&\rho(|\nabla\varphi(x,f(x))|^2)\Big(\pd{\varphi}y(x,f(x))-f'(x)\pd{\varphi}x(x,f(x))\Big)=0,
\quad&&l_0<x<l_1,
\end{align}
where $\rho$ is a known function given by
\begin{align*}
\rho(q^2)=\Big(1-\frac{\gamma-1}{2}q^2\Big)^{1/(\gamma-1)},\quad
0\le q\le\sqrt{2/{(\gamma-1)}}\quad(\gamma>1),
\end{align*}
$q_0\in C^{1}([0,f(l_0)])$ and it satisfies
\begin{align*}
\sqrt{2/{(\gamma+1)}}<\inf_{(0,f(l_0))}q_0\le\sup_{(0,f(l_0))}q_0<\sqrt{2/{(\gamma-1)}},
\end{align*}
which guarantees that the flow is away from the sonic and vacuum state at the inlet.

It is noted that \eqref{int-phy1} is a quasilinear hyperbolic equation in the supersonic region
and is singular at sonic and vacuum points, where two characteristics coincide.
For the special case that the nozzle is straight,
it is shown that smooth supersonic flows
to the problem \eqref{int-phy1}--\eqref{int-phy5}
are away from vacuum in any bounded region.
That is to say, for a supersonic flow in a straight nozzle,
a vacuum never forms and thus the flow exists smoothly and globally
if there is not a sonic point and a shock.
It turns out that a natural choice to exclude the sonic state is that
the streamlines of the incoming flow are rarefactive at the inlet, i.e.
\begin{align}
\label{int-phy6-0}
\Upsilon''(y)\le0,\quad 0<y<f(l_0).
\end{align}
Indeed, the speed of a smooth supersonic flow in a straight nozzle
is not less than the minimal speed of the incoming flow on the inlet if \eqref{int-phy6-0} is satisfied.
Can \eqref{int-phy6-0} also exclude shocks for a supersonic flow in a straight nozzle?
The answer turns to be  no.
If \eqref{int-phy6-0} is satisfied, a sufficient and necessary condition
to exclude shocks for a supersonic flow in a straight nozzle
is proved to be
\begin{align}
\label{int-phy7}
|q'_0(y)|\le\frac{-\Upsilon''(y)}{1+(\Upsilon'(y))^2}\,
\sqrt{\frac{-q_0^2(y)\rho(q_0^2(y))}{\rho(q_0^2(y))+2q_0^2(y)\rho'(q_0^2(y))}}\,,
\quad 0<y<f(l_0),
\end{align}
which also implies \eqref{int-phy6-0}.
Therefore, if $q_0$ satisfies \eqref{int-phy7}, then
the problem \eqref{int-phy1}--\eqref{int-phy5} in a straight nozzle
admits uniquely a smooth supersonic flow which never approaches the sonic and vacuum state
in any bounded region;
further, the flow is shown to be accelerating along each streamline.
If \eqref{int-phy7} is invalid, a sonic point or a shock must form in the flow;
furthermore, if in addition, \eqref{int-phy6-0} holds, then
the flow is away from the sonic state
and its singularity is a shock.
It is noted that the local smooth transonic flow obtained in \cite{WX3}
satisfies \eqref{int-phy7} at the outlet.
So, this local flow can be extended into a global smooth transonic flow without vacuum
if the upper wall of the nozzle after the outlet of the local flow is straight.

For a straight nozzle, \eqref{int-phy7} is a sufficient and necessary condition
to get a global smooth supersonic flow.
Thus, for a general expanding nozzle,
we solve the problem \eqref{int-phy1}--\eqref{int-phy5}
under the assumption \eqref{int-phy7}.
Generally speaking, the geometry of nozzles can create vacuum, sonic state and shocks.
It is shown that for a nozzle whose upper wall is convex,
the problem \eqref{int-phy1}--\eqref{int-phy5} admits uniquely a global smooth supersonic flow.
The convexity of the upper wall of the nozzle is necessary in the following sense:
if the upper wall of the nozzle is a nonconvex perturbation of a straight line,
then there exists a $q_0$ satisfying \eqref{int-phy7}
such that a shock must form for the supersonic flow in the nozzle.
Such a global smooth supersonic flow in a nozzle with convex upper wall
may contain a vacuum,
and its acceleration is continuous in the gas region.
If there is a vacuum,
it is clear that the set of vacuum points is closed
and thus there exists the first vacuum point in the increasing $x$-direction.
It is shown that this first vacuum point must be located at the upper wall
and the set of vacuum points is the closed domain
bounded by the tangent half-line of the upper wall at this point to the downstream
and the upper wall after this point.
Moreover, this tangent half-line is a streamline of the flow.
That is to say, the flow velocity at the boundary between the gas and the vacuum
is along this boundary.
For such a global smooth supersonic flow with vacuum,
the flow speed is globally Lipschitz continuous in the nozzle
and the normal derivatives of the flow speed and the square of the sound speed
on the boundary between the gas and the vacuum both are zero.
Therefore, the vacuum in the flow is not the so called physical vacuum in \cite{YL1,YL2}.
This is natural since there is not an external force on the flow.

The global existence result for smooth supersonic flows in a nozzle with convex upper wall
shows that the geometry of the nozzle
is essential for the formation of vacuum.
For a nozzle whose upper wall is convex but not straight,
if the incoming flow is near the vacuum, then
a vacuum must form in the global smooth supersonic flow.
Furthermore, the first vacuum point will tend to $(x_*,f(x_*))$ with
$x_*=\inf\{x\in[l_0,l_1):f''(x)>0\}$ as $q_0(f(l_0))$ tends to $\sqrt{2/{(\gamma-1)}}$.
It is noted that there are a large class of nozzles
in which a global smooth supersonic flow must contain a vacuum
even if the incoming flow is near the sonic state.
For example, if
\begin{align*}
\left\{
\begin{aligned}
&{\lim_{x\to+\infty}}f''(x)x^{2\gamma/(\gamma+1)}=+\infty,
\quad&&\mbox{ when } l_1=+\infty \mbox{ and }f'(+\infty)<+\infty,
\\
&{\lim_{x\to+l_1}}\frac{f''(x)}{(f'(x))^3}f^{2\gamma/(\gamma+1)}(x)=+\infty,
\quad&&\mbox{ when } l_1\le+\infty \mbox{ and }f'(l_1)=+\infty,
\end{aligned}
\right.
\end{align*}
then each global smooth supersonic flow to the problem \eqref{int-phy1}--\eqref{int-phy5}
admits a vacuum provided that $q_0$ satisfies
\begin{align*}
|q'_0(y)|<\frac{-\Upsilon''(y)}{1+(\Upsilon'(y))^2}\,
\sqrt{\frac{-q_0^2(y)\rho(q_0^2(y))}{\rho(q_0^2(y))+2q_0^2(y)\rho'(q_0^2(y))}}\,,
\quad 0\le y\le f(l_0).
\end{align*}

Before the appearance of a vacuum, the smooth supersonic flow problem is
solved in the potential plane, where the domain is a rectangle
and thus it is convenient to calculate the effect of the boundary condition.
For a global smooth supersonic flow in a straight nozzle,
a necessary condition of the incoming flow at the inlet is derived
by a characteristic method.
The method for the existence of smooth supersonic flows before vacuum formation
is a fixed point argument.
We use the method of characteristics to get
a priori estimates for suitable semilinear problems and then use the fixed point theorem
to get a solution to the quasilinear problem.
However, the coordinates transformation
between the physical coordinates and the potential-stream coordinates
may be invalid in the appearance of vacuum.
So, we solve the smooth supersonic flow problem with vacuum as
a free boundary problem in the physical plane,
where the boundary between the gas and the vacuum is free.
It is proved that this free boundary is a streamline of the flow
by studying the extension of the gas
from the potential level set where the first vacuum forms.

The paper is arranged as follows.
In $\S\, 2$ we formulate the global smooth supersonic flow problem in an infinite expanding nozzle
with convex upper wall. Moreover,
a necessary condition of the incoming flow at the inlet is derived
for the existence of global smooth supersonic flows in a straight nozzle.
In $\S\, 3$ we solve the maximal smooth supersonic flow problem
before vacuum formation in a nozzle with convex upper wall.
In particular, there exists uniquely a global smooth supersonic flow without vacuum
if the nozzle is straight. Properties of smooth supersonic flows before vacuum formation
are also studied.
In $\S\, 4$ we investigate properties of the first vacuum point
for a smooth supersonic flow and give some sufficient conditions for formation of a vacuum.
Subsequently, we prove the global well-posedness
for the smooth supersonic flow problem in an infinite nozzle with convex upper wall in $\S\, 5$.

\section{Formulation of the smooth supersonic flow problem}

In this section, we formulate the smooth supersonic flow problem in an infinite expanding nozzle.
The focus is to characterize the geometry of the nozzles
and the boundary conditions at their inlets to ensure the global existence
of smooth supersonic flows in these nozzles.
One of the main motivations is to extend the local smooth transonic flows
in a de Laval nozzle obtained in \cite{WX3} globally.

\subsection{Governing equations}

Steady isentropic compressible Euler flows satisfy
\begin{align}
\label{euler-1}
&\pd{}{x}(\rho u)+\pd{}{y}(\rho v)=0,
\\
\label{euler-2}
&\pd{}{x}(P+\rho u^2)+\pd{}{y}(\rho uv)=0,
\\
\label{euler-3}
&\pd{}{x}(\rho uv)+\pd{}{y}(P+\rho v^2)=0,
\end{align}
where $(u,v)$, $P$ and $\rho$ represent the velocity, pressure and
density of the flow, respectively. The flow is assumed to be
isentropic so that $P=P(\rho)$ is a smooth function. In
particular, for a polytropic gas with the adiabatic exponent
$\gamma>1$,
\begin{align}
\label{euler-4}
P(\rho)=\frac1\gamma\rho^\gamma
\end{align}
is the normalized pressure. Assume further that the flow is
irrotational, i.e.
\begin{align}
\label{euler-5}
\pd{u}{y}=\pd{v}{x}.
\end{align}
Then the density $\rho$ is expressed in terms of the speed
$q$ according to the Bernoulli law (\cite{CF})
\begin{align}
\label{Bernoulli}
\rho(q^2)=\Big(1-\frac{\gamma-1}{2}q^2\Big)^{1/(\gamma-1)},\quad
q=\sqrt{u^2+v^2},\quad0\le q\le c^*=\sqrt{2/{(\gamma-1)}}.
\end{align}
The sound speed $c$ is defined as
$$
c^2=P'(\rho)=\rho^{\gamma-1}=1-\frac{\gamma-1}{2}q^2,\quad0\le q\le c^*.
$$
At the sonic state, the sound speed is
$c_*=\sqrt{2/{(\gamma+1)}}$,
which is called the sonic speed in the sense that
the flow is subsonic ($q<c$) when $0<q<c_*$, sonic ($q=c$) when $q=c_*$
and supersonic ($q>c$) when $q>c_*$.

Define a velocity potential $\varphi$ and a stream
function $\psi$, respectively, by
$$
\pd\varphi x=u,\quad\pd\varphi y=v, \quad\pd\psi x=-\rho
v,\quad\pd\psi y=\rho u,
$$
which are
$$
\pd\varphi x=q\cos\theta,\quad\pd\varphi y=q\sin\theta,\quad
\pd\psi x=-\rho q\sin\theta,\quad\pd\psi y=\rho q\cos\theta
$$
in terms of polar coordinators in the velocity space,
where $\theta$, which is called a flow angle,
is the angle of the velocity inclination to
the $x$-axis. Direct calculations show that the system
\eqref{euler-1}--\eqref{Bernoulli} can be reduced to the full potential equation
\begin{align}
\label{potentialequation}
\mbox{div}(\rho(|\nabla\varphi|^2)\nabla\varphi)=0
\end{align}
in the physical coordinates $(x,y)$
and the Chaplygin equations
\begin{align}
\label{Chaplygin}
\pd\theta\psi+\frac{\rho(q^2)+2q^2\rho'(q^2)}{q\rho^2(q^2)}
\pd{q}\varphi=0,
\quad\frac{1}{q}\pd{q}\psi-\frac{1}{\rho(q^2)}\pd\theta\varphi=0
\end{align}
in the potential-stream coordinates $(\varphi,\psi)$. Note that
$$
\pd{(\varphi,\psi)}{(x,y)}=\pd\varphi x\pd\psi y-
\pd\varphi y\pd\psi x=\rho q^2.
$$
So the coordinates transformation
between the two coordinate systems
is valid at least in the absence of stagnation and vacuum points. Eliminating
$\theta$ from \eqref{Chaplygin}
yields the following second-order quasilinear equation
\begin{align}
\label{eq2or}
\pd{^2A(q)}{\varphi^2}+\pd{^2B(q)}{\psi^2}=0
\end{align}
with
\begin{align*}
A(q)=\int_{c_*}^q\frac{\rho(s^2)+2s^2\rho'(s^2)}{s\rho^2(s^2)}ds,\quad
B(q)=\int_{c_*}^q\frac{\rho(s^2)}{s}ds, \qquad 0<q<c^*.
\end{align*}
Here, $B(\cdot)$ is strictly increasing in $(0,c^*)$,
while $A(\cdot)$ is strictly increasing in $(0,c_*]$ and
strictly decreasing in $[c_*,c^*)$.
It can be checked easily that both \eqref{potentialequation} and
\eqref{eq2or} are elliptic in the subsonic region $(q<c_*)$
and hyperbolic in the supersonic region $(q>c_*)$,
while singular at the sonic state $(q=c_*)$ and the vacuum $(q=c^*)$.

As shown in \cite{WX3}, in the supersonic region without vacuum, \eqref{eq2or} can be rewritten as
\begin{align*}
Q_{\varphi\varphi}-(b(Q)Q_{\psi})_\psi=0,
\end{align*}
or equivalently as
$$
\left\{
\begin{aligned}
&W_\varphi+b^{1/2}(Q)W_\psi
=\frac14b^{-1}(Q)p(Q)W(W+Z),
\\
&Z_\varphi-b^{1/2}(Q)Z_\psi
=-\frac14b^{-1}(Q)p(Q)Z(W+Z),
\end{aligned}
\right.
$$
where
\begin{align*}
Q=A(q),\quad
W=Q_\varphi-b^{1/2}(Q)Q_\psi,\quad
Z=-Q_\varphi-b^{1/2}(Q)Q_\psi
\end{align*}
and
\begin{gather*}
b(s)=\Big(\frac{\gamma+1}{2}q^2-1\Big)^{-1}
\Big(1-\frac{\gamma-1}{2}q^2\Big)^{2/(\gamma-1)+1}\Big|_{q=A^{-1}_+(s)}>0,
\quad s<0,
\\
p(s)=(\gamma+1)q^4\Big(\frac{\gamma+1}{2}q^2-1\Big)^{-3}
\Big(1-\frac{\gamma-1}{2}q^2\Big)^{3/(\gamma-1)+1}
\Big|_{q=A^{-1}_+(s)}>0,\quad s<0
\end{gather*}
with $A^{-1}_+$ being the inverse function of $A(\cdot)$ lying in $[c_*,c^*)$.
It is noted that
\begin{align}
\label{a8-8}
\lim_{s\to 0^-}s b^{-1}(s)p(s)=-\frac12,\quad
\lim_{s\to-\infty}s b^{-1}(s)p(s)=-(\gamma+1).
\end{align}

\subsection{Smooth transonic flows in finite de Laval nozzles}

Since one of the main motivations of this paper lies in
the global existence of smooth transonic flows,
we first recall the smooth transonic flow problem in finite de Laval nozzles
investigated in \cite{WX3} and show boundary conditions of
such local smooth transonic flows at an outlet.
This will yield the formulation of inlet boundary conditions
for smooth supersonic flows in an expanding nozzle.

It was shown in \cite{WX3} that a smooth transonic flow of Meyer type
with nonexceptional points is unstable
for a $C^1$ small perturbation on the wall (even if the wall is still smooth)
and thus we sought ones whose sonic points are exceptional.
For such a flow, its sonic curve must be located at the throat of the nozzle
and where the potential becomes a constant identically.
The smooth transonic flow problem in \cite{WX3} was described physically as follows:
for the upper wall $\Gamma_{\text{\rm up}}:y=f(x)\,(l_-\le x\le l_+)$
and for the inlet $\Gamma_{\text{\rm in}}:x=g(y)\,(0\le y\le f(l_-),g(f(l_-))=l_-)$,
seek a smooth transonic flow of Meyer type in the nozzle $\Omega$
bounded by the upper wall $\Gamma_{\text{\rm up}}$, the lower wall $x$-axis,
the inlet $\Gamma_{\text{\rm in}}$
and the outlet $\Gamma_{\text{\rm out}}:x=\Upsilon(y)\,(0\le y\le f(l_+),
\Upsilon(f(l_+))=l_+)$
which is a free boundary,
such that its velocity is along the normal direction at the inlet and the outlet,
it satisfies the slip condition on the walls
and its sonic curve is located at the throat where
the potential is equal to zero identically.
Mathematically, the problem was formulated as
\begin{align}
\label{phytran1}
&\mbox{div}(\rho(|\nabla\varphi|^2)\nabla\varphi)=0,\quad&&(x,y)\in\Omega,
\\
\label{phytran2}
&\varphi(g(y),y)=C_{\text{\rm in}},\quad&&0<y<f(l_-),
\\
\label{phytran3}
&\pd{\varphi}y(x,0)=0,\quad&&g(0)<x<\Upsilon(0),
\\
\label{phytran4}
&\pd{\varphi}y(x,f(x))-f'(x)\pd{\varphi}x(x,f(x))=0,\quad&&l_-<x<l_+,
\\
\label{phytran5}
&\varphi(\Upsilon(y),y)=C_{\text{\rm out}},\quad&&0<y<f(l_+),
\\
\label{phytran6}
&\varphi(0,y)=0,\,|\nabla\varphi(0,y)|=c_*,\quad&&0<y<f(0),
\\
\label{phytran7}
&|\nabla\varphi(x,y)|<c_*,\quad&&(x,y)\in\Omega\mbox{ and }x<0,
\\
\label{phytran8}
&|\nabla\varphi(x,y)|>c_*,\quad&&(x,y)\in\Omega\mbox{ and }x>0,
\end{align}
where $C_{\text{\rm in}}$ and $C_{\text{\rm out}}$
are free constants, and $\Gamma_{\text{\rm out}}$ is a free boundary.

\vskip1mm
\hskip8cm
\setlength{\unitlength}{0.6mm}
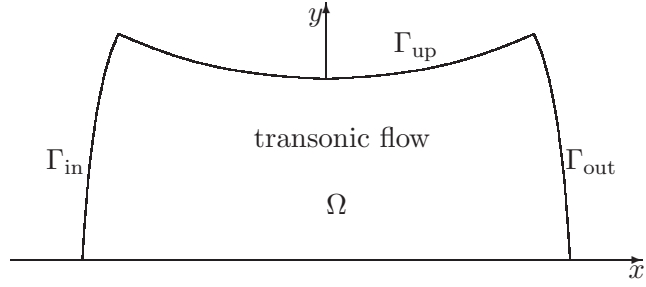
\begin{picture}(250,65)
\put(-70,0){\vector(1,0){140}}
\put(0,40){\vector(0,1){17}}
\put(67,-4){$x$} \put(-4,53){$y$}
\put(0,10){$\Omega$}

\put(-60,48){\cbezier(14,2)(30,-5)(40,-7)(60,-8)}
\put(60,48){\cbezier(-14,2)(-30,-5)(-40,-7)(-60,-8)}

\put(-40,48){\qbezier(-6,2)(-12,-10)(-14,-48)}

\put(40,48){\qbezier(6,2)(12,-10)(14,-48)}

\put(15,46){$\Gamma_{\text{\rm up}}$}
\put(-62,20){$\Gamma_{\text{\rm in}}$}
\put(53,20){$\Gamma_{\text{out}}$}
\put(-16,25){transonic flow}
\end{picture}
\vskip6mm

The smooth transonic flow problem was solved locally in the potential plane in \cite{WX3}.
It follows from \eqref{phytran6} that the mass flux of the flow is
$m=f(0) c_*^{1+2/(\gamma-1)}$.
Denote the flow speed at the inlet
and at the upper wall by
$$
q(g(y),y)={\mathscr Q}_{\text{\rm in}}(y),\quad0\le y\le f(l_{-})
\quad\mbox{ and }\quad
q(x,f(x))={\mathscr Q}_{\text{\rm up}}(x),\quad l_{-}\le x\le l_+,
$$
respectively. The stream function at the inlet
and the potential function at the upper wall can be given by
\begin{gather*}
\psi(g(y),y)=\Psi_{\text{\rm in}}(y)=\int_{0}^{y}
{{\mathscr Q}_{\text{\rm in}}(s)\rho({\mathscr Q}_{\text{\rm in}}^2(s))\sqrt{1+(g'(s))^2}}ds,\quad
0\le y\le f(l_{-})
\end{gather*}
and
\begin{gather*}
\varphi(x,f(x))=\Phi_{\text{\rm up}}(x)=\int_{0}^{x}
{{\mathscr Q}_{\text{\rm up}}(s)\sqrt{1+(f'(s))^2}}ds,\quad  l_{-}\le x\le l_+,
\end{gather*}
respectively.
Then, in the potential plane, the problem \eqref{phytran1}--\eqref{phytran8}
is changed to
\begin{align}
\label{p-eq}
&\pd{^2A(q)}{\varphi^2}+\pd{^2B(q)}{\psi^2}=0,
\quad&&(\varphi,\psi)\in(\zeta_{-},\zeta_+)\times(0,m),
\\
\label{p-inbc}
&\pd{A(q)}{\varphi}(\zeta_-,\psi)=
\frac{g''(y)(1+(g'(y))^2)^{-3/2}}{{\mathscr Q}_{\text{\rm in}}(y)\rho({\mathscr Q}_{\text{\rm in}}^2(y))}
\Big|_{y=Y_{\text{\rm in}}(\psi)},
\quad&&\psi\in(0,m),
\\
\label{p-lbbc}
&\pd{q}{\psi}(\varphi,0)=0,
\quad&&\varphi\in(\zeta_{-},\zeta_+),
\\
\label{p-ubbc}
&\pd{B(q)}{\psi}(\varphi,m)=
\frac{f''(x)(1+(f'(x))^2)^{-3/2}}{{\mathscr Q}_{\text{\rm up}}(x)}
\Big|_{x=X_{\text{\rm up}}(\varphi)},
\quad&&\varphi\in(\zeta_{-},\zeta_+),
\\
\label{p-outbc}
&q(0,\psi)=c_*,\quad&&\psi\in(0,m),
\\
\label{p-outbc1}
&0<q(\varphi,\psi)<c_*,\quad&&(\varphi,\psi)\in(\zeta_{-},0)\times(0,m),
\\
\label{p-outb2}
&c_*<q(\varphi,\psi)<c^*,\quad&&(\varphi,\psi)\in(0,\zeta_+)\times(0,m),
\\
\label{Qinq}
&{\mathscr Q}_{\text{\rm in}}(y)=q(0,\psi)\Big|_{\psi=\Psi_{\text{\rm in}}(y)},
\quad&&y\in(0,f_{-}(l_{-})),
\\
\label{Qubq}
&{\mathscr Q}_{\text{\rm up}}(x)=q(\varphi,m)\Big|_{\varphi=\Phi_{\text{\rm up}}(x)},
\quad&&x\in(l_{-},l_+),
\end{align}
where $\zeta_{\pm}=\Phi_{\text{\rm up}}(l_{\pm})$,
$Y_{\text{\rm in}}$ and $X_{\text{\rm up}}$ are the inverse functions of
$\Psi_{\text{\rm in}}$ and $\Phi_{\text{\rm up}}$, respectively.

The main results in \cite{WX3} are the following existence and uniqueness.

\begin{theorem}
\cite{WX3}
\label{theoremm5}
Assume that $f\in C^{4}([l_-,l_+])\,(-1\le l_-<0<l_+\le1)$ satisfies $f(0)>0$,
\begin{gather*}
\delta_{1,-}(-x)^{\lambda_{-}}\le f''(x)\le\delta_{2,-}(-x)^{\lambda_{-}},
\quad
|f^{(3)}(x)|\le\delta_{3,-}(-x)^{\lambda_-/4+1/2},
\quad
|f^{(4)}(x)|\le\delta_{4,-},\quad l_-\le x\le0,
\\
\delta_{1,+}x^{\lambda_{+}}\le f''(x)\le\delta_{2,+}x^{\lambda_{+}},
\quad |f^{(3)}(x)|\le\delta_{3,+}x^{\lambda_+-1},
\quad0\le x\le{l_+},
\end{gather*}
and $g\in C^{3,\alpha}([0,f(l_-)])\,(0<\alpha<1)$ satisfies
\begin{gather*}
g'(0)=0,\quad
g(f(l_-))=l_-,\quad
g'(f(l_-))=-f'(l_-),
\quad
\sup_{(0,f(l_-))}\Big|\Big(\frac{g'}{\sqrt{1+(g')^2}}\Big)''\Big|
\le\epsilon(l_-)(-l_-)^{3\lambda_{-}/2},
\\
\frac{1}{2}\le\frac{f(l_-)\sqrt{(f'(l_-))^2+1}}{-f'(l_-)}\Big(\frac{g'(y)}{\sqrt{1+(g'(y))^2}}\Big)'
\le\frac{3}{2},
\quad 0\le y\le f(l_-),
\end{gather*}
where $\lambda_{\pm}$, $\delta_{i,\pm}\,(i=1,2,3)$ and $\delta_{4,-}$
are positive constants
such that $\lambda_{\pm}>2$ and $\delta_{1,\pm}\le\delta_{2,\pm}$,
while $0<\epsilon(s)\le1$ for $-1\le s<0$ and $\lim_{s\to0^-}\epsilon(s)=0$.
Then, there exist
two positive constants $\delta_{0,-}$ (depending only on $\gamma$, $m$,
$\lambda_-$, $\delta_{1,-}$, $\delta_{2,-}$ and $\epsilon(\cdot)$)
and $\delta_{0,+}$ (depending only on $\gamma$, $m$,
$\lambda_+$, $\delta_{1,+}$, $\delta_{2,+}$ and $\delta_{3,+}$)
such that if $-\delta_{0,-}\le l_-<0<l_+\le\delta_{0,+}$,
then the problem \eqref{p-eq}--\eqref{Qubq}
admits a unique solution $q\in C^{1,1}([\zeta_{-},\zeta_+]\times[0,m])$ satisfying

{\rm(i)} On $[\zeta_{-},0]\times[0,m]$, $q$ satisfies
\begin{gather*}
c_*-C_{2,-}(-\varphi)^{{\lambda_-}/2+1}\le q(\varphi,\psi)
\le c_*-C_{1,-}(-\varphi)^{{\lambda_-}/2+1},
\\
\Big|\pd{q}\varphi(\varphi,\psi)\Big|\le C_{3,-}(-\varphi)^{{\lambda_-}/4+1/2},\quad
\Big|\pd{q}\psi(\varphi,\psi)\Big|\le C_{3,-}(-\varphi)^{{\lambda_-}/2+1},
\\
\Big|\pd{^2q}{\varphi^2}(\varphi,\psi)\Big|\le C_{4,-},\quad
\Big|\pd{^2q}{\varphi\partial\psi}(\varphi,\psi)\Big|\le C_{4,-}(-\varphi)^{{\lambda_-}/4+1/2},
\quad\Big|\pd{^2q}{\psi^2}(\varphi,\psi)\Big|\le C_{4,-}(-\varphi)^{{\lambda_-}/2+1},
\end{gather*}
where $C_{i,-}\,(i=1,2,3,4)\,(C_{1,-}\le C_{2,-})$
are positive constants depending only on $\gamma$, $m$,
$\lambda_-$, $\delta_{1,-}$ and $\delta_{2,-}$,
and $C_{3,-}$ also on $\delta_{3,-}$, while
$C_{4,-}$ also on $\delta_{3,-}$, $\delta_{4,-}$ and $\alpha$.

{\rm(ii)} On $[0,\zeta_{+}]\times[0,m]$, $Q=A(q)$ satisfies
\begin{gather*}
-C_{2,+}\varphi^{\lambda_++2}\le Q(\varphi,\psi)\le-C_{1,+}\varphi^{\lambda_++2},
\\
-C_{2,+}\varphi^{\lambda_++1}\le Q_\varphi(\varphi,\psi)
\le-C_{1,+}\varphi^{\lambda_++1},
\quad
|Q_\psi(\varphi,\psi)|\le C_{2,+}\varphi^{3{\lambda_+}/2+1},
\\
|Q_{\varphi\varphi}(\varphi,\psi)|\le C_{3,+}\varphi^{\lambda_+},
\quad
|Q_{\varphi\psi}(\varphi,\psi)|\le C_{3,+}\varphi^{5{\lambda_+}/4+1/2},
\quad
|Q_{\psi\psi}(\varphi,\psi)|\le C_{3,+}\varphi^{3{\lambda_+}/2+1},
\end{gather*}
where $C_{i,+}\,(i=1,2,3)\,(C_{1,+}\le C_{2,+})$
are positive constants depending only on  $\gamma$, $m$,
$\lambda_+$, $\delta_{1,+}$ and $\delta_{2,+}$,
while $C_{3,+}$ also on $\delta_{3,+}$.
\end{theorem}

\begin{theorem}
\cite{WX3}
\label{theoremm6}
Under the assumptions of Theorem \ref{theoremm5},
the transonic flow problem \eqref{phytran1}--\eqref{phytran8}
admits a unique smooth solution $\varphi\in C^{2,1}(\overline\Omega)$,
which satisfies Theorem \ref{theoremm5} (i)
for $q=|\nabla\varphi|$ in $\big\{(x,y)\in\Omega:x<0\big\}$,
while satisfies Theorem \ref{theoremm5} (ii)
for $Q=A(|\nabla\varphi|)$ in $\big\{(x,y)\in\Omega:x>0\big\}$.
\end{theorem}

We now characterize the outlet $\Gamma_{\text{\rm out}}:x=\Upsilon(y)\,(0\le y\le f(l_+))$,
which is a free boundary, and the conditions for the transonic flow at it.
If $l_+$ is small, Theorem \ref{theoremm5} shows that at the outlet of the de Laval nozzle,
the smooth transonic flow satisfies $Q(\zeta_+,\cdot)\in C^{1,1}([0,m])$ and
\begin{align*}
Q(\zeta_+,\psi)<0,\quad
-Q_\varphi(\zeta_+,\psi)
-b^{1/2}(Q(\zeta_+,\psi))|Q_\psi(\zeta_+,\psi)|>0,\quad 0\le\psi\le m.
\end{align*}
Therefore, the outlet satisfies
$\Upsilon\in C^{2,1}([0,f(l_+)])$ and
\begin{align}
\label{xc1-1}
\Upsilon(f(l_+))=l_+,\quad\Upsilon'(0)=0,\quad\Upsilon'(f(l_+))=-f'(l_+),\quad
\sup_{(0,f(l_+))}\Upsilon''<0,
\end{align}
and the flow speed at the outlet in the potential plane satisfies
$q_0\in C^{1,1}([0,f(l_+)])$ and
\begin{align}
\label{xc1-2}
c_*<q_0(y)<c^*,
\quad
|q'_0(y)|<\frac{-\Upsilon''(y)}{1+(\Upsilon'(y))^2}\,
\sqrt{\frac{-q_0^2(y)\rho(q_0^2(y))}{\rho(q_0^2(y))+2q_0^2(y)\rho'(q_0^2(y))}}\,,
\quad 0\le y\le f(l_+)
\end{align}
with
$$
q_0(y)=|\nabla\varphi(\Upsilon(y),y)|,\quad y\in[0,f(l_+)].
$$

\subsection{Formulation of the smooth supersonic flow problem}

Assume that $f\in C^{2}([l_0,l_1))\,(0<l_0<l_1\le+\infty)$ satisfies
\begin{align}
\label{sss-0}
f(l_0)>0,\quad\lim_{x\to l_1^-}(x+f(x))=+\infty,\quad
f'(l_0)\ge0,\quad
\inf_{(l_0,l_1)}f''\ge0.
\end{align}
Give an inlet $\Gamma_{\text{\rm in}}:x=\Upsilon(y)\,(0\le y\le f(l_0),\Upsilon(f(l_0))=l_0)$.
We study the smooth supersonic flow problem in the unbounded nozzle $\Omega$ bounded by
the upper wall $\Gamma_{\text{\rm up}}:y=f(x)\,(l_0\le x<l_1)$,
the lower wall $x$-axis
and the inlet $\Gamma_{\text{\rm in}}:x=\Upsilon(y)\,(0\le y\le f(l_0))$.
At the inlet, the flow velocity is prescribed
and it is assumed to be along the normal direction of the inlet for convenience.
On the walls, the flow satisfies the slip condition.
Here, a smooth supersonic flow means that
its velocity is continuous and its speed is Lipschitz continuous.

\vskip5mm
\hskip40mm
\setlength{\unitlength}{0.6mm}
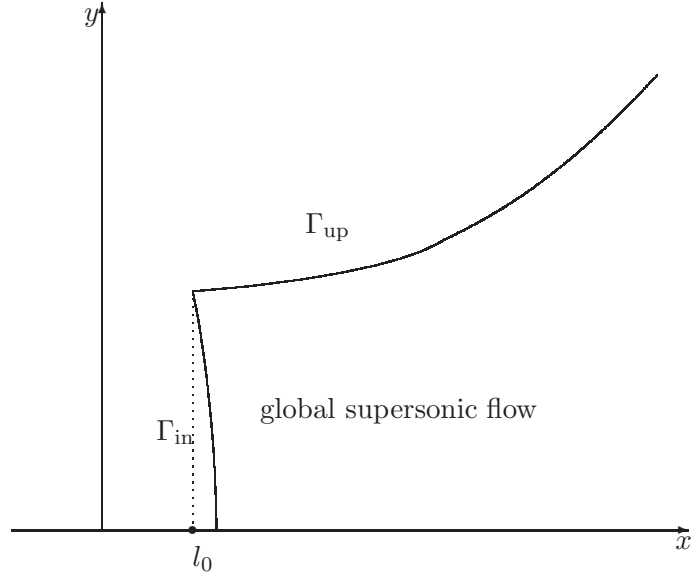
\begin{picture}(250,120)
\put(-20,0){\vector(1,0){150}}
\put(0,0){\vector(0,1){117}}
\put(127,-4){$x$} \put(-4,113){$y$}

\put(60,48){\cbezier(-40,5)(-20,6.5)(5,10)(15,16)}
\put(60,48){\cbezier(15,16)(25,21)(40,28)(63,53)}

\put(60,0){\cbezier(-40,53)(-37.5,42)(-34.8,20)(-34.7,0)}

\put(45,66){$\Gamma_{\text{\rm up}}$}
\put(12,20){$\Gamma_{\text{\rm in}}$}

\put(20,0){\qbezier[30](0,0)(0,26.5)(0,53)}
\put(20,0){\circle*{1.5}}
\put(20,-8){$l_0$}
\put(35,25){global supersonic flow}

\end{picture}
\vskip10mm

The smooth supersonic flow problem can be formulated as follows
\begin{align}
\label{phy1}
&\mbox{div}(\rho(|\nabla\varphi|^2)\nabla\varphi)=0,\quad&&(x,y)\in\Omega,
\\
\label{phy2}
&\varphi(\Upsilon(y),y)=0,\quad&&0<y<f(l_0),
\\
\label{phy3}
&|\nabla\varphi(\Upsilon(y),y)|=q_0(y),\quad&&0<y<f(l_0),
\\
\label{phy4}
&\rho(|\nabla\varphi(x,0)|^2)\pd{\varphi}y(x,0)=0,\quad&&\Upsilon(0)<x<l_1,
\\
\label{phy5}
&\rho(|\nabla\varphi(x,f(x))|^2)\Big(\pd{\varphi}y(x,f(x))-f'(x)\pd{\varphi}x(x,f(x))\Big)=0,
\quad&&l_0<x<l_1,
\end{align}
where $q_0$ is a given function defined on $[0,f(l_0)]$.
Motivated by conditions \eqref{xc1-1} and \eqref{xc1-2}
for local smooth transonic flows obtained in Theorem \ref{theoremm6},
we will solve the problem \eqref{phy1}--\eqref{phy5}
under the assumptions that
$\Upsilon\in C^{2}([0,f(l_0)])$ satisfies
\begin{align}
\label{xc4}
\Upsilon(f(l_0))=l_0,\quad
\Upsilon'(0)=0,\quad\Upsilon'(f(l_0))=-f'(l_0),\quad
\sup_{(0,f(l_0))}\Upsilon''\le0,
\end{align}
and $q_0\in C^{1}([0,f(l_0)])$ satisfies
\begin{gather}
\label{phy-8-q}
c_1\le q_0(y)\le c_2,
\quad 0<y<f(l_0),
\\
\label{phy-8}
|q'_0(y)|\le\frac{-\Upsilon''(y)}{1+(\Upsilon'(y))^2}\,
\sqrt{\frac{-q_0^2(y)\rho(q_0^2(y))}{\rho(q_0^2(y))+2q_0^2(y)\rho'(q_0^2(y))}}\,,
\quad 0<y<f(l_0)
\end{gather}
and the compatible condition
\begin{align}
\label{phy-9}
q'_0(0)=0,\quad
q'_0(f(l_0))=\frac{f''(l_0)}{1+(f'(l_0))^2}q_0(f(l_0)),
\end{align}
where $c_*<c_1\le c_2<c^*$ are two constants.

\begin{definition}
A function $\varphi$ is called to be a solution to the problem \eqref{phy1}--\eqref{phy5},
if $\varphi\in C^1(\overline\Omega)$ with $|\nabla\varphi|\in C^{0,1}(\overline\Omega)$
and
$$
c_*<|\nabla\varphi(x,y)|\le c^*,\quad (x,y)\in\Omega,
$$
and $\varphi$ satisfies \eqref{phy1} in the distribution sense
and satisfies \eqref{phy2}--\eqref{phy5} pointwisely.
\end{definition}

\begin{remark}
Since the density of a supersonic flow is zero in the vacuum region,
the uniqueness of the solution to the problem \eqref{phy1}--\eqref{phy5}
means that the density and momentum are the same.
\end{remark}

\begin{remark}
It is proved in Theorem \ref{phys-thm81} that
for the smooth supersonic flow to the problem \eqref{phy1}--\eqref{phy5},
its speed belongs to $C^1(\Omega)$;
furthermore, its momentum also belongs to $C^1(\Omega)$
in the case $1<\gamma\le2$.
\end{remark}

\begin{remark}
For a supersonic flow to the problem \eqref{phy1}--\eqref{phy5}, \eqref{phy1} is hyperbolic and
the $x$-axis is regarded as the time axis. So, one can also consider the supersonic flow problem
in a finite nozzle.
\end{remark}

\begin{remark}
The inequality in \eqref{xc4} means that streamlines are rarefactive at the inlet,
which is a natural choice to exclude the sonic state in the flow,
as seen in Lemma \ref{nsc1}.
\end{remark}

\begin{remark}
As we will show in Remark \ref{nsc2} and Theorem \ref{ni2-thm-2},
conditions \eqref{phy-8} and \eqref{sss-0} are optimal in the following sense:
\eqref{phy-8} is sufficient and necessary
for the global existence of smooth supersonic flows in a straight nozzle;
while if the upper wall of the nozzle is a nonconvex perturbation of a straight line,
then there is a $q_0\in C^{1}([0,f(l_0)])$ satisfying \eqref{phy-8-q}--\eqref{phy-9}
such that a shock must form for the supersonic flow to the problem \eqref{phy1}--\eqref{phy5}.
\end{remark}

If there is no vacuum, the smooth supersonic flow problem \eqref{phy1}--\eqref{phy5}
can be treated more conveniently in the potential plane.
It is assumed that $(\Upsilon(0),0)$ is transformed into $(0,0)$
in the coordinates transformation
from the physical plane to the potential plane.
The mass flux and the stream function at the inlet can be given by
$$
m=\int_{0}^{f(l_0)}
q_0(y)\rho(q_0^2(y))\sqrt{1+(\Upsilon'(y))^2}dy
$$
and
\begin{gather*}
\psi(\Upsilon(y),y)=\Psi_{\text{\rm in}}(y)=\int_{0}^{y}
q_0(s)\rho(q_0^2(s))\sqrt{1+(\Upsilon'(s))^2}ds,\quad
0\le y\le f(l_{0}),
\end{gather*}
respectively.
If the flow speed at the upper wall is
$$
q(x,f(x))={\mathscr Q}_{\text{\rm up}}(x),\quad l_{0}\le x< l_1,
$$
then the potential function at the upper wall is
\begin{gather*}
\varphi(x,f(x))=\Phi_{\text{\rm up}}(x)=
\int_{l_0}^{x}
{{\mathscr Q}_{\text{\rm up}}(s)}\sqrt{1+(f'(s))^2}ds,\quad l_{0}\le x< l_1.
\end{gather*}
Therefore, in the potential plane,
the smooth supersonic flow problem \eqref{phy1}--\eqref{phy5}
can be formulated as follows
\begin{align}
\label{newsupp-eq}
&Q_{\varphi\varphi}-(b(Q)Q_{\psi})_\psi=0,
\quad&&(\varphi,\psi)\in(0,\zeta)\times(0,m),
\\
\label{newsupp-inbc1}
&Q(0,\psi)=Q_0(\psi),\quad&&\psi\in(0,m),
\\
\label{newsupp-inbc2}
&Q_\varphi(0,\psi)=G_0(\psi),\quad&&\psi\in(0,m),
\\
\label{newsupp-lbbc}
&Q_{\psi}(\varphi,0)=0,\quad&&\varphi\in(0,\zeta),
\\
\label{newsupp-ubbc}
&Q_{\psi}(\varphi,m)=-\frac{f''(X_{\text{\rm up}}(\varphi))(1+(f'(X_{\text{\rm up}}(\varphi)))^2)^{-3/2}}
{b(Q(\varphi,m)))A^{-1}_+(Q(\varphi,m))},
\quad&&\varphi\in(0,\zeta),
\\
\label{newsupp-q}
&{\mathscr Q}_{\text{\rm up}}(x)=A^{-1}_+(Q(\Phi_{\text{\rm up}}(x),m)),
\quad&&x\in(l_0,l),
\end{align}
where $\zeta=\Phi_{\text{\rm up}}(l)\,(l_0<l\le l_1)$,
$X_{\text{\rm up}}$ is the inverse function of $\Phi_{\text{\rm up}}$
and a solution $Q$ is a negative function.
Moreover,
\begin{align}
\label{xc5}
Q_0(\psi)=A(q_0(y))\Big|_{y=Y_{\text{\rm in}}(\psi)},\quad
G_0(\psi)=\frac{\Upsilon''(y)(1+(\Upsilon'(y))^2)^{-3/2}}{q_0(y)\rho(q_0^2(y))}
\Big|_{y=Y_{\text{\rm in}}(\psi)},
\quad 0<\psi<m
\end{align}
with $Y_{\text{\rm in}}$ being the inverse function of
$\Psi_{\text{\rm in}}$.
It follows from \eqref{phy-8-q}--\eqref{phy-9} and \eqref{xc5} that
$Q_0\in C^{1}([0,m])$ and $G_0\in C([0,m])$
satisfy
\begin{align}
\label{sss-1}
-\mu_2\le Q_0(\psi)\le-\mu_1,\quad
0\le-G_0(\psi)\pm b^{1/2}(Q_0(\psi))Q'_0(\psi)\le\mu_3,\quad\psi\in(0,m)
\end{align}
with three constants $0<\mu_1\le\mu_2$ and $\mu_3>0$,
and
\begin{align}
\label{sss-2}
Q'_0(0)=0,\quad
Q'_0(m)=-\frac{f''(l_0)(1+(f'(l_0)^2)^{-3/2}}
{b(Q_0(m))A^{-1}_+(Q_0(m))}.
\end{align}
It can be verified easily that the problem \eqref{newsupp-eq}--\eqref{newsupp-q} is equivalent to
\begin{align}
\label{newsupp-neq1}
&W_\varphi+b^{1/2}(Q)W_\psi
=\frac14b^{-1}(Q)p(Q)W(W+Z),
\quad&&(\varphi,\psi)\in(0,\zeta)\times(0,m),
\\
\label{newsupp-neq2}
&Z_\varphi-b^{1/2}(Q)Z_\psi
=-\frac14b^{-1}(Q)p(Q)Z(W+Z),
\quad&&(\varphi,\psi)\in(0,\zeta)\times(0,m),
\\
\label{newsupp-ninbc1}
&W(0,\psi)=W_0(\psi)=G_0(\psi)-b^{1/2}(Q_0(\psi))Q'_0(\psi),\quad&&\psi\in(0,m),
\\
\label{newsupp-ninbc2}
&Z(0,\psi)=Z_0(\psi)=-G_0(\psi)-b^{1/2}(Q_0(\psi))Q'_0(\psi),\quad&&\psi\in(0,m),
\\
\label{newsupp-nlbbc}
&W(\varphi,0)+Z(\varphi,0)=0,\quad&&\varphi\in(0,\zeta),
\\
\label{newsupp-nubbc}
&W(\varphi,m)+Z(\varphi,m)=\frac{2f''(X_{\text{\rm up}}(\varphi))(1+(f'(X_{\text{\rm up}}(\varphi)))^2)^{-3/2}}
{b^{1/2}(Q(\varphi,m))A^{-1}_+(Q(\varphi,m))},
\quad&&\varphi\in(0,\zeta),
\\
\label{newsupp-q1}
&Q_\varphi=\frac{1}{2}(W-Z),
\quad&&(\varphi,\psi)\in(0,\zeta)\times(0,m),
\\
\label{newsupp-q3}
&Q(0,\psi)=Q_0(\psi),\quad&&\psi\in(0,m),
\\
\label{newsupp-q4}
&{\mathscr Q}_{\text{\rm up}}(x)=A^{-1}_+(Q(\Phi_{\text{\rm up}}(x),m)),
\quad&&x\in(l_0,l).
\end{align}
It follows from \eqref{sss-1} that $Q_0\in C^{1}([0,m])$ and $W_0,Z_0\in C([0,m])$
satisfy
\begin{gather}
\label{sss-3-0}
-\mu_2\le Q_0(\psi)\le-\mu_1,\quad\psi\in(0,m),
\\
\label{sss-3}
-\mu_3\le W_0(\psi)\le0,\quad
0\le Z_0(\psi)\le\mu_3,\quad\psi\in(0,m).
\end{gather}

\subsection{Necessary condition at the inlet for a global smooth supersonic flow in a straight nozzle}

We first study the formation of the singularity of supersonic flows to the Cauchy problem
\begin{align}
\label{cv1-a1}
&Q_{\varphi\varphi}-(b(Q)Q_{\psi})_\psi=0,\quad &&(\varphi,\psi)\in(0,+\infty)\times\mathbb R,
\\
\label{cv1-a2}
&Q(0,\psi)=Q_0(\psi),\quad&&\psi\in\mathbb R,
\\
\label{cv1-a3}
&Q_\varphi(0,\psi)=G_0(\psi),\quad&&\psi\in\mathbb R.
\end{align}

Recall the following lemma for invariant quantities on characteristics.

\begin{lemma}
\label{lemma-sc4}
(\cite{WX3})
Let $\Sigma_{+}\,(\Sigma_{-})$ be a positive (negative) characteristic
for a $C^1$ supersonic flow without vacuum.
Then, ${\mathscr H}(Q)\mp \theta$ is invariant on $\Sigma_{\pm}$, where
$$
{\mathscr H}(s)=\int_{s}^{0}b^{1/2}(s)dt,\quad s\le 0.
$$
\end{lemma}

\begin{remark}
Note that ${\mathscr H}(-\infty)=\lim_{s\to-\infty}{\mathscr H}(s)$ exists and is bounded.
\end{remark}

The following three lemmas will be needed.

\vskip8mm

\hskip45mm
\setlength{\unitlength}{0.6mm}
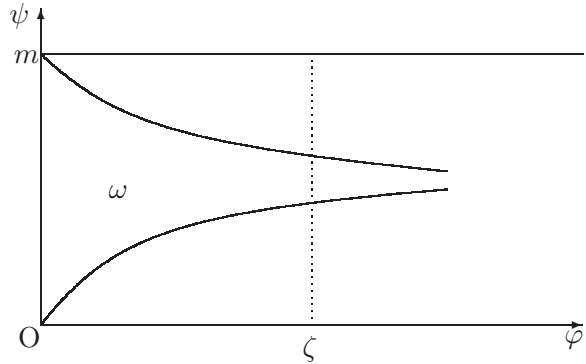
\begin{picture}(250,65)
\put(0,0){\vector(1,0){120}}
\put(0,0){\vector(0,1){70}}
\put(116,-4){$\varphi$} \put(-7,67){$\psi$}
\put(0,60){\qbezier(0,0)(55,0)(120,0)}

\put(0,0){\cbezier(0,0)(16,20)(30,25)(90,30)}

\put(0,60){\cbezier(0,0)(16,-15)(30,-20)(90,-26)}

\put(15,28){$\omega$}

\put(60,0){\qbezier[40](0,0)(0,30)(0,60)}
\put(58,-7){$\zeta$}



\put(-6,58){$m$}

\put(-5,-5){O}

\end{picture}

\vskip10mm

\begin{lemma}
\label{aaasss-thm2}
For given $\zeta>0$, $Q_0\in C^{1}([0,m])$ with $\sup_{(0,m)}Q_0<0$ and $G_0\in C([0,m])$,
assume that $Q\in C^{1}(\overline\omega)$ solves
\begin{align*}
&Q_{\varphi\varphi}-(b(Q)Q_{\psi})_\psi=0,
\quad&&(\varphi,\psi)\in\omega,
\\
&Q(\zeta,\psi)=Q_0(\psi),\quad&&\psi\in(0,m),
\\
&Q_\varphi(\zeta,\psi)=G_0(\psi),\quad&&\psi\in(0,m),
\end{align*}
where $\omega$ is the domain bounded by the line $\varphi=0$, the line $\varphi=\zeta$,
the positive characteristic from $(0,0)$
and the negative characteristic from $(0,m)$.
Then
\begin{align*}
-Q_\varphi(\varphi,\psi)\pm b^{1/2}(Q(\varphi,\psi))Q_\psi(\varphi,\psi)
\le\max\Big\{0,\sup_{(0,m)}\big(-G_0+b^{1/2}(Q_0)|Q'_0|\big)\Big\},
\quad(\varphi,\psi)\in\omega.
\end{align*}
If $\sup_{(0,m)}G_0\le0$ additionally, then
$\sup_\omega Q=\sup_{(0,m)}Q_0$.
\end{lemma}

\Proof
For $(\varphi,\psi)\in\overline\omega$, set
\begin{align*}
W(\varphi,\psi)=Q_\varphi(\varphi,\psi)-b^{1/2}(Q(\varphi,\psi))Q_\psi(\varphi,\psi),
\quad
Z(\varphi,\psi)=-Q_\varphi(\varphi,\psi)-b^{1/2}(Q(\varphi,\psi))Q_\psi(\varphi,\psi).
\end{align*}
Then, $(W,Z)\in C(\overline\omega)
\times C(\overline\omega)$ solving
\begin{align*}
&W_\varphi+b^{1/2}(Q)W_\psi
=\frac14b^{-1}(Q)p(Q)W(W+Z),
&&(\varphi,\psi)\in\omega,
\\
&Z_\varphi-b^{1/2}(Q)Z_\psi
=-\frac14b^{-1}(Q)p(Q)Z(W+Z),
&&(\varphi,\psi)\in\omega,
\\
&W(0,\psi)=W_0(\psi)=G_0(\psi)-b^{1/2}(Q_0(\psi))Q'_0(\psi),
&&\psi\in(0,m),
\\
&Z(0,\psi)=Z_0(\psi)=-G_0(\psi)-b^{1/2}(Q_0(\psi))Q'_0(\psi),
&&\psi\in(0,m).
\end{align*}
Denote
\begin{align}
\label{acv4}
M=\max\Big\{0,-\inf_{(0,m)}W_0,\sup_{(0,m)}Z_0\Big\}
=\max\Big\{0,\sup_{(0,m)}\big(-G_0+b^{1/2}(Q_0)|Q'_0|\big)\Big\}.
\end{align}
It suffices to prove that
\begin{align}
\label{acv5}
W(\varphi,\psi)\ge-M,\quad Z(\varphi,\psi)\le M,
\quad(\varphi,\psi)\in\omega.
\end{align}
If \eqref{acv5} is invalid, then one of the following two cases occurs:
\\

(i) There exist two constants $\varphi_+>0$ and $M_+>M$,
as well as a positive characteristic
$$
\Sigma_+:\Psi'_+(\varphi)=b^{1/2}(Q(\varphi,\Psi_+(\varphi))),
\quad0<\varphi<\varphi_+
$$
such that $W(\varphi_+,\Psi_+(\varphi_+))=-M_+$ and
$Z(\cdot,\Psi_+(\cdot))\le M_+$ on $[0,\varphi_+]$.

(ii) There exist two constants $\varphi_->0$ and $M_->M$,
as well as a negative characteristic
$$
\Sigma_-:\Psi'_-(\varphi)=-b^{1/2}(Q(\varphi,\Psi_-(\varphi))),
\quad0<\varphi<\varphi_-
$$
such that $Z(\varphi_-,\Psi_-(\varphi_-))=M_-$ and
$W(\cdot,\Psi_-(\cdot))\ge-M_-$ on $[0,\varphi_-]$.

Without loss of generality, it is assumed that (i) occurs.
On $\Sigma_+$, $W$ satisfies
\begin{align}
\label{acv6}
\frac{d}{d\varphi}W(\varphi,\Psi_{+}(\varphi))
=R(\varphi)W(\varphi,\Psi_{+}(\varphi))
(W(\varphi,\Psi_{+}(\varphi))+Z(\varphi,\Psi_{+}(\varphi))),
\quad0<\varphi<\varphi_+,
\end{align}
where
\begin{align*}
R(\varphi)=\frac14b^{-1}(Q(\varphi,\Psi_{+}(\varphi)))
p(Q(\varphi,\Psi_{+}(\varphi))),\quad
0\le\varphi\le\varphi_+.
\end{align*}
Since $W(\varphi_+,\Psi_+(\varphi_+))=-M_+<0$, \eqref{acv6} shows
$W(\cdot,\Psi_+(\cdot))<0$ on $[0,\varphi_+]$.
Rewrite \eqref{acv6} as
\begin{align*}
&\Big(W(\varphi,\Psi_{+}(\varphi))
\mbox{exp}\Big\{-\int_{0}^{\varphi}R(s)W(s,\Psi_+(s))ds\Big\}\Big)'
\\
=&-Z(\varphi,\Psi_{+}(\varphi))\Big(\mbox{exp}\Big\{-\int_{0}^{\varphi}
R(s)W(s,\Psi_+(s))ds\Big\}\Big)',
\quad0<\varphi<\varphi_+,
\end{align*}
which, together with $Z(\cdot,\Psi_+(\cdot))\le M_+$
and $W(\cdot,\Psi_+(\cdot))<0$ on $[0,\varphi_+]$, leads to
\begin{align*}
&\Big(W(\varphi,\Psi_{+}(\varphi))
\mbox{exp}\Big\{-\int_{0}^{\varphi}R(s)W(s,\Psi_+(s))ds\Big\}\Big)'
\\
\ge&-M_+\Big(\mbox{exp}\Big\{-\int_{0}^{\varphi}
R(s)W(s,\Psi_+(s))ds\Big\}\Big)',
\quad0<\varphi<\varphi_+.
\end{align*}
Therefore,
\begin{align*}
W_0(\Psi_{+}(0))+M_+\le(W(\varphi_+,\Psi_{+}(\varphi_+))+M_+)
\mbox{exp}\Big\{-\int_{0}^{\varphi_+}R(s)W(s,\Psi_+(s))ds\Big\}=0,
\end{align*}
which contradicts \eqref{acv4} and $M_+>M$.
Hence \eqref{acv5} holds.

Below it is assumed that $\sup_{(0,m)}G_0\le0$ additionally, which
implies $\theta_\psi(0,\cdot)\ge 0$ in $(0,m)$.
Fix $(\tilde\varphi,\tilde\psi)\in\omega$.
Assume that
$$
\tilde\Sigma_{\pm}:
\tilde\Psi'_{\pm}(\varphi)=\pm b^{1/2}(Q(\varphi,\tilde\Psi_{\pm}(\varphi))),
\quad
0<\varphi<\tilde\varphi,
\quad\tilde\Psi_{\pm}(\tilde\varphi)=\tilde\psi
$$
are the positive and negative characteristics approaching $(\tilde\varphi,\tilde\psi)$.
It follows from Lemma \ref{lemma-sc4} that
\begin{align*}
{\mathscr H}(Q_0(\tilde\Psi_{\pm}(0)))\mp \theta(0,\tilde\Psi_{\pm}(0))
={\mathscr H}(Q(\tilde\varphi,\tilde\psi))\mp \theta(\tilde\varphi,\tilde\psi).
\end{align*}
Therefore,
\begin{align*}
2{\mathscr H}(Q(\tilde\varphi,\tilde\psi))
={\mathscr H}(Q_0(\tilde\Psi_{+}(0)))+{\mathscr H}(Q_0(\tilde\Psi_{-}(0)))
-\theta(0,\tilde\Psi_{+}(0))+\theta(0,\tilde\Psi_{-}(0)),
\end{align*}
which, together with $\theta_\psi(0,\cdot)\ge 0$ in $(0,m)$, leads to
\begin{align*}
2{\mathscr H}(Q(\tilde\varphi,\tilde\psi))
\ge{\mathscr H}(Q_0(\Psi_{+}(0)))+{\mathscr H}(Q_0(\Psi_{-}(0)))
\ge2{\mathscr H}\Big(\sup_{(0,m)}Q_0\Big).
\end{align*}
Then, $\sup_\omega Q=\sup_{(0,m)}Q_0$ follows from
the arbitrariness of $(\tilde\varphi,\tilde\psi)\in\omega$.
$\hfill\Box$\vskip 4mm

\begin{lemma}
\label{nsc1}
Assume that $Q\in C^1([0,\zeta]\times[0,m])$ is a solution to
the problem \eqref{newsupp-eq}--\eqref{newsupp-q}.
If $\inf_{(l_0,l_1)}f'=f'(l_0)$
and $\sup_{(0,f(l_0))}\Upsilon''\le0$,
then $\inf_{(0,\zeta)\times(0,m)}A_+^{-1}(Q)=\inf_{(0,f(l_0))}q_0$.
\end{lemma}

\Proof
Note that
$$
\theta(\varphi,0)=0=\inf_{(0,m)}\theta(0,\cdot),\quad
\theta(\varphi,m)\ge\theta(0,m)=\sup_{(0,m)}\theta(0,\cdot),\quad
0\le\varphi\le m.
$$
It follows from the method of characteristics and Lemma \ref{lemma-sc4} that
\begin{align*}
Q(\varphi,0)\le\sup_{(0,m)}Q_0,\quad
Q(\varphi,m)\le\sup_{(0,m)}Q_0,\quad
0\le\varphi\le m.
\end{align*}
Then, a similar argument as the end of the proof of Lemma \ref{aaasss-thm2} leads
$\sup_{(0,\zeta)\times(0,m)}Q=\sup_{(0,m)}Q_0$.
$\hfill\Box$\vskip 4mm

\begin{lemma}
\label{ni}
Assume that $Q\in C^{1}((0,+\infty)\times\mathbb R)$ solves \eqref{cv1-a1}
and satisfies
\begin{align}
\label{cv1-1}
-M(\varphi+1)^2\le Q(\varphi,\psi)<0,\quad (\varphi,\psi)\in(0,+\infty)\times\mathbb R
\end{align}
with $M>0$ being a constant.
Then,
\begin{align*}
-Q_{\varphi}(\varphi,\psi)\pm b^{1/2}(Q(\varphi,\psi))Q_{\psi}(\varphi,\psi)\ge0,
\quad (\varphi,\psi)\in(0,+\infty)\times\mathbb R.
\end{align*}
\end{lemma}

\Proof
For $(\varphi,\psi)\in(0,+\infty)\times\mathbb R$, set
\begin{align*}
W(\varphi,\psi)=Q_\varphi(\varphi,\psi)-b^{1/2}(Q(\varphi,\psi))Q_\psi(\varphi,\psi),
\quad
Z(\varphi,\psi)=-Q_\varphi(\varphi,\psi)-b^{1/2}(Q(\varphi,\psi))Q_\psi(\varphi,\psi).
\end{align*}
Then, $(W,Z)\in C((0,+\infty)\times\mathbb R)
\times C((0,+\infty)\times\mathbb R)$ solves
\begin{align}
\label{cv2-2}
&W_\varphi+b^{1/2}(Q)W_\psi
=\frac14b^{-1}(Q)p(Q)W(W+Z),
&&\quad (\varphi,\psi)\in(0,+\infty)\times\mathbb R,
\\
\label{cv3}
&Z_\varphi-b^{1/2}(Q)Z_\psi
=-\frac14b^{-1}(Q)p(Q)Z(W+Z),
&&\quad (\varphi,\psi)\in(0,+\infty)\times\mathbb R.
\end{align}
It follows from \eqref{cv1-1} that
each characteristic is unbounded.
Set
$$
h(s)=\mbox{exp}\Big\{\frac14\int_{-1}^sb^{-1}(t)p(t)dt\Big\},\quad s<0.
$$
It follows from \eqref{a8-8} and $\gamma>1$ that
\begin{align}
\label{a8-8-0}
\lim_{s\to 0^-}(-s)^{1/4}h(s)=0,\quad
\lim_{s\to-\infty}(-s)^{1/2}h(s)=0.
\end{align}
By \eqref{cv1-1} and \eqref{a8-8-0},
there exists a positive constant $M_0$ depending only on $\gamma$ and $M$
such that
\begin{align}
\label{cv6}
b^{-1}(Q(\varphi,\psi))p(Q(\varphi,\psi))\ge 4M_0(\varphi+1)^{-1}h(Q(\varphi,\psi)),
\quad (\varphi,\psi)\in(0,+\infty)\times\mathbb R
\end{align}
and
\begin{align}
\label{cv6-z}
b^{-3/2}(Q(\varphi,\psi))p(Q(\varphi,\psi))\ge 4M_0h(Q(\varphi,\psi)),
\quad (\varphi,\psi)\in(0,+\infty)\times\mathbb R.
\end{align}

There are four cases to be discussed.

(i)
Assume that
$$
\Psi'_+(\varphi)=b^{1/2}(Q(\varphi,\Psi_+(\varphi))),
\quad\varphi>\varphi_+>0
$$
is a positive characteristic.
Then, it follows from \eqref{cv2-2} that
\begin{align*}
(W(\varphi,\Psi_+(\varphi)))'=&
\frac14b^{-1}(Q(\varphi,\Psi_+(\varphi)))p(Q(\varphi,\Psi_+(\varphi)))
W(\varphi,\Psi_+(\varphi))
\\
&\qquad\cdot
(W(\varphi,\Psi_+(\varphi))-(Q(\varphi,\Psi_+(\varphi)))'),\quad
\varphi>\varphi_+
\end{align*}
and thus
\begin{align}
\label{cv4}
&(h(Q(\varphi,\Psi_+(\varphi)))W(\varphi,\Psi_+(\varphi)))'
\nonumber
\\
=&
\frac14b^{-1}(Q(\varphi,\Psi_+(\varphi)))p(Q(\varphi,\Psi_+(\varphi)))
h(Q(\varphi,\Psi_+(\varphi)))W^2(\varphi,\Psi_+(\varphi)),\quad
\varphi>\varphi_+.
\end{align}
Substituting \eqref{cv6} into \eqref{cv4} yields
\begin{align*}
(h(Q(\varphi,\Psi_+(\varphi)))W(\varphi,\Psi_+(\varphi)))'\ge
M_0(\varphi+1)^{-1}(h(Q(\varphi,\Psi_+(\varphi)))W(\varphi,\Psi_+(\varphi)))^2,\quad
\varphi>\varphi_+,
\end{align*}
which leads to
$W(\cdot,\Psi_+(\cdot))\le0$ on $[\varphi_+,+\infty)$.

(ii)
Assume that
$$
\Psi'_-(\varphi)=-b^{1/2}(Q(\varphi,\Psi_-(\varphi))),
\quad\varphi>\varphi_->0
$$
is a negative characteristic.
Then, it follows from \eqref{cv3} and \eqref{cv6} that
\begin{align*}
&(h(Q(\varphi,\Psi_-(\varphi)))Z(\varphi,\Psi_-(\varphi)))'
\\
=&-\frac14b^{-1}(Q(\varphi,\Psi_-(\varphi)))p(Q(\varphi,\Psi_-(\varphi)))
h(Q(\varphi,\Psi_-(\varphi)))Z^2(\varphi,\Psi_-(\varphi))
\\
\le&
-M_0(\varphi+1)^{-1}(h(Q(\varphi,\Psi_-(\varphi)))Z(\varphi,\Psi_-(\varphi)))^2,\quad
\varphi>\varphi_-,
\end{align*}
which leads to
$Z(\cdot,\Psi_-(\cdot))\ge0$ on $[\varphi_-,+\infty)$.

(iii)
Assume that
$$
\Phi'_+(\psi)=b^{-1/2}(Q(\Phi_+(\psi),\psi)),
\quad\psi>\psi_+
$$
is a positive characteristic.
Then, it follows from \eqref{cv2-2} and \eqref{cv6-z} that
\begin{align*}
&(h(Q(\Phi_+(\psi),\psi))W(\Phi_+(\psi),\psi))'
\\
=&\frac14b^{-3/2}(Q(\Phi_+(\psi),\psi))p(Q(\Phi_+(\psi),\psi))
h(Q(\Phi_+(\psi),\psi))W^2(\Phi_+(\psi),\psi)
\\
\ge&
M_0(h(Q(\Phi_+(\psi),\psi))W(\Phi_+(\psi),\psi))^2,\quad
\psi>\psi_+,
\end{align*}
which leads to
$W(\Phi_+(\cdot),\cdot)\le0$ on $[\psi_+,+\infty)$.

(iv)
Assume that
$$
\Phi'_-(\psi)=-b^{-1/2}(Q(\Phi_-(\psi),\psi)),
\quad\psi<\psi_-
$$
is a negative characteristic.
Then, it follows from \eqref{cv3}  and \eqref{cv6-z} that
\begin{align*}
&(h(Q(\Phi_-(\psi),\psi))Z(\Phi_-(\psi),\psi))'
\\
=&-\frac14b^{-3/2}(Q(\Phi_-(\psi),\psi))p(Q(\Phi_-(\psi),\psi))
h(Q(\Phi_-(\psi),\psi))Z^2(\Phi_-(\psi),\psi)
\\
\le&
-M_0(h(Q(\Phi_-(\psi),\psi))Z(\Phi_-(\psi),\psi))^2,\quad
\psi<\psi_-,
\end{align*}
which leads to
$Z(\Phi_-(\cdot),\cdot)\le0$ on $(-\infty,\psi_-]$.

Summing up, one gets that $W\le0$ and $Z\ge0$ in $(0,+\infty)\times\mathbb R$.
$\hfill\Box$\vskip 4mm

\begin{remark}
\label{ni-r}
Lemma \ref{ni} remains true if
\eqref{cv1-1} is relaxed by
\begin{align*}
0<-Q(\varphi,\psi)<
\left\{
\begin{aligned}
&\lambda_1(\varphi+1)^{4/(3-\gamma)-\lambda_2},
&&\quad\mbox{if }1<\gamma<3,
\\
&\lambda_3(\varphi+1)^{\lambda_4},
&&\quad\mbox{if }\gamma=3,
\\
&+\infty,
&&\quad\mbox{if }\gamma>3,
\end{aligned}
\right.
\quad(\varphi,\psi)\in(0,+\infty)\times\mathbb R,
\end{align*}
where $\lambda_i\,(i=1,2,3,4)$ are positive constants with
$\lambda_2<4/(3-\gamma)$.
\end{remark}

It follows from Lemmas \ref{aaasss-thm2} and \ref{ni} that

\begin{proposition}
\label{ni2}
Assume that
$Q_0\in C^{1}(\mathbb R)\cap L^\infty(\mathbb R)$ with $\sup_{\mathbb R}Q_0<0$ and
$Q'_0\in L^\infty(\mathbb R)$,
and $G_0\in C(\mathbb R)\cap L^\infty(\mathbb R)$.

{\rm(i)} There is not any vacuum point for a $C^1$ supersonic flow to the problem \eqref{cv1-a1}--\eqref{cv1-a3}.
More precisely,
if $Q\in C^{1}(\mathbb R\times[0,\zeta])$ is a local solution to the problem \eqref{cv1-a1}--\eqref{cv1-a3}
for some positive constant $\zeta$, then
\begin{align*}
-Q_\varphi(\varphi,\psi)
\le\max\Big\{0,\sup_{\mathbb R}\big(-G_0+b^{1/2}(Q_0)|Q'_0|\big)\Big\},
\quad(\varphi,\psi)\in\mathbb R\times[0,\zeta].
\end{align*}

{\rm(ii)} If $\sup_{\mathbb R}G_0\le0$,
then $C^1$ supersonic flows to the problem \eqref{cv1-a1}--\eqref{cv1-a3}
are uniformly away from the sonic state.
More precisely,
if $Q\in C^{1}(\mathbb R\times[0,\zeta])$ is a local solution to the problem \eqref{cv1-a1}--\eqref{cv1-a3}
for some positive constant $\zeta$, then
$\sup_{\mathbb R\times(0,\zeta)}Q=\sup_{\mathbb R}Q_0$.

{\rm(iii)} If $\sup_{\mathbb R}(G_0+b^{1/2}(Q_0)|Q'_0|)>0$,
then there is not a global $C^1$ solution to the problem \eqref{cv1-a1}--\eqref{cv1-a3}.
In particular, if $\sup_{\mathbb R}G_0\le0$ additionally,
then there must be a shock for the global supersonic flow to the problem \eqref{cv1-a1}--\eqref{cv1-a3}
and the flow is uniformly away from the sonic and vacuum state before the shock forms.
\end{proposition}

\begin{remark}
$\sup_{\mathbb R}G_0\le0$ is equivalent to that $\theta(0,\cdot)$ is nondecreasing in $\mathbb R$.
\end{remark}

\begin{remark}
Proposition \ref{ni2} shows that
$\sup_{\mathbb R}(G_0+b^{1/2}(Q_0)|Q'_0|)\le0$
is a necessary condition for the existence of
global $C^1$ solutions to the problem \eqref{cv1-a1}--\eqref{cv1-a3}.
As will be shown in Theorem \ref{ni3}, it is also a sufficient condition.
\end{remark}

Applying Proposition \ref{ni2} to the supersonic flow problem in a straight nozzle gives the following theorem.

\begin{theorem}
\label{ni2-thm}
Assume that
$f\in C^2([l_0,+\infty))$ satisfies
\begin{align}
\label{aaasss-0}
f(l_0)>0,\quad
f'(l_0)\ge0,\quad
f''(x)=0\mbox{ for each } x\ge l_0,
\end{align}
$\Upsilon\in C^{2}([0,f(l_0)])$ satisfies \eqref{xc4}
and $q_0\in C^{1}([0,f(l_0)])$ satisfy \eqref{phy-8-q} and \eqref{phy-9}.

{\rm(i)} A smooth supersonic flow to the problem \eqref{phy1}--\eqref{phy5}
never approaches the sonic and vacuum state in any bounded region.
In particular, the speed of a smooth supersonic flow is always not less than $\inf_{(0,f(l_0))}q_0$.

{\rm(ii)} If \eqref{phy-8} is invalid,
then there is not a global smooth supersonic flow to the problem \eqref{phy1}--\eqref{phy5}.
More precisely, there must be a shock for the global supersonic flow to the problem \eqref{phy1}--\eqref{phy5}
and the flow is uniformly away from the sonic and vacuum state before the shock forms.
\end{theorem}

\begin{remark}
\label{nsc2}
Under the assumptions of Theorem \ref{ni2-thm}, \eqref{phy-8}
is necessary for the existence of global smooth supersonic flows.
As will be shown in Theorem \ref{phys-thm51}, it is also a sufficient condition.
Furthermore, there is not any sonic and vacuum point for such global supersonic flows.
\end{remark}

\subsection{On the global smooth transonic flow problem}

As a direct application of the results on the smooth supersonic flow problem \eqref{phy1}--\eqref{phy5},
one can get the global extension of the smooth transonic flow obtained
in Theorems \ref{theoremm5} and \ref{theoremm6}.
Assume that $f$ is a function defined on $[l_-,l_1)\,(l_-<0<l_1\le+\infty)$
satisfying $\lim_{x\to l_1^-}(x+f(x))=+\infty$.
Let $\Omega$ be the unbounded domain bounded
by the upper wall $\Gamma_{\text{\rm up}}:y=f(x)\,(l_-\le x< l_1)$,
the lower wall $x$-axis
and the inlet $\Gamma_{\text{\rm in}}:x=g(y)\,(0\le y\le f(l_-),g(f(l_-))=l_-)$.
Then, the global smooth transonic problem is formulated as follows
\begin{align}
\label{sonic1}
&\mbox{div}(\rho(|\nabla\varphi|^2)\nabla\varphi)=0,\quad&&(x,y)\in\Omega,
\\
\label{sonic2}
&\varphi(g(y),y)=C_{\text{\rm in}},\quad&&0<y<f(l_-),
\\
\label{sonic3}
&\rho(|\nabla\varphi(x,0)|^2)\pd{\varphi}y(x,0)=0,\quad&&g(0)<x<l_1,
\\
\label{sonic4}
&\rho(|\nabla\varphi(x,f(x))|^2)
\Big(\pd{\varphi}y(x,f(x))-f'(x)\pd{\varphi}x(x,f(x))\Big)=0,\quad&&l_-<x<l_1,
\\
\label{sonic5}
&\varphi(0,y)=0,\,|\nabla\varphi(0,y)|=c_*,\quad&&0<y<f(0),
\\
\label{sonic6}
&|\nabla\varphi(x,y)|<c_*,\quad&&(x,y)\in\Omega\mbox{ and }x<0,
\\
\label{sonic7}
&|\nabla\varphi(x,y)|>c_*,\quad&&(x,y)\in\Omega\mbox{ and }x>0,
\end{align}
where $C_{\text{\rm in}}$ is a free constant.

\vskip5mm

\hskip55mm
\setlength{\unitlength}{0.6mm}
\begin{picture}(250,120)
\put(-60,0){\vector(1,0){190}}
\put(0,52){\vector(0,1){65}}
\put(127,-4){$x$} \put(-4,113){$y$}

\put(60,48){\qbezier(-40,5)(-45,4.15)(-60,4)}
\put(60,48){\cbezier(-40,5)(-20,6.5)(5,10)(15,16)}
\put(60,48){\cbezier(15,16)(25,21)(40,28)(63,53)}

\put(35,59){$\Gamma_{\text{\rm up}}$}
\put(-20,25){smooth transonic flow}
\put(20,40){$\rho>0$}

\put(0,48){\cbezier(-40,10)(-25,5.3)(-15,4.3)(0,4)}

\put(0,0){\cbezier(-40,58)(-45,48)(-48,23)(-48,0)}
\put(-58,20){$\Gamma_{\text{\rm in}}$}

\put(105,80){$\rho=0$}
\put(72,62.5){\circle*{1.5}}
\put(75,55){\qbezier(-3,7)(26,17)(49,27)}

\end{picture}
\vskip5mm

For the readers' convenience,
we state the main results for the problem \eqref{sonic1}--\eqref{sonic7}.
The first one is the following existence and uniqueness.

\begin{theorem}
\label{phys-thmm1}
Assume that $f\in C^{4}([l_-,l_+])\cap C^2([l_-,l_1))\,
(-\delta_{0,-}\le l_-<0<l_+<l_1\le+\infty)$
and $g\in C^{3,\alpha}([0,f(l_-)])\,(0<\alpha<1)$ satisfy the assumptions of Theorem \ref{theoremm5}
and
\begin{align*}
\lim_{x\to l_1^-}(x+f(x))=+\infty,\quad
f''(x)\ge0\mbox{ for each } l_+\le x<l_1.
\end{align*}
Then the problem \eqref{sonic1}--\eqref{sonic7} admits uniquely a global smooth solution
$\varphi\in C^{1}(\overline\Omega)$ with $|\nabla\varphi|\in C^{0,1}(\overline\Omega)$,
which possesses the properties of Theorem \ref{theoremm5}
near the sonic state. Moreover,
the global smooth transonic flow belongs to one and only one of the following two cases:
\\

{\rm\bf Case I} Global smooth transonic flow without vacuum.
In this case, $\varphi\in C^{2}(\overline\Omega)$ and
$|\nabla\varphi|<c^*$ on $\overline\Omega$.
If $f''(1+(f')^2)^{-3/2}\in L^\infty(l_+,l_1)$
additionally, then $|\nabla\varphi|$ is globally Lipschitz continuous in $\Omega$.

{\rm\bf Case II} Global smooth transonic flow with vacuum.
In this case, set
$$
x_0=\sup\big\{l\in(0,l_1):|\nabla\varphi(x,f(x))|<c^*\mbox{ for each }0<x<l\big\}
$$
and
$$
\Omega_{\text{\rm v}}=\big\{(x,y)\in\Omega:x>x_0,y>f(x_0)+f'(x_0)(x-x_0)\big\}.
$$
Then

{\rm (i)} $\overline\Omega_{\text{\rm v}}$ is the set of vacuum points.
Moreover, $\varphi\in C^{2}(\overline\Omega\setminus\overline\Omega_{\text{\rm v}})$
and
\begin{align*}
\pd\varphi x(x,y)=\frac{c^*}{\sqrt{1+(f'(x_0))^2}},\quad
\pd\varphi y(x,y)=\frac{f'(x_0)c^*}{\sqrt{1+(f'(x_0))^2}},\quad
(x,y)\in\partial\Omega_{\text{\rm v}}\cap\Omega.
\end{align*}

{\rm (ii)}
$|\nabla\varphi(x,f(x))|=c^*-O\big((x_0-x)^2\big)$ as $x\to x_0^-$.

{\rm (iii)} $|\nabla\varphi|$ is globally Lipschitz continuous in $\Omega$.

{\rm (iv)} $|\nabla\varphi|\in C^1(\overline\Omega\setminus{(x_0,f(x_0))})$
and
for each $(\tilde x,\tilde y)\in\partial\Omega_{\text{\rm v}}\cap\Omega$,
\begin{align*}
{\lim_{\stackrel{(x,y)\to(\tilde x,\tilde y)}
{(x,y)\in\Omega\setminus\overline\Omega_{\text{\rm v}}}}}
\nabla|\nabla\varphi(x,y)|=(0,0),
\quad
{\lim_{\stackrel{(x,y)\to(\tilde x,\tilde y)}
{(x,y)\in\Omega\setminus\overline\Omega_{\text{\rm v}}}}}
(c^*-|\nabla\varphi(x,y)|)^{1/2}\nabla\arctan\frac{\varphi_y(x,y)}{\varphi_x(x,y)}=(0,0).
\end{align*}
In particular, if $1<\gamma\le2$, then
$\rho(|\nabla\varphi|^2)\nabla\varphi\in C^1(\overline\Omega\setminus{(x_0,f(x_0))})$.

{\rm (v)} If $1<\gamma<3$ and $f''(x_0)>0$ additionally,
then for any $\varepsilon>0$,
\begin{align*}
{\lim_{\stackrel{(x,y)\to(\tilde x,\tilde y)}
{(x,y)\in\Omega\setminus\overline\Omega_{\text{\rm v}}}}}
(c^*-|\nabla\varphi(x,y)|)^{-(\gamma+1+\varepsilon)/(4\gamma-4)}
\nabla |\nabla\varphi(x,y)|\cdot(-1,f'(x_0))=+\infty,
\quad (\tilde x,\tilde y)\in\partial\Omega_{\text{\rm v}}\cap\Omega.
\end{align*}
\end{theorem}

\begin{remark}
If $f\in C^{2,1}([l_+,l_1))$ additionally, then
the acceleration of the flow is Lipschitz continuous in the gas region.
\end{remark}

The following theorem gives some sufficient conditions to determine
whether a vacuum forms or not for global smooth transonic flows.

\begin{theorem}
\label{phys-thmm2}
{\rm (i)} Set
\begin{align*}
\check x=\inf\Big\{x\in(0,l_1):\arctan f'(x)\ge
\int_{-\infty}^{h(x)}b^{1/2}(s)ds\Big\},
\end{align*}
where
$$
h(x)=-\frac{1}{f(0) c_*^{2/(\gamma-1)}}\int_{0}^{x}\arctan f'(s)
{\sqrt{1+(f'(s))^2}}ds,\quad x\in(0,l_1).
$$
Then there is no vacuum on $\{(x,f(x)):l_-\le x\le \check x\}$
for the global smooth transonic flow to the problem \eqref{sonic1}--\eqref{sonic7}.
If
$$
f''(x)=0,\quad x>\check x,
$$
then there is no vacuum for the global smooth transonic flow to the problem \eqref{sonic1}--\eqref{sonic7}.

{\rm (ii)} If
$\arctan f'(l_1)\ge\int_{-\infty}^0 b^{1/2}(s)ds$,
then the global smooth transonic flow to the problem \eqref{sonic1}--\eqref{sonic7}
admits a vacuum.
Furthermore, the first vacuum point in the increasing $x$-direction must form at the upper wall before
$$
\hat x=\inf\Big\{x\in(0,l_1):\arctan f'(x)\ge\int_{-\infty}^0 b^{1/2}(s)ds\Big\}.
$$

{\rm (iii)} If
\begin{align*}
\left\{
\begin{aligned}
&{\lim_{x\to+\infty}}f''(x)x^{2\gamma/(\gamma+1)}=+\infty,
\quad&&\mbox{ when } l_1=+\infty \mbox{ and }f'(+\infty)<+\infty,
\\
&{\lim_{x\to+l_1}}\frac{f''(x)}{(f'(x))^3}f^{2\gamma/(\gamma+1)}(x)=+\infty,
\quad&&\mbox{ when } l_1\le+\infty \mbox{ and }f'(l_1)=+\infty,
\end{aligned}
\right.
\end{align*}
then the global smooth transonic flow to the problem \eqref{sonic1}--\eqref{sonic7}
admits a vacuum.
\end{theorem}

\section{Smooth supersonic flows before vacuum formation}

In this section, we consider the smooth supersonic flow problem \eqref{phy1}--\eqref{phy5}
before vacuum formation.
For the existence, we will solve the problem \eqref{newsupp-eq}--\eqref{newsupp-q}, which
is equivalent to the problem \eqref{phy1}--\eqref{phy5} when there is no vacuum,
by a fixed point argument.
Since \eqref{newsupp-eq} is singular at the vacuum,
a cut-off technique is needed.
For $0<\varepsilon<1/{\mu_2}$, let $H_\varepsilon\in C^\infty(\mathbb R)$ be a cut-off function such that
$$
H_\varepsilon(s)
\left\{
\begin{aligned}
&=s,\quad&&\mbox{if }s\ge-\frac1\varepsilon,
\\
&=-\frac2\varepsilon,\quad&&\mbox{if }s\le-\frac2\varepsilon
\end{aligned}
\right.
\quad
\mbox{ and }
\quad
0\le H'_\varepsilon(s)\le2\mbox{ for } s\in\mathbb R.
$$

\subsection{Semilinear problems}

Give ${\mathscr Q}\in C([l_0,l_1))$ satisfying
\begin{align}
\label{sss-5-0}
c_*\le{\mathscr Q}(x)\le c^*,\quad l_0<x<l_1.
\end{align}
Fix $\zeta>0$. Assume that $\tilde Q\in C^{0,1}([0,\zeta]\times[0,m])$ satisfies
\begin{gather}
\label{sss-5-1}
\tilde Q(\varphi,\psi)\le-\mu_1,\quad
(\varphi,\psi)\in(0,\zeta)\times(0,m)
\end{gather}
and
\begin{gather}
\label{sss-5}
-\beta\le\tilde Q_\varphi(\varphi,\psi)\le0,\quad
|\tilde Q_\psi(\varphi,\psi)|\le\beta,\quad
(\varphi,\psi)\in(0,\zeta)\times(0,m),
\end{gather}
where $\beta$ is a positive constant to be determined.
Define
\begin{align*}
&\tilde W(\varphi,\psi)=
\pd{H_\varepsilon(\tilde Q)}{\varphi}(\varphi,\psi)-b^{1/2}(H_\varepsilon(\tilde Q(\varphi,\psi)))
\pd{H_\varepsilon(\tilde Q)}{\psi}(\varphi,\psi),
&&\quad(\varphi,\psi)\in(0,\zeta)\times(0,m),
\\
&\tilde Z(\varphi,\psi)=
-\pd{H_\varepsilon(\tilde Q)}{\varphi}(\varphi,\psi)-b^{1/2}(H_\varepsilon(\tilde Q(\varphi,\psi)))
\pd{H_\varepsilon(\tilde Q)}{\psi}(\varphi,\psi),
&&\quad(\varphi,\psi)\in(0,\zeta)\times(0,m).
\end{align*}
By \eqref{sss-5-1}, \eqref{sss-5} and the definition of $H_\varepsilon$,
there exists a positive constant $M_1=M_1(\gamma,\mu_1)$
such that
\begin{align}
\label{a-13-0}
|\tilde W(\varphi,\psi)|\le M_1\beta,\quad
|\tilde Z(\varphi,\psi)|\le M_1\beta,
\quad(\varphi,\psi)\in(0,\zeta)\times(0,m).
\end{align}
Here and hereafter, a parenthesis after a generic constant means that
this constant depends only on the variables in the parentheses.
Consider the problem
\begin{align}
\label{b-1}
&W_\varphi+b^{1/2}(H_\varepsilon(\tilde Q))W_\psi
=\frac14b^{-1}(H_\varepsilon(\tilde Q))p(H_\varepsilon(\tilde Q))W(W+\tilde Z),
\quad&&(\varphi,\psi)\in(0,\zeta)\times(0,m),
\\
\label{b-2}
&Z_\varphi-b^{1/2}(H_\varepsilon(\tilde Q))Z_\psi
=-\frac14b^{-1}(H_\varepsilon(\tilde Q))p(H_\varepsilon(\tilde Q))Z(\tilde W+Z),
\quad&&(\varphi,\psi)\in(0,\zeta)\times(0,m),
\\
\label{b-3}
&W(0,\psi)=W_0(\psi),\quad&&\psi\in(0,m),
\\
\label{b-4}
&Z(0,\psi)=Z_0(\psi),\quad&&\psi\in(0,m),
\\
\label{b-5}
&W(\varphi,0)+Z(\varphi,0)=0,\quad&&\varphi\in(0,\zeta),
\\
\label{b-6}
&W(\varphi,m)+Z(\varphi,m)=\tilde h(\varphi),
\quad&&\varphi\in(0,\zeta),
\\
\label{b-7}
&Q_\varphi=\frac{1}{2}(W-Z),
\quad&&(\varphi,\psi)\in(0,\zeta)\times(0,m),
\\
\label{b-9}
&Q(0,\psi)=Q_0(\psi),\quad&&\psi\in(0,m),
\end{align}
where
\begin{align*}
\tilde h(\varphi)=\frac{2f''(X_{\text{\rm up}}(\varphi))(1+(f'(X_{\text{\rm up}}(\varphi)))^2)^{-3/2}}
{b^{1/2}(H_\varepsilon(\tilde Q(\varphi,m)))A^{-1}_+(\tilde Q(\varphi,m))},
\quad\varphi\in(0,\zeta).
\end{align*}
It follows from \eqref{sss-0}, \eqref{sss-5-1} and the definition of $H_\varepsilon$ that
there exists a positive constant $M_2=M_2(\gamma,\zeta,\varepsilon,f)$
such that
\begin{align}
\label{a-13}
0\le\tilde h(\varphi)\le M_2,\quad
\varphi\in(0,\zeta).
\end{align}
Indeed, $M_2$ depends only on $\gamma$, $\varepsilon$
and $\sup_{(l_0,X_{\text{\rm up}}(\zeta))}f''(1+(f')^2)^{-3/2}$.

\begin{lemma}
\label{sss-prop1}
Assume that $W_0,Z_0\in L^\infty(0,m)$ satisfy \eqref{sss-3}.
Then the problem \eqref{b-1}--\eqref{b-6} admits a unique weak solution
$(W,Z)\in L^\infty((0,\zeta)\times(0,m))\times L^\infty((0,\zeta)\times(0,m))$.
Furthermore, the solution satisfies
\begin{align}
\label{sss-prop1-00}
-\mu_4(1+\mu_3)\le W(\varphi,\psi)\le 0,\quad
0\le Z(\varphi,\psi)\le \mu_4(1+\mu_3),\quad
(\varphi,\psi)\in(0,\zeta)\times(0,m),
\end{align}
where $\mu_4=\mu_4(\gamma,m,\mu_1,\zeta,\varepsilon,f)$ is a positive constant.
\end{lemma}

\Proof
Since the system \eqref{b-1}, \eqref{b-2} is strictly hyperbolic
according to \eqref{sss-5-1} and the definition of $H_\varepsilon$,
the classical theory yields the local existence and the uniqueness of the weak solution
to the problem \eqref{b-1}--\eqref{b-6}.
By \eqref{sss-3} and \eqref{a-13}, one can get the global existence.
Furthermore, the solution satisfies
\begin{align}
\label{sss-prop1-0}
W(\varphi,\psi)\le 0,\quad
Z(\varphi,\psi)\ge0,\quad
(\varphi,\psi)\in(0,\zeta)\times(0,m).
\end{align}

Assume that
$$
\Sigma_{\pm}:
\Psi'_{\pm}(\varphi)=\pm b^{1/2}(H_\varepsilon(\tilde Q(\varphi,\Psi_{\pm}(\varphi)))),
\quad
0<\Psi_{\pm}(\varphi)<m,\quad\hat\varphi_{\pm}<\varphi<\check\varphi_{\pm}
\quad(0\le\hat\varphi_{\pm}<\check\varphi_{\pm}\le\zeta)
$$
are a positive and a negative characteristics, respectively.
By \eqref{b-1} and \eqref{b-2}, $W$ and $Z$ satisfy
\begin{align}
\label{sss-prop1-1}
\frac{d}{d\varphi}W(\varphi,\Psi_{+}(\varphi))\ge
\tilde R_+(\varphi)W(\varphi,\Psi_{+}(\varphi)),
\quad\hat\varphi_+<\varphi<\check\varphi_+
\end{align}
and
\begin{align}
\label{sss-prop1-2}
\frac{d}{d\varphi}Z(\varphi,\Psi_{-}(\varphi))\le
\tilde R_-(\varphi)Z(\varphi,\Psi_{-}(\varphi)),
\quad\hat\varphi_-<\varphi<\check\varphi_-
\end{align}
on $\Sigma_{\pm}$, respectively, where
\begin{align}
\label{sss-prop1-4}
\tilde R_\pm(\varphi)=
\frac{d}{d\varphi}F(H_\varepsilon(\tilde Q(\varphi,\Psi_{\pm}(\varphi)))),\quad
\varphi\in(\hat\varphi_{\pm},\check\varphi_{\pm})
\end{align}
with
$$
F(s)=\frac14\int_s^{-1}b^{-1}(t)p(t)dt,\quad s<0.
$$
If follows from \eqref{sss-prop1-0}--\eqref{sss-prop1-4} that
\begin{align}
\label{sss-prop1-5}
W(\check\varphi,\Psi_{+}(\check\varphi))
\ge M_3W(\hat\varphi,\Psi_{+}(\hat\varphi)),
\quad
Z(\check\varphi,\Psi_{-}(\check\varphi))
\le M_3Z(\hat\varphi,\Psi_{-}(\hat\varphi))
\end{align}
with
$M_3=\mbox{exp}\big\{{F(-{2}/{\varepsilon})}-{F(-\mu_1)}\big\}>1$
depending only on $\gamma$, $\mu_1$ and $\varepsilon$.

\vskip10mm

\hskip45mm
\setlength{\unitlength}{0.6mm}
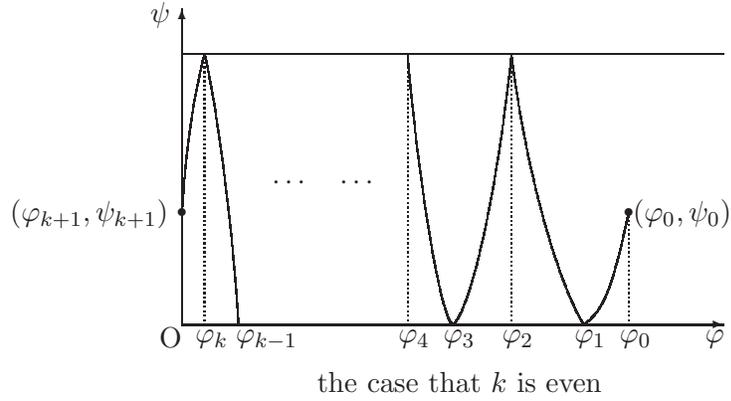
\begin{picture}(250,65)
\put(0,0){\vector(1,0){120}}
\put(0,0){\vector(0,1){70}}
\put(116,-4){$\varphi$} \put(-7,67){$\psi$}
\put(0,60){\qbezier(0,0)(55,0)(120,0)}

\put(0,25){\cbezier(0,0)(1,15)(2.5,25)(5,35)}

\put(5,60){\cbezier(0,0)(2.5,-10)(6.5,-40)(7.5,-60)}

\put(20,30){$\cdots\quad\cdots$}

\put(60,0){\cbezier(0,0)(-3,2)(-8,35)(-10,60)}
\put(60,0){\cbezier(0,0)(4,3)(10,30)(13,60)}

\put(73,60){\cbezier(0,0)(3,-25)(10,-50)(16,-60)}

\put(89,0){\cbezier(0,0)(4,3)(6,6)(10,25)}

\put(99,0){\qbezier[25](0,0)(0,12.5)(0,25)}
\put(73,0){\qbezier[60](0,0)(0,30)(0,60)}
\put(50,0){\qbezier[60](0,0)(0,30)(0,60)}
\put(5,0){\qbezier[60](0,0)(0,30)(0,60)}

\put(99,25){\circle*{1.5}}
\put(100,23){$(\varphi_0,\psi_0)$}

\put(97,-4){$\varphi_0$}
\put(87,-4){$\varphi_1$}
\put(71,-4){$\varphi_2$}
\put(58,-4){$\varphi_3$}
\put(48,-4){$\varphi_4$}

\put(12,-4){$\varphi_{k-1}$}
\put(3,-4){$\varphi_k$}
\put(-5,-5){O}

\put(0,25){\circle*{1.5}}
\put(-38,23){$(\varphi_{k+1},\psi_{k+1})$}

\put(30,-15){the case that $k$ is even}

\end{picture}

\vskip20mm

\hskip45mm
\setlength{\unitlength}{0.6mm}
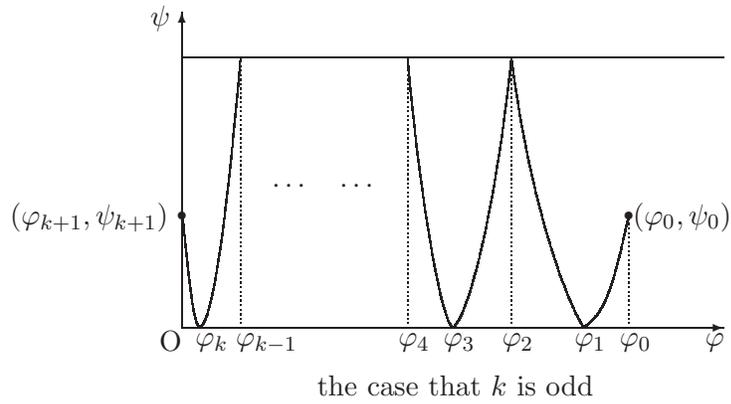
\begin{picture}(250,65)
\put(0,0){\vector(1,0){120}}
\put(0,0){\vector(0,1){70}}
\put(116,-4){$\varphi$} \put(-7,67){$\psi$}
\put(0,60){\qbezier(0,0)(55,0)(120,0)}

\put(0,25){\cbezier(0,0)(1,-6)(2,-24)(4,-25)}

\put(4,0){\cbezier(0,0)(3,2)(6,25)(9,60)}

\put(20,30){$\cdots\quad\cdots$}

\put(60,0){\cbezier(0,0)(-3,2)(-8,35)(-10,60)}
\put(60,0){\cbezier(0,0)(4,3)(10,30)(13,60)}

\put(73,60){\cbezier(0,0)(3,-25)(10,-50)(16,-60)}

\put(89,0){\cbezier(0,0)(4,3)(6,6)(10,25)}

\put(99,0){\qbezier[25](0,0)(0,12.5)(0,25)}
\put(73,0){\qbezier[60](0,0)(0,30)(0,60)}
\put(50,0){\qbezier[60](0,0)(0,30)(0,60)}
\put(13,0){\qbezier[60](0,0)(0,30)(0,60)}

\put(99,25){\circle*{1.5}}
\put(100,23){$(\varphi_0,\psi_0)$}

\put(97,-4){$\varphi_0$}
\put(87,-4){$\varphi_1$}
\put(71,-4){$\varphi_2$}
\put(58,-4){$\varphi_3$}
\put(48,-4){$\varphi_4$}

\put(12,-4){$\varphi_{k-1}$}
\put(3,-4){$\varphi_k$}
\put(-5,-5){O}

\put(0,25){\circle*{1.5}}
\put(-38,23){$(\varphi_{k+1},\psi_{k+1})$}

\put(30,-15){the case that $k$ is odd}

\end{picture}
\vskip12mm

We now estimate the lower bound of $W$ by the method of characteristics.
Fix $(\varphi_0,\psi_0)\in(0,\zeta)\times(0,m)$.
Let $\psi=\Psi_1(\varphi)$ be the positive characteristic from $(\varphi_0,\psi_0)$,
which approaches either $\{0\}\times[0,m]$
or $(0,\zeta)\times\{0\}$ at a point
$(\varphi_1,\psi_1)$.
If $\varphi_1>0$, then there exists a negative characteristic
$\psi=\Psi_2(\varphi)$ from $(\varphi_1,\psi_1)$,
which approaches either $\{0\}\times[0,m]$
or $(0,\zeta)\times\{m\}$ at a point
$(\varphi_2,\psi_2)$.
If $\varphi_2>0$, then there exists a positive characteristic
$\psi=\Psi_3(\varphi)$ from $(\varphi_2,\psi_2)$,
which approaches either $\{0\}\times[0,m]$
or $(0,\zeta)\times\{0\}$ at a point
$(\varphi_3,\psi_3)$.
Therefore, there exists a nonnegative integer $k$ such that
$$
\varphi_0>\varphi_1>\cdots>\varphi_k>\varphi_{k+1}=0,
$$
$$
\psi_j=\left\{
\begin{aligned}
&0,\quad1\le j\le k\mbox{ and $j$ is odd},
\\
&m,\quad1\le j\le k\mbox{ and $j$ is even},
\end{aligned}
\right.
\quad0\le\psi_{k+1}\le m,
$$
$$
\left\{
\begin{aligned}
&\Psi'_{j}(\varphi)=(-1)^{j-1}b^{1/2}(H_\varepsilon(\tilde Q(\varphi,\Psi_{j}(\varphi)))),
\quad\varphi_{j}<\varphi<\varphi_{j-1},
\\
&\Psi_{j}(\varphi_j)=\psi_{j},\quad \Psi_{j}(\varphi_{j-1})=\psi_{j-1},
\end{aligned}
\right.
\quad1\le j\le k+1.
$$
By \eqref{sss-5-1},
there exists an integer $K=K(\gamma,m,\mu_1,\zeta)$ such that $k\le K$.
For $1\le j\le k+1$, it follows from \eqref{sss-prop1-5} that
\begin{align}
\label{refin1}
W(\varphi_{j-1},\psi_{j-1})\ge M_3W(\varphi_j,\psi_{j})\mbox{ if $j$ is odd},
\quad
Z(\varphi_{j-1},\psi_{j-1})\le M_3Z(\varphi_j,\psi_{j})\mbox{ if $j$ is even}.
\end{align}
If $k$ is even, one can get from \eqref{refin1}, \eqref{b-3}--\eqref{b-6} that
\begin{align}
\label{a-11}
W(\varphi_{0},\psi_{0})
\ge&M_3W(\varphi_1,\psi_{1})
=-M_3Z(\varphi_1,\psi_{1})
\nonumber
\\
\ge&
-M_3^2Z(\varphi_2,\psi_{2})
=M_3^2W(\varphi_2,\psi_{2})-M_3^2\tilde h(\varphi_2,\psi_{2})
\nonumber
\\
\ge&
M_3^{k+1}W_0(\psi_{k+1})
-\sum_{2\le2j\le k}M_3^{2j}
\tilde h(\varphi_{2j},\psi_{2j}).
\end{align}
Similarly, if $k$ is odd, then
\begin{align}
\label{a-12}
W(\varphi_{0},\psi_{0})\ge
-M_3^{k+1}Z_0(\psi_{k+1})
-\sum_{2\le2j\le k}M_3^{2j}
\tilde h(\varphi_{2j},\psi_{2j}).
\end{align}
It follows from \eqref{a-11}, \eqref{a-12}, \eqref{sss-3} and \eqref{a-13} that
in both cases,
\begin{align*}
W(\varphi_{0},\psi_{0})\ge-\mu_3 M_3^{K+1}-KM_3^{K}M_2,
\end{align*}
which yields the lower bound of $W$ in \eqref{sss-prop1-00} by choosing
$\mu_4=M_3^{K}(M_3+KM_2)$.
The upper bound of $Z$ in \eqref{sss-prop1-00} can be obtained similarly.
$\hfill\Box$\vskip 4mm

\begin{lemma}
\label{sss-prop2}
Let $(W,Z)\in L^\infty((0,\zeta)\times(0,m))\times L^\infty((0,\zeta)\times(0,m))$
be the weak solution to the problem \eqref{b-1}--\eqref{b-6} with
$W_0,Z_0\in L^\infty(0,m)$ satisfying \eqref{sss-3}.
Assume that $Q_0\in L^{\infty}(0,m)$ satisfies \eqref{sss-3-0} and
\begin{align}
\label{zz-2a}
Q_0(\psi)\ge-\frac1\varepsilon,\quad 0<\psi<m.
\end{align}
Then the problem \eqref{b-7}, \eqref{b-9} admits a unique weak solution
$Q\in L^\infty((0,\zeta)\times(0,m))$.
Furthermore, the solution $Q\in C^{0,1}([0,\zeta]\times[0,m])$ and it satisfies
\begin{gather}
\label{sss-prop2-1-0}
-\mu_2-\mu_4(1+\mu_3)\zeta\le Q(\varphi,\psi)\le-\mu_1,\quad
(\varphi,\psi)\in(0,\zeta)\times(0,m),
\\
\label{sss-prop2-1}
-\mu_4(1+\mu_3)\le Q_\varphi(\varphi,\psi)\le0,\quad
(\varphi,\psi)\in(0,\zeta)\times(0,m),
\\
\label{sss-prop2-2}
|Q_\psi(\varphi,\psi)|\le \mu_5(1+\mu_3^2)(1+\beta\varphi),\quad
(\varphi,\psi)\in(0,\zeta)\times(0,m),
\end{gather}
where $\mu_5=\mu_5(\gamma,m,\mu_1,\zeta,\varepsilon,f)$ is a positive constant.
\end{lemma}

\Proof
The existence and uniqueness of the weak solution to the problem \eqref{b-7}, \eqref{b-9} are clear.
Moreover, \eqref{sss-prop2-1} follows from \eqref{b-7} and \eqref{sss-prop1-00},
while \eqref{sss-prop2-1-0} follows from \eqref{sss-3-0} and \eqref{sss-prop2-1}, directly.
It remains to prove \eqref{sss-prop2-2}.
Rewrite \eqref{b-1} as
\begin{align}
\label{ab-1}
W_\psi=-b^{-1/2}(H_\varepsilon(\tilde Q))W_\varphi+
\frac14b^{-3/2}(H_\varepsilon(\tilde Q))p(H_\varepsilon(\tilde Q))W(W+\tilde Z),
\quad(\varphi,\psi)\in(0,\zeta)\times(0,m).
\end{align}
Fix $\hat\varphi\in(0,\zeta)$ and $\psi\in(0,m)$. Integrating \eqref{ab-1}
over $(0,\hat\varphi)$, one gets from \eqref{sss-5}, \eqref{a-13-0},
\eqref{sss-prop1-00} and the definition of $H_\varepsilon$
that
\begin{align}
\label{zz-1}
&\Big|\frac{d}{d\psi}\int_{0}^{\hat\varphi}W(\varphi,\psi)d\varphi\Big|
\nonumber
\\
=&\Big|-b^{-1/2}(H_\varepsilon(\tilde Q(\hat\varphi,\psi)))W(\hat\varphi,\psi)
+b^{-1/2}(H_\varepsilon(\tilde Q(0,\psi)))W(0,\psi)
\nonumber
\\
&\qquad
-\frac12\int_{0}^{\hat\varphi}
b^{-3/2}(H_\varepsilon(\tilde Q(\varphi,\psi)))p(H_\varepsilon(\tilde Q(\varphi,\psi)))
H'_\varepsilon(\tilde Q(\varphi,\psi))
\tilde Q_\varphi(\varphi,\psi)W(\varphi,\psi)d\varphi
\nonumber
\\
&\qquad
+\frac14\int_{0}^{\hat\varphi}
b^{-3/2}(H_\varepsilon(\tilde Q(\varphi,\psi)))p(H_\varepsilon(\tilde Q(\varphi,\psi)))
W(\varphi,\psi)(W(\varphi,\psi)+\tilde Z(\varphi,\psi))d\varphi\Big|
\nonumber
\\
\le&\mu_4(1+\mu_3)\big(b^{-1/2}(H_\varepsilon(\tilde Q(\hat\varphi,\psi)))
+b^{-1/2}(H_\varepsilon(\tilde Q(0,\psi)))\big)
\nonumber
\\
&\qquad
-\frac12\mu_4(1+\mu_3)\int_{0}^{\hat\varphi}
b^{-3/2}(H_\varepsilon(\tilde Q(\varphi,\psi)))p(H_\varepsilon(\tilde Q(\varphi,\psi)))
H'_\varepsilon(\tilde Q(\varphi,\psi))
\tilde Q_\varphi(\varphi,\psi)d\varphi
\nonumber
\\
&\qquad
+\frac14{\mu_4(1+\mu_3)}\int_{0}^{\hat\varphi}
b^{-3/2}(H_\varepsilon(\tilde Q(\varphi,\psi)))p(H_\varepsilon(\tilde Q(\varphi,\psi)))
(|W(\varphi,\psi)|+|\tilde Z(\varphi,\psi)|)d\varphi
\nonumber
\\
\le&\mu_4(1+\mu_3)\big(b^{-1/2}(H_\varepsilon(\tilde Q(\hat\varphi,\psi)))
+b^{-1/2}(H_\varepsilon(\tilde Q(0,\psi)))\big)
\nonumber
\\
&\qquad
+\mu_4(1+\mu_3)\big(b^{-1/2}(H_\varepsilon(\tilde Q(\hat\varphi,\psi)))
-b^{-1/2}(H_\varepsilon(\tilde Q(0,\psi)))\big)
\nonumber
\\
&\qquad
+\frac{1}4\mu_4(1+\mu_3)(\mu_4(1+\mu_3)+M_1\beta)\int_{0}^{\hat\varphi}
b^{-3/2}(H_\varepsilon(\tilde Q(\varphi,\psi)))p(H_\varepsilon(\tilde Q(\varphi,\psi)))d\varphi.
\end{align}
A similar argument shows
\begin{align}
\label{zz-2}
&\Big|\frac{d}{d\psi}\int_{0}^{\hat\varphi}Z(\varphi,\psi)d\varphi\Big|
\nonumber
\\
\le&\mu_4(1+\mu_3)\big(b^{-1/2}(H_\varepsilon(\tilde Q(\hat\varphi,\psi)))
+b^{-1/2}(H_\varepsilon(\tilde Q(0,\psi)))\big)
\nonumber
\\
&\qquad
+\mu_4(1+\mu_3)\big(b^{-1/2}(H_\varepsilon(\tilde Q(\hat\varphi,\psi)))
-b^{-1/2}(H_\varepsilon(\tilde Q(0,\psi)))\big)
\nonumber
\\
&\qquad
+\frac{1}4\mu_4(1+\mu_3)(\mu_4(1+\mu_3)+M_1\beta)\int_{0}^{\hat\varphi}
b^{-3/2}(H_\varepsilon(\tilde Q(\varphi,\psi)))p(H_\varepsilon(\tilde Q(\varphi,\psi)))d\varphi.
\end{align}
It follows from \eqref{zz-2a} and \eqref{sss-3} that
\begin{align}
\label{zz-2aa}
|Q'_0(\psi)|\le M_4\mu_3,
\end{align}
where $M_4=M_4(\gamma,\varepsilon)$ is a positive constant.
By \eqref{b-7}, \eqref{zz-1}--\eqref{zz-2aa},
there exists a positive constant $\mu_5=\mu_5(\gamma,m,\mu_1,\zeta,\varepsilon,f)$
such that
\begin{align*}
|Q_\psi(\hat\varphi,\psi)|
=\Big|Q'_0(\psi)+\frac{1}{2}\frac{d}{d\psi}\int_{0}^{\hat\varphi}(W(\varphi,\psi)-Z(\varphi,\psi))d\varphi\Big|
\le\mu_5(1+\mu_3^2)(1+\beta\hat\varphi),
\end{align*}
which implies \eqref{sss-prop2-2}.
$\hfill\Box$\vskip 4mm

\begin{remark}
Precisely, $\mu_4$ and $\mu_5$ depend only on $\gamma$, $m$, $\mu_1$, $\zeta$, $\varepsilon$
and $\sup_{(l_0,X_{\text{\rm up}}(\zeta))}f''(1+(f')^2)^{-3/2}$.
\end{remark}

\subsection{Existence and uniqueness}

\begin{theorem}
\label{sss-thm1}
Assume that $f\in C^{2}([l_0,l_1))$ satisfies \eqref{sss-0},
and $Q_0\in C^{1}([0,m])$ and $G_0\in C([0,m])$
satisfy \eqref{sss-1} and \eqref{sss-2}.
Then the problem \eqref{newsupp-eq}--\eqref{newsupp-q} admits uniquely a solution
$Q\in C^{1}([0,\zeta)\times[0,m])$,
where either
$\zeta=+\infty$
or
$0<\zeta<+\infty$ with $\lim_{\varphi\to\zeta^-}\min_{[0,m]}Q(\varphi,\cdot)=-\infty$.
Furthermore, the solution satisfies
\begin{align*}
-Q_\varphi(\varphi,\psi)\pm b^{1/2}(Q(\varphi,\psi))Q_\psi(\varphi,\psi)\ge0,\quad
(\varphi,\psi)\in(0,\zeta)\times(0,m).
\end{align*}
\end{theorem}

\Proof
The uniqueness result follows from the uniqueness
of the classical solution to the problem \eqref{phy1}--\eqref{phy5},
which is a standard strictly hyperbolic problem,
since the problem \eqref{newsupp-eq}--\eqref{newsupp-q}
and the problem \eqref{phy1}--\eqref{phy5} are equivalent for supersonic flows without vacuum.

For the existence,
we claim that for any $0<\varepsilon<1/\mu_2$ and any $\hat\zeta>0$,
there exists $\hat\zeta_\varepsilon\in(0,\hat\zeta]$ such that
the problem \eqref{newsupp-neq1}--\eqref{newsupp-q4} with $\zeta=\hat\zeta_\varepsilon$
admits a weak solution
$(W,Z,Q)\in L^\infty((0,\hat\zeta_\varepsilon)\times(0,m))\times
L^\infty((0,\hat\zeta_\varepsilon)\times(0,m))
\times C^{0,1}([0,\hat\zeta_\varepsilon]\times[0,m])$
satisfying
\begin{align*}
Q(\varphi,\psi)>-\frac1\varepsilon,\quad
W(\varphi,\psi)\le0,\quad Z(\varphi,\psi)\ge0,\quad
(\varphi,\psi)\in(0,\hat\zeta_\varepsilon)\times(0,m)
\end{align*}
and either
$\hat\zeta_\varepsilon=\hat\zeta$
or
$0<\hat\zeta_\varepsilon<\hat\zeta$ with $\min_{[0,m]}Q(\hat\zeta_\varepsilon,\cdot)=-1/\varepsilon$.
Assuming the claim, we complete the proof of the theorem first.
Since $\hat\zeta>0$ is arbitrary,
there exists $0<\zeta_\varepsilon\le+\infty$ such that
the problem \eqref{newsupp-neq1}--\eqref{newsupp-q4} with $\zeta=\zeta_\varepsilon$
admits a weak solution
$(W,Z,Q)\in L^{\infty}((0,\zeta_\varepsilon)\times(0,m))\times
L^{\infty}((0,\zeta_\varepsilon)\times(0,m))
\times C^{0,1}([0,\zeta_\varepsilon)\times[0,m])$
satisfying
\begin{align*}
Q(\varphi,\psi)>-\frac1\varepsilon,\quad
W(\varphi,\psi)\le0,\quad Z(\varphi,\psi)\ge0,\quad
(\varphi,\psi)\in(0,\zeta_\varepsilon)\times(0,m)
\end{align*}
and either
$\zeta_\varepsilon=+\infty$
or
$0<\zeta_\varepsilon<+\infty$ with
$\lim_{\varphi\to\zeta_\varepsilon^-}\min_{[0,m]}Q(\varphi,\cdot)=-1/\varepsilon$.
Equivalently, $Q$ solves the problem \eqref{newsupp-eq}--\eqref{newsupp-q}
with $\zeta=\zeta_\varepsilon$.
The regularity of initial and boundary data
and \eqref{sss-2} yield that
$Q\in C^{1}([0,\zeta_\varepsilon)\times[0,m])$.
Then, the theorem follows from the arbitrariness of $\varepsilon\in(0,1/\mu_2)$.

We now turn to the proof of the claim. Fix $0<\varepsilon<1/\mu_2$ and $\hat\zeta>0$.
Let $\mu_4$ and $\mu_5$ be the constants determined in Lemmas \ref{sss-prop1} and \ref{sss-prop2}
with $\zeta=\hat\zeta$,
which depend only on $\gamma$, $m$, $\mu_1$, $\hat\zeta$, $\varepsilon$ and $f$.
The proof consists of four steps.

{\bf Step I}\quad Local existence.

Choose
\begin{align}
\label{sss-5-00}
\beta=\max\big\{\mu_4(1+\mu_3),2\mu_5(1+\mu_3^2)\big\}.
\end{align}
For ${\mathscr Q}\in C([l_0,l_1))$ satisfying \eqref{sss-5-0}
and $\tilde Q\in C^{0,1}([0,\hat\zeta]\times[0,m])$
satisfying  \eqref{sss-5-1} and \eqref{sss-5} with $\zeta=\hat\zeta$ and \eqref{sss-5-00},
Lemmas \ref{sss-prop1} and \ref{sss-prop2} show that
the problem \eqref{b-1}--\eqref{b-9} with $\zeta=\hat\zeta$ admits a unique weak solution
$(W,Z,Q)\in L^\infty((0,\hat\zeta)\times(0,m))\times L^\infty((0,\hat\zeta)\times(0,m))
\times C^{0,1}([0,\hat\zeta]\times[0,m])$
and $(W,Z,Q)$ satisfies
\begin{gather}
\label{sss-thm1-1}
-\mu_4(1+\mu_3)\le W(\varphi,\psi)\le 0,\quad
0\le Z(\varphi,\psi)\le\mu_4(1+\mu_3),\quad
(\varphi,\psi)\in(0,\hat\zeta)\times(0,m),
\\
\label{sss-thm1-2-0}
-\mu_2-\mu_4(1+\mu_3)\hat\zeta\le Q(\varphi,\psi)\le-\mu_1,\quad
(\varphi,\psi)\in(0,\hat\zeta)\times(0,m),
\\
\label{sss-thm1-2}
-\mu_4(1+\mu_3)\le Q_\varphi(\varphi,\psi)\le0,\quad
(\varphi,\psi)\in(0,\hat\zeta)\times(0,m),
\\
\label{sss-thm1-3}
|Q_\psi(\varphi,\psi)|\le \mu_5(1+\mu_3^2)(1+\beta\varphi),\quad
(\varphi,\psi)\in(0,\hat\zeta)\times(0,m).
\end{gather}
Let
$\hat\zeta_0=\min\big\{\hat\zeta,1/{\beta}\big\}$.
Then, it follows from \eqref{sss-thm1-2} and \eqref{sss-thm1-3} that
\begin{align}
\label{sss-thm1-4}
-\beta\le Q_\varphi(\varphi,\psi)\le0,\quad
|Q_\psi(\varphi,\psi)|\le\beta,\quad
(\varphi,\psi)\in(0,\hat\zeta_0)\times(0,m).
\end{align}
Set
$$
\hat{\mathscr Q}(x)=\left\{
\begin{aligned}
&A^{-1}_+(Q(\Phi_{\text{\rm up}}(x),m)),\quad&&l_0\le x\le X_{\text{\rm up}}(\hat\zeta_0),
\\
&A^{-1}_+(Q(\hat\zeta_0,m)),\quad&& X_{\text{\rm up}}(\hat\zeta_0)<x<l_1,
\end{aligned}
\right.
$$
where $\Phi_{\text{\rm up}}$ and $X_{\text{\rm up}}$ are defined in $\S\, 2.3$.
Define
\begin{align*}
{\mathscr S}_0=&\Big\{({\mathscr Q},Q)\in C([l_0,l_1))\times
C^{0,1}([0,\hat\zeta_0]\times[0,m]):
c_*\le{\mathscr Q}(x)\le c^*\mbox{ for }l_0<x<l_1,
\\
&\qquad-\mu_2-\mu_4(1+\mu_3)\hat\zeta\le Q(\varphi,\psi)\le-\mu_1,
-\beta\le Q_\varphi(\varphi,\psi)\le0
\\
&\qquad\mbox{ and }
|Q_\psi(\varphi,\psi)|\le\beta\mbox{ for }
(\varphi,\psi)\in(0,\hat\zeta_0)\times(0,m)\Big\}
\end{align*}
with the following norm
$$
\|({\mathscr Q},Q)\|_{{\mathscr S}_0}=\max\Big\{
\|{\mathscr Q}\|_{L^\infty(l_0,l_1)},
\|Q\|_{L^\infty((0,\hat\zeta_0)\times(0,m))}\Big\},\quad
({\mathscr Q},Q)\in {\mathscr S}_0.
$$
By \eqref{sss-thm1-2-0} and \eqref{sss-thm1-4}, $(\hat{\mathscr Q},Q)\in {\mathscr S}_0$.
Therefore, one can define a mapping as follows
$$
{\mathscr J}_0:{\mathscr S}_0\longrightarrow{\mathscr S}_0,
\quad({\mathscr Q},\tilde Q)\longmapsto (\hat{\mathscr Q},Q).
$$
By \eqref{sss-thm1-4},
${\mathscr J}_0$ is compact.
Furthermore, one can show that ${\mathscr J}_0$ is continuous
(see for example \cite{WX3} Theorem 5.1).
Then, it follows from the Schauder fixed point theorem that ${\mathscr J}_0$ admits a fixed point $({\mathscr Q},Q)$.
That is to say, there exists
$(W,Z,Q)\in L^\infty((0,\hat\zeta_0)\times(0,m))\times L^\infty((0,\hat\zeta_0)\times(0,m))
\times C^{0,1}([0,\hat\zeta_0]\times[0,m])$, which solves
\begin{align*}
&W_\varphi+b^{1/2}(H_\varepsilon(Q))W_\psi
=\frac14b^{-1}(H_\varepsilon(Q))p(H_\varepsilon(Q))W(W-{H'_\varepsilon(Q)}{Q_\varphi}
&&
\nonumber
\\
&\qquad\qquad\qquad\qquad\qquad\qquad
-b^{1/2}(H_\varepsilon(Q))
{H'_\varepsilon(Q)}{Q_\psi}),
\quad&&(\varphi,\psi)\in(0,\hat\zeta_0)\times(0,m),
\\
&Z_\varphi-b^{1/2}(H_\varepsilon(Q))Z_\psi
=-\frac14b^{-1}(H_\varepsilon(Q))p(H_\varepsilon(Q))Z({H'_\varepsilon(Q)}{Q_\varphi}
&&
\nonumber
\\
&\qquad\qquad\qquad\qquad\qquad\qquad
-b^{1/2}(H_\varepsilon(Q))
{H'_\varepsilon(Q)}{Q_\psi}+Z),
\quad&&(\varphi,\psi)\in(0,\hat\zeta_0)\times(0,m),
\\
&W(0,\psi)=W_0(\psi),\quad&&\psi\in(0,m),
\\
&Z(0,\psi)=Z_0(\psi),\quad&&\psi\in(0,m),
\\
&W(\varphi,0)+Z(\varphi,0)=0,\quad&&\varphi\in(0,\hat\zeta_0),
\\
&W(\varphi,m)+Z(\varphi,m)=\frac{2f''(X_{\text{\rm up}}(\varphi))(1+(f'(X_{\text{\rm up}}(\varphi)))^2)^{-3/2}}
{b^{1/2}(H_\varepsilon(Q(\varphi,m)))A^{-1}_+(Q(\varphi,m))},
\quad&&\varphi\in(0,\hat\zeta_0),
\\
&Q_\varphi=\frac{1}{2}(W-Z),
\quad&&(\varphi,\psi)\in(0,\hat\zeta_0)\times(0,m),
\\
&Q(0,\psi)=Q_0(\psi),\quad&&\psi\in(0,m),
\\
&{\mathscr Q}(x)=A^{-1}_+(Q(\Phi_{\text{\rm up}}(x),m)),
\quad&&x\in(l_0,X_{\text{\rm up}}(\hat\zeta_0)).
\end{align*}
Furthermore,
\eqref{sss-thm1-1} and \eqref{sss-thm1-2-0} lead to
\begin{gather}
\label{sss-thm1-1-1}
-\mu_4(1+\mu_3)\le W(\varphi,\psi)\le 0,\quad
0\le Z(\varphi,\psi)\le\mu_4(1+\mu_3),\quad
(\varphi,\psi)\in(0,\hat\zeta_0)\times(0,m),
\\
\label{sss-thm1-2-0-1}
-\mu_2-\mu_4(1+\mu_3)\hat\zeta\le Q(\varphi,\psi)\le-\mu_1,\quad
(\varphi,\psi)\in(0,\hat\zeta_0)\times(0,m).
\end{gather}
Let
$$
\hat\zeta_{0,\varepsilon}=\sup\Big\{\tilde\varphi\in(0,\hat\zeta_0):
Q(\varphi,\psi)>-\frac1\varepsilon
\mbox{ for each }(\varphi,\psi)\in(0,\tilde\varphi)\times(0,m)\Big\}.
$$
Then, for $(\varphi,\psi)\in(0,\hat\zeta_{0,\varepsilon})\times(0,m)$,
\begin{align*}
W(\varphi,\psi)=
Q_\varphi(\varphi,\psi)-b^{1/2}(Q(\varphi,\psi))
Q_\psi(\varphi,\psi),
\quad
Z(\varphi,\psi)=
-Q_\varphi(\varphi,\psi)-b^{1/2}(Q(\varphi,\psi))
Q_\psi(\varphi,\psi),
\end{align*}
and $(W,Z,Q)$ solves
\begin{align*}
&W_\varphi+b^{1/2}(Q)W_\psi
=\frac14b^{-1}(Q)p(Q)W(W+Z),
\quad&&(\varphi,\psi)\in(0,\hat\zeta_{0,\varepsilon})\times(0,m),
\\
&Z_\varphi-b^{1/2}(Q)Z_\psi
=-\frac14b^{-1}(Q)p(Q)Z(W+Z),
\quad&&(\varphi,\psi)\in(0,\hat\zeta_{0,\varepsilon})\times(0,m),
\\
&W(0,\psi)=W_0(\psi),\quad&&\psi\in(0,m),
\\
&Z(0,\psi)=Z_0(\psi),\quad&&\psi\in(0,m),
\\
&W(\varphi,0)+Z(\varphi,0)=0,\quad&&\varphi\in(0,\hat\zeta_{0,\varepsilon}),
\\
&W(\varphi,m)+Z(\varphi,m)=\frac{2f''(X_{\text{\rm up}}(\varphi))(1+(f'(X_{\text{\rm up}}(\varphi)))^2)^{-3/2}}
{b^{1/2}(Q(\varphi,m))A^{-1}_+(Q(\varphi,m))},
\quad&&\varphi\in(0,\hat\zeta_{0,\varepsilon}),
\\
&Q_\varphi=\frac{1}{2}(W-Z),
\quad&&(\varphi,\psi)\in(0,\hat\zeta_{0,\varepsilon})\times(0,m),
\\
&Q(0,\psi)=Q_0(\psi),\quad&&\psi\in(0,m),
\\
&{\mathscr Q}(x)=A^{-1}_+(Q(\Phi_{\text{\rm up}}(x),m)),
\quad&&x\in(l_0,X_{\text{\rm up}}(\hat\zeta_{0,\varepsilon})).
\end{align*}
It follows from \eqref{sss-thm1-1-1} and \eqref{sss-thm1-2-0-1} that
\begin{gather}
\label{zz-8-0}
-\mu_4(1+\mu_3)\le W(\varphi,\psi)\le 0,\quad
0\le Z(\varphi,\psi)\le\mu_4(1+\mu_3),\quad
(\varphi,\psi)\in(0,\hat\zeta_{0,\varepsilon})\times(0,m),
\\
\label{zz-8}
-\mu_2-\mu_4(1+\mu_3)\hat\zeta\le Q(\varphi,\psi)\le-\mu_1,\quad
\quad
(\varphi,\psi)\in(0,\hat\zeta_{0,\varepsilon})\times(0,m).
\end{gather}
If $\hat\zeta_{0,\varepsilon}=\hat\zeta$ or
$\min_{[0,m]}Q(\hat\zeta_{0,\varepsilon},\cdot)=-1/\varepsilon$,
then the claim is proved.
So, without loss of generality, it is assumed that
$\hat\zeta_{0,\varepsilon}<\hat\zeta$ and
$\min_{[0,m]}Q(\hat\zeta_{0,\varepsilon},\cdot)>-1/\varepsilon$,
which also leads to $\hat\zeta_{0,\varepsilon}=\hat\zeta_0$.

{\bf Step II}\quad Extension.

Set
\begin{align*}
\hat\mu_2=\mu_2+\mu_4(1+\mu_3)\hat\zeta,
\quad
\hat\mu_3=\mu_4(1+\mu_3),\quad
\hat\beta=\max\big\{\mu_4(1+\hat\mu_3),2\mu_5(1+\hat\mu_3^2)\big\}.
\end{align*}
Let
$\hat\zeta_1=\min\big\{\hat\zeta,\hat\zeta_0+{1}/{\hat\beta}\big\}$.
Define
\begin{align*}
{\mathscr S}_1=&\Big\{({\mathscr Q},Q)\in C([X_{\text{\rm up}}(\hat\zeta_0),l_1))\times
C^{0,1}([\hat\zeta_0,\hat\zeta_1]\times[0,m]):
c_*\le{\mathscr Q}(x)\le c^*\mbox{ for }X_{\text{\rm up}}(\hat\zeta_0)<x<l_1,
\\
&\qquad-\hat\mu_2-\mu_4(1+\hat\mu_3)\hat\zeta\le Q(\varphi,\psi)\le-\mu_1,
-\hat\beta\le Q_\varphi(\varphi,\psi)\le0
\\
&\qquad\mbox{ and }
|Q_\psi(\varphi,\psi)|\le\hat\beta\mbox{ for }
(\varphi,\psi)\in(\hat\zeta_0,\hat\zeta_1)\times(0,m)\Big\}.
\end{align*}
For given $({\mathscr Q},\tilde Q)\in{\mathscr S}_1$, consider the problem
\begin{align}
\label{zb-1}
&W_\varphi+b^{1/2}(H_\varepsilon(\tilde Q))W_\psi
=\frac14b^{-1}(H_\varepsilon(\tilde Q))p(H_\varepsilon(\tilde Q))W(W+\tilde Z),
\quad&&(\varphi,\psi)\in(\hat\zeta_0,\hat\zeta_1)\times(0,m),
\\
\label{zb-2}
&Z_\varphi-b^{1/2}(H_\varepsilon(\tilde Q))Z_\psi
=-\frac14b^{-1}(H_\varepsilon(\tilde Q))p(H_\varepsilon(\tilde Q))Z(\tilde W+Z),
\quad&&(\varphi,\psi)\in(\hat\zeta_0,\hat\zeta_1)\times(0,m),
\\
\label{zb-3}
&W(\hat\zeta_0,\psi)=W_1(\psi),\quad&&\psi\in(0,m),
\\
\label{zb-4}
&Z(\hat\zeta_0,\psi)=Z_1(\psi),\quad&&\psi\in(0,m),
\\
\label{zb-5}
&W(\varphi,0)+Z(\varphi,0)=0,\quad&&\varphi\in(\hat\zeta_0,\hat\zeta_1),
\\
\label{zb-6}
&W(\varphi,m)+Z(\varphi,m)=\tilde h(\varphi),
\quad&&\varphi\in(\hat\zeta_0,\hat\zeta_1),
\\
\label{zb-7}
&Q_\varphi=\frac{1}{2}(W-Z),
\quad&&(\varphi,\psi)\in(\hat\zeta_0,\hat\zeta_1)\times(0,m),
\\
\label{zb-9}
&Q(\hat\zeta_0,\psi)=Q_1(\psi),\quad&&\psi\in(0,m),
\end{align}
where
\begin{align*}
&\tilde W(\varphi,\psi)=
\pd{H_\varepsilon(\tilde Q)}{\varphi}(\varphi,\psi)-b^{1/2}(H_\varepsilon(\tilde Q(\varphi,\psi)))
\pd{H_\varepsilon(\tilde Q)}{\psi}(\varphi,\psi),
&&\quad(\varphi,\psi)\in(\hat\zeta_0,\hat\zeta_1)\times(0,m),
\\
&\tilde Z(\varphi,\psi)=
-\pd{H_\varepsilon(\tilde Q)}{\varphi}(\varphi,\psi)-b^{1/2}(H_\varepsilon(\tilde Q(\varphi,\psi)))
\pd{H_\varepsilon(\tilde Q)}{\psi}(\varphi,\psi),
&&\quad(\varphi,\psi)\in(\hat\zeta_0,\hat\zeta_1)\times(0,m)
\end{align*}
and
\begin{align*}
W_1(\psi)=W(\hat\zeta_0,\psi),\quad
Z_1(\psi)=Z(\hat\zeta_0,\psi),\quad
Q_1(\psi)=Q(\hat\zeta_0,\psi),\quad\psi\in(0,m).
\end{align*}
It follows from \eqref{zz-8-0} and \eqref{zz-8} that $Q_1\in C^{0,1}([0,m])$,
$W_1,Z_1\in L^{\infty}(0,m)$ such that
\begin{align}
\label{zz-8-1}
-\hat\mu_2\le Q_1(\psi)\le-\mu_1,\quad
-\hat\mu_3\le W_1(\psi)\le0,\quad
0\le Z_1(\psi)\le\hat\mu_3,\quad\psi\in(0,m).
\end{align}
Then, Lemmas \ref{sss-prop1} and \ref{sss-prop2} show that
the problem \eqref{zb-1}--\eqref{zb-9} admits a unique weak solution
$(W,Z,Q)\in L^\infty((\hat\zeta_0,\hat\zeta_1)\times(0,m))\times L^\infty((\hat\zeta_0,\hat\zeta_1)\times(0,m))
\times C^{0,1}([\hat\zeta_0,\hat\zeta_1]\times[0,m])$
and it satisfies
\begin{gather}
\label{zz-9}
-\mu_4(1+\hat\mu_3)\le W(\varphi,\psi)\le 0,\quad
0\le Z(\varphi,\psi)\le\mu_4(1+\hat\mu_3),\quad
(\varphi,\psi)\in(\hat\zeta_0,\hat\zeta_1)\times(0,m),
\\
\label{zz-10}
-\hat\mu_2-\mu_4(1+\hat\mu_3)\hat\zeta\le Q(\varphi,\psi)\le-\mu_1,\quad
(\varphi,\psi)\in(\hat\zeta_0,\hat\zeta_1)\times(0,m),
\\
\label{zz-11}
-\hat\beta\le-\mu_4(1+\hat\mu_3)\le Q_\varphi(\varphi,\psi)\le0,\quad
(\varphi,\psi)\in(\hat\zeta_0,\hat\zeta_1)\times(0,m),
\\
\label{zz-12}
|Q_\psi(\varphi,\psi)|\le \mu_5(1+\hat\mu_3^2)(1+\hat\beta(\varphi-\hat\zeta_0))\le\hat\beta,\quad
(\varphi,\psi)\in(\hat\zeta_0,\hat\zeta_1)\times(0,m).
\end{gather}
Therefore, one can define a mapping as follows
$$
{\mathscr J}_1:{\mathscr S}_1\longrightarrow{\mathscr S}_1,
\quad({\mathscr Q},\tilde Q)\longmapsto (\hat{\mathscr Q},Q),
$$
where
$$
\hat{\mathscr Q}(x)=\left\{
\begin{aligned}
&A^{-1}_+(Q(\Phi_{\text{\rm up}}(x),m)),\quad&&X_{\text{\rm up}}(\hat\zeta_0)\le x
\le X_{\text{\rm up}}(\hat\zeta_1),
\\
&A^{-1}_+(Q(\hat\zeta_1,m)),\quad&& X_{\text{\rm up}}(\hat\zeta_1)<x<l_1.
\end{aligned}
\right.
$$
Then, the Schauder fixed point theorem shows that
${\mathscr J}_1$ admits a fixed point $({\mathscr Q},Q)$.
That is to say, there exists
$(W,Z,Q)\in L^\infty((\hat\zeta_0,\hat\zeta_1)\times(0,m))
\times L^\infty((\hat\zeta_0,\hat\zeta_1)\times(0,m))
\times C^{0,1}([\hat\zeta_0,\hat\zeta_1]\times[0,m])$
satisfying \eqref{zz-9}--\eqref{zz-12} and solving
\begin{align*}
&W_\varphi+b^{1/2}(H_\varepsilon(Q))W_\psi
=\frac14b^{-1}(H_\varepsilon(Q))p(H_\varepsilon(Q))W(W-{H'_\varepsilon(Q)}{Q_\varphi}
&&
\\
&\qquad\qquad\qquad\qquad\qquad\qquad
-b^{1/2}(H_\varepsilon(Q))
{H'_\varepsilon(Q)}{Q_\psi}),
\quad&&(\varphi,\psi)\in(\hat\zeta_0,\hat\zeta_1)\times(0,m),
\\
&Z_\varphi-b^{1/2}(H_\varepsilon(Q))Z_\psi
=-\frac14b^{-1}(H_\varepsilon(Q))p(H_\varepsilon(Q))Z({H'_\varepsilon(Q)}{Q_\varphi}
&&
\\
&\qquad\qquad\qquad\qquad\qquad\qquad
-b^{1/2}(H_\varepsilon(Q))
{H'_\varepsilon(Q)}{Q_\psi}+Z),
\quad&&(\varphi,\psi)\in(\hat\zeta_0,\hat\zeta_1)\times(0,m),
\\
&W(\hat\zeta_0,\psi)=W_1(\psi),\quad&&\psi\in(0,m),
\\
&Z(\hat\zeta_0,\psi)=Z_1(\psi),\quad&&\psi\in(0,m),
\\
&W(\varphi,0)+Z(\varphi,0)=0,\quad&&\varphi\in(\hat\zeta_0,\hat\zeta_1),
\\
&W(\varphi,m)+Z(\varphi,m)=\frac{2f''(X_{\text{\rm up}}(\varphi))(1+(f'(X_{\text{\rm up}}(\varphi)))^2)^{-3/2}}
{b^{1/2}(H_\varepsilon(Q(\varphi,m)))A^{-1}_+(Q(\varphi,m))},
\quad&&\varphi\in(\hat\zeta_0,\hat\zeta_1),
\\
&Q_\varphi=\frac{1}{2}(W-Z),
\quad&&(\varphi,\psi)\in(\hat\zeta_0,\hat\zeta_1)\times(0,m),
\\
&Q(\hat\zeta_0,\psi)=Q_1(\psi),\quad&&\psi\in(0,m),
\\
&{\mathscr Q}(x)=A^{-1}_+(Q(\Phi_{\text{\rm up}}(x),m)),
\quad&&x\in(X_{\text{\rm up}}(\hat\zeta_0),X_{\text{\rm up}}(\hat\zeta_1)).
\end{align*}
Let
$$
\hat\zeta_{1,\varepsilon}=\sup\Big\{\tilde\varphi\in(\hat\zeta_0,\hat\zeta_1):
Q(\varphi,\psi)>-\frac1\varepsilon
\mbox{ for each }(\varphi,\psi)\in(\hat\zeta_0,\tilde\varphi)\times(0,m)\Big\}.
$$
Then, $(W,Z,Q)$ solves
\begin{align*}
&W_\varphi+b^{1/2}(Q)W_\psi
=\frac14b^{-1}(Q)p(Q)W(W+Z),
\quad&&(\varphi,\psi)\in(\hat\zeta_0,\hat\zeta_{1,\varepsilon})\times(0,m),
\\
&Z_\varphi-b^{1/2}(Q)Z_\psi
=-\frac14b^{-1}(Q)p(Q)Z(W+Z),
\quad&&(\varphi,\psi)\in(\hat\zeta_0,\hat\zeta_{1,\varepsilon})\times(0,m),
\\
&W(\hat\zeta_0,\psi)=W_1(\psi),\quad&&\psi\in(0,m),
\\
&Z(\hat\zeta_0,\psi)=Z_1(\psi),\quad&&\psi\in(0,m),
\\
&W(\varphi,0)+Z(\varphi,0)=0,\quad&&\varphi\in(\hat\zeta_0,\hat\zeta_{1,\varepsilon}),
\\
&W(\varphi,m)+Z(\varphi,m)=\frac{2f''(X_{\text{\rm up}}(\varphi))(1+(f'(X_{\text{\rm up}}(\varphi)))^2)^{-3/2}}
{b^{1/2}(Q(\varphi,m))A^{-1}_+(Q(\varphi,m))},
\quad&&\varphi\in(\hat\zeta_0,\hat\zeta_{1,\varepsilon}),
\\
&Q_\varphi=\frac{1}{2}(W-Z),
\quad&&(\varphi,\psi)\in(\hat\zeta_0,\hat\zeta_{1,\varepsilon})\times(0,m),
\\
&Q(\hat\zeta_0,\psi)=Q_1(\psi),\quad&&\psi\in(0,m),
\\
&{\mathscr Q}(x)=A^{-1}_+(Q(\Phi_{\text{\rm up}}(x),m)),
\quad&&x\in(X_{\text{\rm up}}(\hat\zeta_0),X_{\text{\rm up}}(\hat\zeta_{1,\varepsilon})).
\end{align*}
Combining the above existence result and Step I shows that
$(W,Z,Q)\in L^{\infty}((0,\hat\zeta_{1,\varepsilon})\times(0,m))\times
L^{\infty}((0,\hat\zeta_{1,\varepsilon})\times(0,m))
\times C^{0,1}([0,\hat\zeta_{1,\varepsilon}]\times[0,m])$, which satisfies
\begin{gather}
\label{zz-16}
-\frac1\varepsilon<Q(\varphi,\psi)\le-\mu_1,\quad
-Q_\varphi(\varphi,\psi)\pm b^{1/2}(Q(\varphi,\psi))Q_\psi(\varphi,\psi)\ge0,\quad
(\varphi,\psi)\in(0,\hat\zeta_{1,\varepsilon})\times(0,m)
\end{gather}
and solves
\begin{align}
\label{zz-16-1}
&W_\varphi+b^{1/2}(Q)W_\psi
=\frac14b^{-1}(Q)p(Q)W(W+Z),
\quad&&(\varphi,\psi)\in(0,\hat\zeta_{1,\varepsilon})\times(0,m),
\\
&Z_\varphi-b^{1/2}(Q)Z_\psi
=-\frac14b^{-1}(Q)p(Q)Z(W+Z),
\quad&&(\varphi,\psi)\in(0,\hat\zeta_{1,\varepsilon})\times(0,m),
\\
&W(0,\psi)=W_0(\psi),\quad&&\psi\in(0,m),
\\
&Z(0,\psi)=Z_0(\psi),\quad&&\psi\in(0,m),
\\
&W(\varphi,0)+Z(\varphi,0)=0,\quad&&\varphi\in(0,\hat\zeta_{1,\varepsilon}),
\\
\label{zz-16-2}
&W(\varphi,m)+Z(\varphi,m)=\frac{2f''(X_{\text{\rm up}}(\varphi))(1+(f'(X_{\text{\rm up}}(\varphi)))^2)^{-3/2}}
{b^{1/2}(Q(\varphi,m))A^{-1}_+(Q(\varphi,m))},
\quad&&\varphi\in(0,\hat\zeta_{1,\varepsilon}),
\\
\label{zz-16-3}
&Q_\varphi=\frac{1}{2}(W-Z),
\quad&&(\varphi,\psi)\in(0,\hat\zeta_{1,\varepsilon})\times(0,m),
\\
\label{zz-16-4}
&Q(0,\psi)=Q_0(\psi),\quad&&\psi\in(0,m),
\\
&{\mathscr Q}(x)=A^{-1}_+(Q(\Phi_{\text{\rm up}}(x),m)),
\quad&&x\in(l_0,X_{\text{\rm up}}(\hat\zeta_{1,\varepsilon})).
\nonumber
\end{align}

{\bf Step III}\quad Estimate.

Note that $(W,Z)$ solves the problem \eqref{zz-16-1}--\eqref{zz-16-2} with
$Q$ satisfying \eqref{zz-16}.
Then, Lemma \ref{sss-prop1} shows that $W$ and $Z$ satisfy
\begin{gather}
\label{zz-13}
-\mu_4(1+\mu_3)\le W(\varphi,\psi)\le 0,\quad
0\le Z(\varphi,\psi)\le\mu_4(1+\mu_3),\quad
(\varphi,\psi)\in(0,\hat\zeta_{1,\varepsilon})\times(0,m),
\end{gather}
which, together with \eqref{zz-16-3}, \eqref{zz-16-4} and \eqref{sss-1}, also leads to
\begin{gather}
\label{zz-14}
-\mu_2-\mu_4(1+\mu_3)\hat\zeta\le Q(\varphi,\psi)\le-\mu_1,\quad
(\varphi,\psi)\in(0,\hat\zeta_{1,\varepsilon})\times(0,m).
\end{gather}

{\bf Step IV}\quad Iteration.

If $\hat\zeta_{1,\varepsilon}=\hat\zeta$ or
$\min_{[0,m]}Q(\hat\zeta_{1,\varepsilon},\cdot)=-1/\varepsilon$,
then the claim is proved.
Otherwise, $\hat\zeta_{1,\varepsilon}<\hat\zeta$ and
$\min_{[0,m]}Q(\hat\zeta_{1,\varepsilon},\cdot)>-1/\varepsilon$,
which also leads to
$\hat\zeta_{1,\varepsilon}=\hat\zeta_1$.
Set
\begin{align*}
W_2(\psi)=W(\hat\zeta_{1},\psi),\quad
Z_2(\psi)=Z(\hat\zeta_{1},\psi),\quad
Q_2(\psi)=Q(\hat\zeta_{1},\psi),\quad\psi\in(0,m).
\end{align*}
It follows from \eqref{zz-13} and \eqref{zz-14} that
$Q_2\in C^{0,1}([0,m])$ and $W_2,Z_2\in L^{\infty}(0,m)$
satisfy
\begin{align*}
-\hat\mu_2\le Q_2(\psi)\le-\mu_1,\quad
-\hat\mu_3\le W_2(\psi)\le0,\quad
0\le Z_2(\psi)\le\hat\mu_3,\quad\psi\in(0,m).
\end{align*}
Note that $(W_2,Z_2,Q_2)$ satisfies the same condition as \eqref{zz-8-1}
for $(W_1,Z_1,Q_1)$.
Then, repeating Step II, one gets $\hat\zeta_{2,\varepsilon}\in(\zeta_1,\hat\zeta]$ such that
the problem \eqref{newsupp-neq1}--\eqref{newsupp-q4}
with $\zeta=\hat\zeta_{2,\varepsilon}$ admits a weak solution
$(W,Z,Q)\in L^{\infty}((0,\hat\zeta_{2,\varepsilon})\times(0,m))\times
L^{\infty}((0,\hat\zeta_{2,\varepsilon})\times(0,m))
\times C^{0,1}([0,\hat\zeta_{2,\varepsilon}]\times[0,m])$
satisfying
\begin{align*}
-\frac1\varepsilon<Q(\varphi,\psi)\le-\mu_1,\quad
-Q_\varphi(\varphi,\psi)\pm b^{1/2}(Q(\varphi,\psi))Q_\psi(\varphi,\psi)\ge0,\quad
(\varphi,\psi)\in(0,\hat\zeta_{2,\varepsilon})\times(0,m)
\end{align*}
and either
$\hat\zeta_{2,\varepsilon}=\min\big\{\hat\zeta,\hat\zeta_0+{2}/{\hat\beta}\big\}$
or
$\hat\zeta_{2,\varepsilon}<\min\big\{\hat\zeta,\hat\zeta_0+{2}/{\hat\beta}\big\}$
with $\min_{[0,m]}Q(\hat\zeta_{2,\varepsilon},\cdot)=-1/\varepsilon$.
Then, $\hat\zeta_{\varepsilon}$ in the claim can be obtained after repeating Step III and Step II
for a finite number of times.
$\hfill\Box$\vskip 4mm

\begin{remark}
The estimates in Lemmas \ref{sss-prop1} and \ref{sss-prop2}
are independent of the value of the solution in the domain.
So, the solution to the problem \eqref{b-1}--\eqref{b-9} can satisfy
the same estimates at each $\varphi$, which guarantees the success of the iteration
in the proof of Theorem \ref{sss-thm1}.
\end{remark}

\begin{remark}
The iteration in the proof of Theorem \ref{sss-thm1} begins from not $(W_1,Z_1,Q_1)$
but $(W_2,Z_2,Q_2)$.
Indeed, the conditions between $(W_1,Z_1,Q_1)$ and $(W_0,Z_0,Q_0)$ are different,
while $(W_2,Z_2,Q_2)$ satisfies the same conditions as for $(W_1,Z_1,Q_1)$.
Note that the estimates in Lemmas \ref{sss-prop1} and \ref{sss-prop2}
depend on the initial data.
To show that $(W_2,Z_2,Q_2)$ satisfies the same conditions as for $(W_1,Z_1,Q_1)$,
the initial data of $(W_2,Z_2,Q_2)$ is taken to be not $(W_1,Z_1,Q_1)$
but $(W_0,Z_0,Q_0)$ (see {Step III}).
\end{remark}

\begin{remark}
It should be noted that for the local existence of sonic-supersonic flows in \cite{WX3},
the iteration scheme in the fixed point argument is
\begin{align*}
&W_\varphi+b^{1/2}(\tilde Q)W_\psi
=\frac14b^{-1}(\tilde Q)p(\tilde Q)\tilde W(W+Z),
\\
&Z_\varphi-b^{1/2}(\tilde Q)Z_\psi
=-\frac14b^{-1}(\tilde Q)p(\tilde Q)\tilde Z(W+Z),
\end{align*}
where
$$
\tilde W=\tilde Q_\varphi-b^{1/2}(\tilde Q)\tilde Q_\psi,\quad
\tilde Z=-\tilde Q_\varphi-b^{1/2}(\tilde Q)\tilde Q_\psi.
$$
This iteration does not work for Theorem \ref{sss-thm1}.
The reason lies in that
the desired estimates in \cite{WX3} need $\tilde W\le0$ and $\tilde Z\ge0$,
while the set of $\tilde Q$ satisfying $\tilde W\le0$ and $\tilde Z\ge0$
is not convex.
Furthermore, another disadvantage is that in Theorem \ref{sss-thm1}
the boundary condition at $\psi=m$ may be large.
\end{remark}

In the physical plane, Theorem \ref{sss-thm1} can be stated as follows.

\begin{theorem}
\label{phys-thm31}
Assume that $f\in C^{2}([l_0,l_1))$ satisfies \eqref{sss-0},
$\Upsilon\in C^{2}([0,f(l_0)])$ satisfies \eqref{xc4} and $q_0\in C^{1}([0,f(l_0)])$
satisfies \eqref{phy-8-q}--\eqref{phy-9}.
Then the problem \eqref{phy1}--\eqref{phy5} admits uniquely a maximal smooth supersonic flow before
vacuum formation.
Furthermore, the solution satisfies one and only one of the following two cases:

{\rm (i)} $\varphi\in C^{2}(\overline\Omega)$ and
$c_1\le|\nabla\varphi|<c^*$ on $\overline\Omega$.

{\rm (ii)} $\varphi\in C^{2}(\overline\Omega_{\zeta}\setminus\Gamma_{\zeta})$ with
\begin{align*}
c_1\le |\nabla\varphi(x,y)|<c^*\mbox{ for each }
(x,y)\in\overline\Omega_{\zeta}\setminus\Gamma_{\zeta},
\quad\sup_{\Omega_\zeta}|\nabla\varphi|=c^*
\end{align*}
for some constant $\zeta>0$, where
$\Omega_{\zeta}=\big\{(x,y)\in\Omega:\varphi(x,y)<\zeta\big\}$
and $\Gamma_{\zeta}=\overline{\partial\Omega_{\zeta}\cap\Omega}$.
\end{theorem}

\subsection{Global smooth supersonic flows without vacuum in straight nozzles}

According to Theorems \ref{phys-thm31} and \ref{ni2-thm},
under the assumptions of Theorem \ref{phys-thm31},
the supersonic flow problem \eqref{phy1}--\eqref{phy5} admits
a smooth solution before vacuum formation,
furthermore, if in addition, the nozzle is straight,
then the smooth supersonic flow is away from vacuum in any bounded region.
Therefore, under the assumptions of Theorem \ref{phys-thm31},
there always exists a global smooth supersonic flow without sonic and vacuum point in a straight nozzle.

\begin{theorem}
\label{phys-thm51}
Assume that $f\in C^2([l_0,+\infty))$ satisfies \eqref{aaasss-0},
$\Upsilon\in C^{2}([0,f(l_0)])$ satisfies \eqref{xc4} and $q_0\in C^{1}([0,f(l_0)])$
satisfies \eqref{phy-8-q}--\eqref{phy-9}.
Then the problem \eqref{phy1}--\eqref{phy5} admits uniquely a smooth supersonic flow,
which is away from vacuum in any bounded region.
Furthermore, the solution $\varphi\in C^{2}(\overline\Omega)$ and
$c_1\le|\nabla\varphi|<c^*$ on $\overline\Omega$.
\end{theorem}

Below we indicate how to solve the Cauchy problem \eqref{cv1-a1}--\eqref{cv1-a3}.

\begin{theorem}
\label{ni3}
Assume that $Q_0\in C^{1}(\mathbb R)$ and $G_0\in C(\mathbb R)$
satisfy
\begin{align*}
-\mu_2\le Q_0(\psi)\le-\mu_1,\quad
0\le-G_0(\psi)\pm b^{1/2}(Q_0(\psi))Q'_0(\psi)\le\mu_3,\quad\psi\in\mathbb R
\end{align*}
with three constants $0<\mu_1\le\mu_2$ and $\mu_3>0$.
Then the Cauchy problem \eqref{cv1-a1}--\eqref{cv1-a3}
admits uniquely a solution
$Q\in C^{1}([0,+\infty)\times\mathbb R)$.
Furthermore, the solution satisfies
\begin{align}
\label{zzc1}
Q(\varphi,\psi)\le-\mu_1,\quad
0\le-Q_\varphi(\varphi,\psi)\pm b^{1/2}(Q(\varphi,\psi))Q_\varphi(\varphi,\psi)\le\mu_3,\quad
(\varphi,\psi)\in(0,+\infty)\times\mathbb R.
\end{align}
\end{theorem}

\Proof
The problem \eqref{cv1-a1}--\eqref{cv1-a3} is equivalent to
\begin{align*}
&W_\varphi+b^{1/2}(Q)W_\psi
=\frac14b^{-1}(Q)p(Q) W(W+Z),
\quad&&(\varphi,\psi)\in(0,+\infty)\times\mathbb R,
\\
&Z_\varphi-b^{1/2}(Q)Z_\psi
=-\frac14b^{-1}(Q)p(Q) Z(W+Z),
\quad&&(\varphi,\psi)\in(0,+\infty)\times\mathbb R,
\\
&W(0,\psi)=W_0(\psi)=G_0(\psi)-b^{1/2}(Q_0(\psi))Q'_0(\psi),\quad&&\psi\in\mathbb R,
\\
&Z(0,\psi)=Z_0(\psi)=-G_0(\psi)-b^{1/2}(Q_0(\psi))Q'_0(\psi),\quad&&\psi\in\mathbb R,
\\
&Q_\varphi=\frac{1}{2}(W-Z),
\quad&&(\varphi,\psi)\in(0,+\infty)\times\mathbb R,
\\
&Q(0,\psi)=Q_0(\psi),\quad&&\psi\in\mathbb R.
\end{align*}
Similar to the proof of Theorem \ref{sss-thm1}, one can prove that
the problem \eqref{cv1-a1}--\eqref{cv1-a3} admits uniquely a solution
$Q\in C^{1}([0,\zeta)\times\mathbb R)$,
where either
$\zeta=+\infty$
or
$0<\zeta<+\infty$ with $\lim_{\varphi\to\zeta^-}\inf_{\mathbb R}Q(\varphi,\cdot)=-\infty$.
Furthermore, the solution satisfies
\begin{align*}
-Q_\varphi(\varphi,\psi)\pm b^{1/2}(Q(\varphi,\psi))Q_\psi(\varphi,\psi)\ge0,\quad
(\varphi,\psi)\in(0,\zeta)\times\mathbb R.
\end{align*}
Proposition \ref{ni2} shows that
\begin{align*}
-Q_\varphi(\varphi,\psi)\pm b^{1/2}(Q(\varphi,\psi))Q_\psi(\varphi,\psi)\le\mu_3,\quad
(\varphi,\psi)\in(0,\zeta)\times\mathbb R.
\end{align*}
Therefore, $\zeta=+\infty$ and $Q\in C^{1}([0,+\infty)\times\mathbb R)$ satisfies \eqref{zzc1}.
By the way, one can use another iteration scheme for the existence,
where the cut-off technique is not needed.
Fix $\tilde\zeta>0$ and define
\begin{align*}
{\mathscr B}=&\Big\{Q\in C^{0,1}([0,\tilde\zeta]\times\mathbb R):
-\mu_2\le Q(0,\psi)\le-\mu_1\mbox{ for }\psi\in\mathbb R,
-\mu_3\le Q_\varphi(\varphi,\psi)\le0
\\
&\qquad\mbox{ and }
|Q_\psi(\varphi,\psi)|\le\beta\mbox{ for }
(\varphi,\psi)\in(0,\tilde\zeta)\times\mathbb R\Big\},
\end{align*}
where $\beta>0$ is a constant to be determined.
For fixed $\tilde Q\in {\mathscr B}$, consider the problem
\begin{align}
\label{aaab-1}
&W_\varphi+b^{1/2}(\tilde Q)W_\psi
=\frac14b^{-1}(\tilde Q)p(\tilde Q) W(W+Z),
\quad&&(\varphi,\psi)\in(0,\tilde\zeta)\times\mathbb R,
\\
\label{aaab-2}
&Z_\varphi-b^{1/2}(\tilde Q)Z_\psi
=-\frac14b^{-1}(\tilde Q)p(\tilde Q) Z(W+Z),
\quad&&(\varphi,\psi)\in(0,\tilde\zeta)\times\mathbb R,
\\
\label{aaab-3}
&W(0,\psi)=W_0(\psi),\quad&&\psi\in\mathbb R,
\\
\label{aaab-4}
&Z(0,\psi)=Z_0(\psi),\quad&&\psi\in\mathbb R,
\\
\label{aaab-5}
&Q_\varphi=\frac{1}{2}(W-Z),
\quad&&(\varphi,\psi)\in(0,\tilde\zeta)\times\mathbb R,
\\
\label{aaab-6}
&Q(0,\psi)=Q_0(\psi),\quad&&\psi\in\mathbb R.
\end{align}
The key estimate, which can be proved by the method of characteristics, is that for each $(w,z)\in L^\infty((0,\tilde\zeta)\times\mathbb R)
\times L^\infty((0,\tilde\zeta)\times\mathbb R)$ with
\begin{align*}
-\mu_3\le w(\varphi,\psi)\le0,\quad
0\le z(\varphi,\psi)\le\mu_3,\quad
(\varphi,\psi)\in(0,\tilde\zeta)\times\mathbb R,
\end{align*}
the problem
\begin{align*}
&W_\varphi+b^{1/2}(\tilde Q)W_\psi
=\frac14b^{-1}(\tilde Q)p(\tilde Q) W(W+z),
\quad&&(\varphi,\psi)\in(0,\tilde\zeta)\times\mathbb R,
\\
&Z_\varphi-b^{1/2}(\tilde Q)Z_\psi
=-\frac14b^{-1}(\tilde Q)p(\tilde Q) Z(w+Z),
\quad&&(\varphi,\psi)\in(0,\tilde\zeta)\times\mathbb R,
\\
&W(0,\psi)=W_0(\psi),\quad&&\psi\in\mathbb R,
\\
&Z(0,\psi)=Z_0(\psi),\quad&&\psi\in\mathbb R
\end{align*}
admits uniquely a weak solution $(W,Z)\in L^\infty((0,\tilde\zeta)\times\mathbb R)
\times L^\infty((0,\tilde\zeta)\times\mathbb R)$,
and the solution satisfies
\begin{align*}
-\mu_3\le W(\varphi,\psi)\le0,\quad
0\le Z(\varphi,\psi)\le\mu_3,\quad
(\varphi,\psi)\in(0,\tilde\zeta)\times\mathbb R.
\end{align*}
Based on this key estimate, the problem \eqref{aaab-1}--\eqref{aaab-4}
admits a unique weak solution
$(W,Z)\in L^\infty((0,\tilde\zeta)\times\mathbb R)
\times L^\infty((0,\tilde\zeta)\times\mathbb R)$ satisfying
\begin{align*}
-\mu_3\le W(\varphi,\psi)\le0,\quad
0\le Z(\varphi,\psi)\le\mu_3,\quad
(\varphi,\psi)\in(0,\tilde\zeta)\times\mathbb R.
\end{align*}
Then, similar to the proof of Lemma \ref{sss-prop2}, one can show that
the problem \eqref{aaab-5}, \eqref{aaab-6} admits a unique weak solution $Q\in C^{0,1}([0,\tilde\zeta]\times\mathbb R)$
and it satisfies
$$
-\mu_3\le Q_\varphi(\varphi,\psi)\le0,\quad
|Q_\psi(\varphi,\psi)|\le \beta,\quad
(\varphi,\psi)\in(0,\tilde\zeta)\times\mathbb R
$$
with a positive constant $\beta=\beta(\gamma,\mu_1,\mu_2,\mu_3,\tilde\zeta)$.
$\hfill\Box$\vskip 4mm

\begin{remark}
The key estimate in Theorem \ref{ni3} is invalid for the problem \eqref{newsupp-eq}--\eqref{newsupp-q}
with a nonhomogeneous boundary condition at $\psi=m$.
So, the new iteration scheme in Theorem \ref{ni3} is unsuitable for Theorem \ref{sss-thm1}.
\end{remark}

\subsection{Properties of smooth supersonic flows before vacuum formation}

In this subsection,
we investigate the properties of smooth supersonic flows before vacuum formation,
including the lower bound estimate of the flow speed in the potential plane
and the Lipschitz continuity of the flow speed in the physical plane.

\begin{proposition}
\label{sss-thm1-lb}
Let $f\in C^{2}([l_0,l_1))$ satisfy \eqref{sss-0},
and $Q_0\in C^{1}([0,m])$ and $G_0\in C([0,m])$
satisfy \eqref{sss-1}, \eqref{sss-2} and
\begin{align}
\label{sss-1-lb}
-G_0(\psi)\pm b^{1/2}(Q_0(\psi))Q'_0(\psi)\ge\mu_0,\quad\psi\in(0,m),
\end{align}
where $\mu_0>0$ is a positive constant.
Assume that $Q\in C^{1}([0,\zeta)\times[0,m])$ is a solution
to the problem \eqref{newsupp-eq}--\eqref{newsupp-q} with some $0<\zeta\le+\infty$.
Then there exists a positive constant $M=M(\gamma,m,\mu_0,\mu_1)$
such that
\begin{gather}
\label{sd-thm3-re0}
Q(\varphi,\psi)\le-M (\varphi+1)^{2/(\gamma+1)},\quad
(\varphi,\psi)\in(0,\zeta)\times(0,m).
\end{gather}
If, in addition, $f''(1+(f')^2)^{-3/2}\in L^\infty(l_0,l_1)$,
then
\begin{gather}
\label{sd-thm3-re1}
-N\le b^{1/2}(Q(\varphi,\psi))Q_\varphi(\varphi,\psi)\le 0,\quad
|b(Q(\varphi,\psi))Q_\psi(\varphi,\psi)|\le N,\quad
(\varphi,\psi)\in(0,\zeta)\times(0,m),
\end{gather}
where $N=N(\gamma,m,\mu_0,\mu_1,\mu_3,\|f''(1+(f')^2)^{-3/2}\|_{L^\infty(l_0,l_1)})$ is a positive constant.
\end{proposition}

\Proof
For $(\varphi,\psi)\in[0,\zeta)\times[0,m]$, set
\begin{align*}
W(\varphi,\psi)=Q_\varphi(\varphi,\psi)-{b^{1/2}(Q(\varphi,\psi))}Q_\psi(\varphi,\psi),\quad
Z(\varphi,\psi)=-Q_\varphi(\varphi,\psi)-{b^{1/2}(Q(\varphi,\psi))}Q_\psi(\varphi,\psi).
\end{align*}
Theorem \ref{sss-thm1} shows that
\begin{gather}
\label{lb-0}
Q(\varphi,\psi)\le-\mu_1,\quad(\varphi,\psi)\in(0,\zeta)\times(0,m)
\end{gather}
and
$(W,Z)\in C([0,\zeta)\times[0,m])\times C([0,\zeta)\times[0,m])$
solves
\begin{align}
\nonumber
&W_\varphi+b^{1/2}(Q)W_\psi
=\frac14b^{-1}(Q)p(Q)W(W+Z),
\quad&&(\varphi,\psi)\in(0,\zeta)\times(0,m),
\\
\nonumber
&Z_\varphi-b^{1/2}(Q)Z_\psi
=-\frac14b^{-1}(Q)p(Q)Z(W+Z),
\quad&&(\varphi,\psi)\in(0,\zeta)\times(0,m),
\\
\label{lbnewsupp-ninbc1}
&W(0,\psi)=W_0(\psi)=G_0(\psi)-b^{1/2}(Q_0(\psi))Q'_0(\psi)\le-\mu_0,\quad&&\psi\in(0,m),
\\
\label{lbnewsupp-ninbc2}
&Z(0,\psi)=Z_0(\psi)=-G_0(\psi)-b^{1/2}(Q_0(\psi))Q'_0(\psi)\ge\mu_0,\quad&&\psi\in(0,m),
\\
\label{lbnewsupp-nlbbc}
&W(\varphi,0)+Z(\varphi,0)=0,\quad&&\varphi\in(0,\zeta),
\\
\label{lbnewsupp-nubbc}
&W(\varphi,m)+Z(\varphi,m)=\frac{2f''(X_{\text{\rm up}}(\varphi))(1+(f'(X_{\text{\rm up}}(\varphi)))^2)^{-3/2}}
{b^{1/2}(Q(\varphi,m))A^{-1}_+(Q(\varphi,m))}\ge0,
\quad&&\varphi\in(0,\zeta).
\end{align}
By the method of characteristics, it is clear that
\begin{align}
\label{lb-1}
W(\varphi,\psi)<0,\quad
Z(\varphi,\psi)>0,\quad
(\varphi,\psi)\in(0,\zeta)\times(0,m).
\end{align}

Assume that
$$
\Sigma_{\pm}:
\Psi'_{\pm}(\varphi)=\pm b^{1/2}(Q(\varphi,\Psi_{\pm}(\varphi))),
\quad
0<\Psi_{\pm}(\varphi)<m,\quad\hat\varphi_{\pm}<\varphi<\check\varphi_{\pm}
\quad(0\le\hat\varphi_{\pm}<\check\varphi_{\pm}<\zeta)
$$
are a positive and a negative characteristics.
It follows from \eqref{lb-0} and \eqref{lb-1} that $W$ and $Z$ satisfy
\begin{align*}
\frac{d}{d\varphi}W(\varphi,\Psi_{+}(\varphi))\le
-\frac{M_1}{Q(\varphi,\Psi_{+}(\varphi))}W^2(\varphi,\Psi_{+}(\varphi)),
\quad\hat\varphi_+<\varphi<\check\varphi_+
\end{align*}
and
\begin{align*}
\frac{d}{d\varphi}Z(\varphi,\Psi_{-}(\varphi))\ge
\frac{M_1}{Q(\varphi,\Psi_{+}(\varphi))}Z^2(\varphi,\Psi_{-}(\varphi)),
\quad\hat\varphi_-<\varphi<\check\varphi_-
\end{align*}
on $\Sigma_{\pm}$, respectively, where $M_1=M_1(\gamma,\mu_1)$ is a positive constant.
Therefore,
\begin{align}
\label{lb-2}
&\frac{d}{d\varphi}\Big(\frac1{W(\varphi,\Psi_{+}(\varphi))}\Big)\ge
\frac{M_1}{Q(\varphi,\Psi_{+}(\varphi))},
&&\quad\hat\varphi_+<\varphi<\check\varphi_+,
\\
\label{lb-3}
&\frac{d}{d\varphi}\Big(\frac1{Z(\varphi,\Psi_{-}(\varphi))}\Big)\le
-\frac{M_1}{Q(\varphi,\Psi_{-}(\varphi))},
&&\quad\hat\varphi_-<\varphi<\check\varphi_-,
\end{align}
which, together with \eqref{lb-0}, yield
\begin{align}
\label{lb-4}
&\frac{d}{d\varphi}\Big(\frac1{W(\varphi,\Psi_{+}(\varphi))}\Big)\ge
-\frac{M_1}{\mu_1},
&&\quad\hat\varphi_+<\varphi<\check\varphi_+,
\\
\label{lb-5}
&\frac{d}{d\varphi}\Big(\frac1{Z(\varphi,\Psi_{-}(\varphi))}\Big)\le
\frac{M_1}{\mu_1},
&&\quad\hat\varphi_-<\varphi<\check\varphi_-.
\end{align}
Direct calculation shows
\begin{align*}
&\frac{d}{d\varphi}\Big(\frac1{(-Q(\varphi,\Psi_{+}(\varphi)))^{(\gamma-1)/2}W(\varphi,\Psi_{+}(\varphi))}\Big)
\\
=&\frac1{(-Q(\varphi,\Psi_{+}(\varphi)))^{(\gamma-1)/2}}
\frac{d}{d\varphi}\Big(\frac1{W(\varphi,\Psi_{+}(\varphi))}\Big)
\\
&\qquad+\frac{\gamma-1}{2}\,
\frac{Q_\varphi(\varphi,\Psi_{+}(\varphi))+b^{1/2}(\varphi)Q_\psi(\varphi,\Psi_{+}(\varphi))}
{(-Q(\varphi,\Psi_{+}(\varphi)))^{(\gamma-1)/2+1}W(\varphi,\Psi_{+}(\varphi))}
\\
=&\frac1{(-Q(\varphi,\Psi_{+}(\varphi)))^{(\gamma-1)/2}}\frac{d}{d\varphi}\Big(\frac1{W(\varphi,\Psi_{+}(\varphi))}\Big)
\\
&\qquad-\frac{\gamma-1}{2}\,
\frac{Z(\varphi,\Psi_{+}(\varphi))}
{(-Q(\varphi,\Psi_{+}(\varphi)))^{(\gamma-1)/2+1}W(\varphi,\Psi_{+}(\varphi))},
\quad\hat\varphi_+<\varphi<\check\varphi_+,
\end{align*}
which, together with \eqref{lb-1}, \eqref{lb-2} and \eqref{lb-0}, leads to
\begin{align}
\label{lb-4-0}
\frac{d}{d\varphi}\Big(\frac1{(-Q(\varphi,\Psi_{+}(\varphi)))^{(\gamma-1)/2}W(\varphi,\Psi_{+}(\varphi))}\Big)
\ge&-\frac{M_1}{(-Q(\varphi,\Psi_{+}(\varphi)))^{(\gamma+1)/2}}
\nonumber
\\
\ge&-{M_2}b^{1/2}(Q(\varphi,\Psi_{+}(\varphi))),
\quad\hat\varphi_+<\varphi<\check\varphi_+,
\end{align}
where $M_2=M_2(\gamma,\mu_1)$ is a positive constant.
Similarly, one gets from \eqref{lb-1}, \eqref{lb-3} and \eqref{lb-0} that
\begin{align}
\label{lb-5-0}
\frac{d}{d\varphi}\Big(\frac1{(-Q(\varphi,\Psi_{-}(\varphi)))^{(\gamma-1)/2}Z(\varphi,\Psi_{-}(\varphi))}\Big)
\le&{M_2}b^{1/2}(Q(\varphi,\Psi_{-}(\varphi))),
\quad\hat\varphi_-<\varphi<\check\varphi_-.
\end{align}

Fix $(\varphi_0,\psi_0)\in(0,\zeta)\times(0,m)$.
Similar to the proof of Lemma \ref{sss-prop1},
there exists two nonnegative integers $k_\pm$ such that
$$
\varphi_0=\varphi_0^\pm>\varphi_1^\pm>\cdots>\varphi_{k_\pm}^\pm>\varphi_{k_\pm+1}^\pm=0,
\quad\psi_{0}^\pm=\psi_0,
\quad0\le\psi_{k_\pm+1}^\pm\le m,
$$
$$
\psi_j^+=\left\{
\begin{aligned}
&0,\quad1\le j\le k_+\mbox{ and $j$ is odd},
\\
&m,\quad1\le j\le k_+\mbox{ and $j$ is even},
\end{aligned}
\right.
\quad\psi_j^-=\left\{
\begin{aligned}
&m,\quad1\le j\le k_-\mbox{ and $j$ is odd},
\\
&0,\quad1\le j\le k_-\mbox{ and $j$ is even},
\end{aligned}
\right.
$$
$$
\left\{
\begin{aligned}
&\frac{d}{d\varphi}\Psi_{j}^\pm(\varphi)=\pm(-1)^{j-1}b^{1/2}(Q(\varphi,\Psi_{j}^\pm(\varphi))),
\quad\varphi_{j}^\pm<\varphi<\varphi_{j-1}^\pm,
\\
&\Psi_{j}^\pm(\varphi_j^\pm)=\psi_{j}^\pm,\quad \Psi_{j}^\pm(\varphi_{j-1}^\pm)=\psi_{j-1}^\pm,
\end{aligned}
\right.
\quad1\le j\le k_\pm+1.
$$
Below we estimate $W(\varphi_0,\psi_0)$ and $Z(\varphi_0,\psi_0)$ by three steps.

{\bf Step I}\quad Rough estimates for $W(\varphi_0,\psi_0)$ and $Z(\varphi_0,\psi_0)$.

For $1\le j\le k_++1$, \eqref{lb-4} and \eqref{lb-5} show
\begin{align*}
&\frac1{W(\varphi_{j-1}^+,\psi_{j-1}^+)}\ge\frac1{W(\varphi_{j}^+,\psi_{j}^+)}-
\frac{M_1}{\mu_1}(\varphi_{j-1}^+-\varphi_{j}^+),
&&\quad\mbox{ if $j$ is odd},
\\
&\frac1{Z(\varphi_{j-1}^+,\psi_{j-1}^+)}\le\frac1{Z(\varphi_{j}^+,\psi_{j}^+)}+
\frac{M_1}{\mu_1}(\varphi_{j-1}^+-\varphi_{j}^+),
&&\quad\mbox{ if $j$ is even},
\end{align*}
which, together with \eqref{lbnewsupp-ninbc1}--\eqref{lbnewsupp-nubbc},
lead to
\begin{align}
\label{lb-5-3}
\frac1{W(\varphi_{0},\psi_{0})}\ge-\frac{1}{\mu_0}-\frac{M_1}{\mu_1}\varphi_{0}.
\end{align}
Similarly, one can show that
\begin{align}
\label{lb-5-4}
\frac1{Z(\varphi_{0},\psi_{0})}\le\frac{1}{\mu_0}+\frac{M_1}{\mu_1}\varphi_{0}.
\end{align}

{\bf Step II}\quad Upper bounds for $k_\pm$.

By \eqref{lb-5-3}, \eqref{lb-5-4} and the arbitrariness of $(\varphi_0,\psi_0)\in(0,\zeta)\times(0,m)$,
one gets that
\begin{align*}
Q_\varphi(\varphi,\psi)\le-\frac{1}
{1/{\mu_0}+{M_1}\varphi/{\mu_1}},
\quad(\varphi,\psi)\in(0,\zeta)\times(0,m),
\end{align*}
which yields
\begin{align}
\label{zxx1}
Q(\varphi,\psi)\le-M_3(1+\ln(\varphi+1)),\quad(\varphi,\psi)\in(0,\zeta)\times(0,m),
\end{align}
where $M_3=M_3(\gamma,\mu_0,\mu_1)$ is a positive constant.
It follows from \eqref{sss-1-lb} that
$\theta_\psi(0,\cdot)\ge\mu_0$ in $(0,m)$,
which, together with $\theta(0,0)=0$, leads to
$f'(l_0)=\tan\theta(0,m)\ge\tan(\mu_0m)$.
Choose a positive constant $M_4=M_4(\gamma,m,\mu_0,\mu_1)\ge0$ such that
\begin{align}
\label{lb-6}
{\mathscr H}(-M_3(1+\ln(M_4+1)))>{\mathscr H}(-\infty)-\arctan f'(l_0).
\end{align}
Note that \eqref{sss-0} yields
\begin{align}
\label{zxx2}
\theta_\varphi(\varphi,m)=\frac{d}{d\varphi}\arctan f'(X_{\text{\rm up}}(\varphi))
=\frac{f''(X_{\text{\rm up}}(\varphi)X'_{\text{\rm up}}(\varphi)}{1+(f'(X_{\text{\rm up}}(\varphi))^2}\ge0,
\quad0\le\varphi<\zeta.
\end{align}
If $0<\zeta\le M_4$, then \eqref{lb-0} yields
$k_\pm\le b^{1/2}(-\mu_1){M_4}/{m}+1$.
In the other case that $M_4<\zeta\le+\infty$,
for each $\tilde\varphi\ge M_4$,
Lemma \ref{lemma-sc4}, \eqref{zxx1}--\eqref{zxx2} show that
the positive characteristic from $(\tilde\varphi,0)$
must not touch $(0,\zeta)\times\{m\}$
and the negative characteristic from $(\tilde\varphi,m)$
must not touch $(0,\zeta)\times\{0\}$.
Therefore, there exists a positive integer $K=K(\gamma,m,\mu_0,\mu_1)$ such that
in both cases,
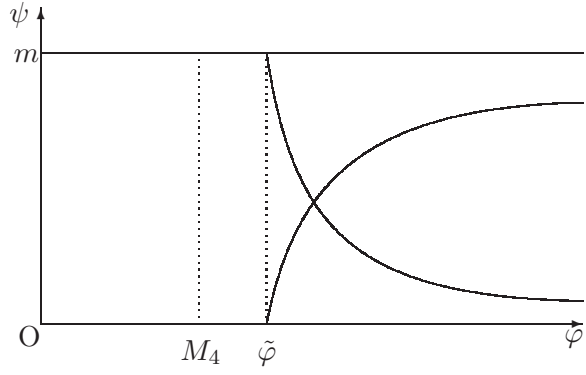
\begin{align}
\label{lb-7}
k_\pm\le K.
\end{align}

\vskip5mm

\hskip45mm
\setlength{\unitlength}{0.6mm}
\begin{picture}(250,65)
\put(0,0){\vector(1,0){120}}
\put(0,0){\vector(0,1){70}}
\put(116,-4){$\varphi$} \put(-7,67){$\psi$}
\put(0,60){\qbezier(0,0)(55,0)(120,0)}

\put(0,0){\cbezier(50,0)(56,30)(70,48)(120,49)}

\put(0,0){\cbezier(50,60)(56,20)(70,6)(120,5)}



\put(-6,58){$m$}

\put(-5,-5){O}

\put(48,-8){$\tilde\varphi$}

\put(50,0){\qbezier[40](0,0)(0,30)(0,60)}

\put(31,-8){$M_4$}

\put(35,0){\qbezier[40](0,0)(0,30)(0,60)}

\end{picture}

\vskip10mm

{\bf Step III}\quad Refined estimates for $W(\varphi_0,\psi_0)$ and $Z(\varphi_0,\psi_0)$.

For $1\le j\le k_++1$, \eqref{lb-4-0} and \eqref{lb-5-0} show
\begin{align*}
&\frac1{(-Q(\varphi_{j-1}^+,\psi_{j-1}^+))^{(\gamma-1)/2}W(\varphi_{j-1}^+,\psi_{j-1}^+)}\ge
\frac1{(-Q(\varphi_{j}^+,\psi_{j}^+))^{(\gamma-1)/2}W(\varphi_{j}^+,\psi_{j}^+)}-{M_2}m,
&&\quad\mbox{ if $j$ is odd},
\\
&\frac1{(-Q(\varphi_{j-1}^+,\psi_{j-1}^+))^{(\gamma-1)/2}Z(\varphi_{j-1}^+,\psi_{j-1}^+)}\le
\frac1{(-Q(\varphi_{j}^+,\psi_{j}^+))^{(\gamma-1)/2}Z(\varphi_{j}^+,\psi_{j}^+)}+{M_2}m,
&&\quad\mbox{ if $j$ is even},
\end{align*}
which, together with \eqref{lbnewsupp-ninbc1}--\eqref{lbnewsupp-nubbc} and
\eqref{lb-7}, yield
\begin{align}
\label{lb-10}
\frac1{(-Q(\varphi_{0},\psi_{0}))^{(\gamma-1)/2}W(\varphi_{0},\psi_{0})}
\ge-\frac{1}{\mu_1^{(\gamma-1)/2}\mu_0}-(K+1)M_2m.
\end{align}
Similarly, one can show that
\begin{align}
\label{lb-11}
\frac1{(-Q(\varphi_{0},\psi_{0}))^{(\gamma-1)/2}Z(\varphi_{0},\psi_{0})}
\le\frac{1}{\mu_1^{(\gamma-1)/2}\mu_0}+(K+1)M_2m.
\end{align}
It follows from \eqref{lb-10}, \eqref{lb-11} and
the arbitrariness of $(\varphi_0,\psi_0)\in(0,\zeta)\times(0,m)$ that
\begin{align*}
(-Q(\varphi,\psi))^{(\gamma-1)/2}Q_\varphi(\varphi,\psi)\le-\frac{1}
{{1}/{(\mu_1^{(\gamma-1)/2}\mu_0)}+(K+1)M_2m},
\quad(\varphi,\psi)\in(0,\zeta)\times(0,m),
\end{align*}
which yields \eqref{sd-thm3-re0} directly.

Below we prove \eqref{sd-thm3-re1}.
Set
\begin{align*}
&{\mathscr U}(\varphi,\psi)
=b^{1/2}(Q(\varphi,\psi))Q_\varphi(\varphi,\psi)-b(Q(\varphi,\psi))Q_\psi(\varphi,\psi),
\quad&&(\varphi,\psi)\in[0,\zeta)\times[0,m],
\\
&{\mathscr V}(\varphi,\psi)
=-b^{1/2}(Q(\varphi,\psi))Q_\varphi(\varphi,\psi)-b(Q(\varphi,\psi))Q_\psi(\varphi,\psi),
\quad&&(\varphi,\psi)\in[0,\zeta)\times[0,m].
\end{align*}
Theorem \ref{sss-thm1} shows
\begin{align}
\label{00sd2}
{\mathscr U}(\varphi,\psi)\le 0,\quad
{\mathscr V}(\varphi,\psi)\ge 0,\quad
(\varphi,\psi)\in(0,\zeta)\times(0,m).
\end{align}
It follows from \eqref{newsupp-neq1}--\eqref{newsupp-q4} and \eqref{00sd2} that
$({\mathscr U},{\mathscr V})\in C([0,\zeta)\times[0,m])
\times C([0,\zeta)\times[0,m])$ solves
\begin{align*}
&{\mathscr U}_\varphi+b^{1/2}(Q){\mathscr U}_\psi
=\frac14b^{-3/2}(Q)p(Q){\mathscr U}({\mathscr U}-{\mathscr V})\ge0,
\quad&&(\varphi,\psi)\in(0,\zeta)\times(0,m),
\\
&{\mathscr V}_\varphi-b^{1/2}(Q){\mathscr V}_\psi
=\frac14b^{-3/2}(Q)p(Q){\mathscr V}({\mathscr U}-{\mathscr V})\le0,
\quad&&(\varphi,\psi)\in(0,\zeta)\times(0,m),
\\
&{\mathscr U}(\varphi,0)+{\mathscr V}(\varphi,0)=0,\quad&&\varphi\in(0,\zeta),
\\
&{\mathscr U}(\varphi,m)+{\mathscr V}(\varphi,m)=
\frac{2f''(X_{\text{\rm up}}(\varphi))(1+(f'(X_{\text{\rm up}}(\varphi)))^2)^{-3/2}}
{A^{-1}_+(Q(\varphi,m))},\quad&&\varphi\in(0,\zeta).
\end{align*}
By $f''(1+(f')^2)^{-3/2}\in L^\infty(l_0,l_1)$
and \eqref{lb-7},
one can prove by the method of characteristics that
\begin{align*}
-N\le{\mathscr U}(\varphi,\psi)\le 0,\quad
0\le {\mathscr V}(\varphi,\psi)\le N,\quad
(\varphi,\psi)\in(0,\zeta)\times(0,m)
\end{align*}
with $N=N(\gamma,m,\mu_0,\mu_1,\mu_3,\|f''(1+(f')^2)^{-3/2}\|_{L^\infty(l_0,l_1)})$
being a positive constant,
which yields \eqref{sd-thm3-re1} directly.
$\hfill\Box$\vskip 4mm

\begin{remark}
If $q_0$ satisfies
\begin{align}
\label{phy-8-positive}
\inf_{(0,f(l_+))}\Big(
\frac{-\Upsilon''}{1+(\Upsilon')^2}\,\sqrt{\frac{-q_0^2\rho(q_0^2)}{\rho(q_0^2)+2q_0^2\rho'(q_0^2)}}
-|q'_0|\Big)>0,
\end{align}
then \eqref{sss-1-lb} holds.
\end{remark}

It can be verified easily that

\begin{remark}
\label{sd-thm3-re}
For the case $0<\zeta<+\infty$,
\eqref{sd-thm3-re1} remains true even if
$f''(1+(f')^2)^{-3/2}\not\in L^\infty(l_0,l_1)$ and \eqref{sss-1-lb} is not satisfied.
Here, the constant $N$
depends only on $\gamma$, $m$, $\mu_1$, $\mu_3$, $\zeta$ and $f$.
\end{remark}

As an application of Proposition \ref{sss-thm1-lb} in the physical plane, one gets

\begin{theorem}
\label{phys-thm53}
Assume that $f\in C^{2}([l_0,l_1))$ satisfies \eqref{sss-0}
and $f''(1+(f')^2)^{-3/2}\in L^\infty(l_0,l_1)$,
$\Upsilon\in C^{2}([0,f(l_0)])$ satisfies \eqref{xc4} and
$q_0\in C^{1}([0,f(l_0)])$
satisfies \eqref{phy-8-q}, \eqref{phy-8-positive} and \eqref{phy-9}.
Let $\varphi\in C^{2}(\overline\Omega)$
be a global supersonic flow without vacuum
to the problem \eqref{phy1}--\eqref{phy5}.
Then $|\nabla\varphi|$ is globally Lipschitz continuous in $\Omega$.
\end{theorem}

\begin{remark}
As will be shown in Theorem \ref{sss-thm3-4},
if $\lim_{\overline{x\to l_1^-}}f''(x)(1+(f'(x))^2)^{-3/2}>0$,
the smooth supersonic flow to the problem \eqref{phy1}--\eqref{phy5}
must admit a vacuum.
\end{remark}

\section{Formation of vacuum in a smooth supersonic flow}

For a smooth supersonic flow with vacuum to the problem \eqref{phy1}--\eqref{phy5},
it is clear that the set of vacuum points is closed
and thus there exists the first vacuum point in the increasing $x$-direction.
In this section, we will study the formation of the first vacuum point.

\subsection{Location of the first vacuum point and behavior of the flow near this point}

The following theorem shows that the first vacuum point of the flow
in the increasing $\varphi$-direction
must form at the upper wall.

\begin{theorem}
\label{sss-thm2}
Assume that $f\in C^{2}([l_0,l_1))$ satisfies \eqref{sss-0},
and $Q_0\in C^{1}([0,m])$ and $G_0\in C([0,m])$
satisfy \eqref{sss-1} and \eqref{sss-2}.
Let $Q\in C^{1}([0,\zeta)\times[0,m])$ with $0<\zeta<+\infty$
be the maximal smooth supersonic flow before vacuum formation
to the problem \eqref{newsupp-eq}--\eqref{newsupp-q}. Then
$q,\theta\in C([0,\zeta]\times[0,m])\cap C^{1}([0,\zeta]\times[0,m]\setminus\{(\zeta,m)\})$ satisfying
\begin{gather}
\label{sd4}
|q(\hat\varphi,\hat\psi)-q(\check\varphi,\check\psi)|\le
M\big(|\hat\varphi-\check\varphi|+|\hat\psi-\check\psi|^{1-1/\gamma}\big),
\quad
(\hat\varphi,\hat\psi),(\check\varphi,\check\psi)\in(0,\zeta)\times(0,m),
\\
\label{sd4-1}
q(\zeta,m)=c^*,\quad q(\varphi,\psi)<c^*,\quad
(\varphi,\psi)\in[0,\zeta]\times[0,m]\setminus\{(\zeta,m)\},
\\
\label{sd5}
q(\varphi,m)\ge c^*-M(\zeta-\varphi)^2,\quad0<\varphi<\zeta,
\\
\label{sd5-1}
|\theta_\varphi(\varphi,\psi)|\le M,
\quad
0\le (c^*-q(\varphi,\psi))^{1/(\gamma-1)+1/2}\theta_\psi(\varphi,\psi)\le  M,\quad
(\varphi,\psi)\in[0,\zeta]\times[0,m]\setminus\{(\zeta,m)\},
\end{gather}
where $M=M(\gamma,m,\mu_1,\mu_3,\zeta,f)$ is a positive constant.
\end{theorem}

\Proof
It follows from Proposition \ref{sss-thm1-lb} and Remark \ref{sd-thm3-re} that
\begin{align}
\label{bbv6}
-M_1\le{\mathscr U}(\varphi,\psi)\le 0,\quad
0\le {\mathscr V}(\varphi,\psi)\le M_1,\quad
(\varphi,\psi)\in(0,\zeta)\times(0,m)
\end{align}
with $M_1=M_1(\gamma,m,\mu_1,\mu_3,\zeta,f)$ being a positive constant.
Theorem \ref{sss-thm1} shows
$c_*<A_+^{-1}(-\mu_1)\le q<c^*$ in $(0,\zeta)\times(0,m)$,
which, together with \eqref{bbv6}, yields that
for each $(\varphi,\psi)\in(0,\zeta)\times(0,m)$,
\begin{gather}
\label{sd6}
0\le (c^*-q(\varphi,\psi))^{-1/2}q_\varphi(\varphi,\psi)\le M_2,
\quad
\big|(c^*-q(\varphi,\psi))^{1/(\gamma-1)}q_\psi(\varphi,\psi)|\le M_2,
\\
\label{sd7}
|\theta_\varphi(\varphi,\psi)|\le M_2,
\quad
0\le (c^*-q(\varphi,\psi))^{1/(\gamma-1)+1/2}\theta_\psi(\varphi,\psi)\le M_2,
\end{gather}
where $M_2=M_2(\gamma,m,\mu_1,\mu_3,\zeta,f)$ is a positive constant.
It follows from \eqref{sd6} that $q\in C([0,\zeta]\times[0,m])$
and it satisfies \eqref{sd4}.
Theorem \ref{sss-thm1} yields
$0<b^{1/2}(Q)\le b^{1/2}(-\mu_1)$ in $(0,\zeta)\times(0,m)$.
Therefore, for each $\max\{0,\zeta-mb^{-1/2}(-\mu_1)\}<\tilde\varphi<\zeta$,
Lemma \ref{aaasss-thm2} shows that
$Q\in C^1(\overline{\tilde\omega})$
with
$$
\tilde\omega=\big\{(\varphi,\psi)\in(\tilde\varphi,\zeta)\times(0,m):\psi<m-b^{1/2}(-\mu_1)(\varphi-\tilde\varphi)\big\}.
$$
Hence $Q\in C^{1}([0,\zeta]\times[0,m]\setminus\{(\zeta,m)\})$,
which, together with $q\in C([0,\zeta]\times[0,m])$ and
$\lim_{\varphi\to\zeta^-}\min_{[0,m]}Q(\varphi,\cdot)=-\infty$, leads to
\eqref{sd4-1} and $q,\theta\in C^{1}([0,\zeta]\times[0,m]\setminus\{(\zeta,m)\})$.
Finally, \eqref{sd5} follows from the first formula in \eqref{sd4-1} and the first estimate in \eqref{sd6},
while $\theta\in C([0,\zeta]\times[0,m])$ and \eqref{sd5-1} follow from \eqref{sd7}.
$\hfill\Box$\vskip 4mm

Below we characterize the boundary conditions of a local smooth supersonic flow
at the potential level set where the first vacuum forms.
Assume that $Q\in C^{1}([0,\zeta)\times[0,m])$ with $0<\zeta<+\infty$
is a maximal smooth supersonic flow before vacuum formation
to the problem \eqref{newsupp-eq}--\eqref{newsupp-q},
which is $\varphi\in C^{2}(\overline\Omega_{\zeta}\setminus\Gamma_{\zeta})$
to the problem \eqref{phy1}--\eqref{phy5} in the physical plane.
Theorems \ref{sss-thm2} and \ref{sss-thm1} show that
the first vacuum point of the flow in the increasing $x$-direction must form at the upper wall,
which is $(X_{\text{\rm up}}(\zeta),f(X_{\text{\rm up}}(\zeta)))$.
For convenience, denote $x_0=X_{\text{\rm up}}(\zeta)$
and $y_0=f(X_{\text{\rm up}}(\zeta))$.
It follows from Theorems \ref{sss-thm1} and \ref{sss-thm2} that
$Q\in C^{1}([0,\zeta]\times[0,m]\setminus\{(\zeta,m)\})$
satisfies
\begin{align*}
-Q_\varphi(\varphi,\psi)\pm b^{1/2}(Q(\varphi,\psi))Q_\psi(\varphi,\psi)\ge0,
\quad (\varphi,\psi)\in[0,\zeta]\times[0,m]\setminus\{(\zeta,m)\},
\end{align*}
$q,\theta\in C([0,\zeta]\times[0,m])$
satisfy $q(\zeta,m)=c^*$ and $\theta(\zeta,m)=\arctan f'(x_0)>0$.
Set
\begin{align*}
\tilde Q_0(\psi)=Q(\zeta,\psi),\quad
\tilde G_0(\psi)=Q_\varphi(\zeta,\psi),
\quad \psi\in[0,m).
\end{align*}
Then, $\tilde Q_0\in C^{1}([0,m))$ and $\tilde G_0\in C([0,m))\cap L^1(0,m)$
satisfying
\begin{gather}
\label{gg-1}
\tilde Q_0(\psi)<0,\quad
-\tilde G_0(\psi)\pm b^{1/2}(\tilde Q_0(\psi))\tilde Q'_0(\psi)\ge0,\quad\psi\in[0,m),
\\
\label{gg-2}
\tilde Q'_0(0)=0,\quad \lim_{\psi\to m^-}\tilde Q_0(\psi)=-\infty,
\quad \int_0^m\tilde G_0(\psi)d\psi=\theta(\zeta,m)=\arctan f'(x_0).
\end{gather}
In the physical plane,
$\varphi\in C^{2}(\overline\Omega_{\zeta}
\setminus\{(x_0,y_0)\})
\cap C^{1}(\overline\Omega_{\zeta})$,
and $\Gamma_{\zeta}$ is a graph of function with respect to $y$
on $[0,y_0]$, which is denoted by $\tilde\Upsilon$.
Set
\begin{align*}
\tilde q_0(y)=|\nabla\varphi(\tilde\Upsilon(y),y)|,\quad
\quad y\in[0,y_0].
\end{align*}
Then, $\tilde\Upsilon\in C^{2}([0,y_0))\cap C^{1}([0,y_0])$
and $\tilde q_0\in C^{1}([0,y_0))\cap C([0,y_0])$ satisfying
\begin{gather}
\label{gg-phy1}
\tilde\Upsilon(y_0)=x_0,\quad
\tilde\Upsilon'(0)=0,\quad
\tilde\Upsilon'(y_0)=-f'(x_0),
\quad
\tilde q'_0(0)=0,\quad
\tilde q_0(y_0)=c^*,
\\[3mm]
\label{gg-phy2}
c_*<\tilde q_0(y)<c^*,\quad
|\tilde q'_0(y)|\le\frac{-\tilde\Upsilon''(y)}{1+(\tilde\Upsilon'(y))^2}\,
\sqrt{\frac{-\tilde q_0^2(y)\rho(\tilde q_0^2(y))}
{\rho(\tilde q_0^2(y))+2\tilde q_0^2(y)\rho'(\tilde q_0^2(y))}}\,,
\quad 0\le y<y_0.
\end{gather}

\subsection{Sufficient conditions for formation of vacuum}

\begin{proposition}
\label{sss-thm3}
Assume that $f\in C^{2}([l_0,l_1))$ satisfies \eqref{sss-0},
and $Q_0\in C^{1}([0,m])$ and $G_0\in C([0,m])$
satisfy \eqref{sss-1} and \eqref{sss-2}.
Let $Q\in C^{1}([0,\zeta)\times[0,m])$ with $\zeta\le+\infty$
be a solution
to the problem \eqref{newsupp-eq}--\eqref{newsupp-q}.
Then
$$
\arctan f'(\check x)-\arctan f'(\hat x)\le
{\mathscr H}(Q(\Phi_{\text{\rm up}}(\check x),m))
-{\mathscr H}(Q(\Phi_{\text{\rm up}}(\hat x),m)),\quad
l_0\le\hat x\le \check x<X_{\text{\rm up}}(\zeta).
$$
\end{proposition}

\Proof
Theorem \ref{sss-thm1} shows
\begin{align*}
W(\varphi,\psi)\le0,
\quad
Z(\varphi,\psi)\ge0,
\quad(\varphi,\psi)\in[0,\zeta)\times[0,m],
\end{align*}
which, together with \eqref{newsupp-nubbc}, leads to
\begin{align*}
Q_\varphi(\varphi,m)=\frac{1}{2}(W(\varphi,m)-Z(\varphi,m))
\le&-\frac{1}{2}(W(\varphi,m)+Z(\varphi,m))
\\
=&-\frac{f''(X_{\text{\rm up}}(\varphi))(1+(f'(X_{\text{\rm up}}(\varphi)))^2)^{-3/2}}
{b^{1/2}(Q(\varphi,m))A^{-1}_+(Q(\varphi,m))},\quad0<\varphi<\zeta.
\end{align*}
Therefore, for any $l_0\le\hat x\le \check x<X_{\text{\rm up}}(\zeta)$,
\begin{align*}
{\mathscr H}(Q(\Phi_{\text{\rm up}}(\hat x),m))
-{\mathscr H}(Q(\Phi_{\text{\rm up}}(\check x),m))
=&\int_{\Phi_{\text{\rm up}}(\hat x)}^{\Phi_{\text{\rm up}}(\check x)}b^{1/2}(Q(\varphi,m))Q_\varphi(\varphi,m)d\varphi
\\
\le&-\int_{\Phi_{\text{\rm up}}(\hat x)}^{\Phi_{\text{\rm up}}(\check x)}
\frac{f''(X_{\text{\rm up}}(\varphi))(1+(f'(X_{\text{\rm up}}(\varphi)))^2)^{-3/2}}{A^{-1}_+(Q(\varphi,m))}d\varphi
\\
=&-\int_{\hat x}^{\check x}\frac{f''(x)}{1+(f'(x))^2}dx
=\arctan f'(\hat x)-\arctan f'(\check x).
\end{align*}
$\hfill\Box$\vskip 4mm

\begin{remark}
\label{sss-thm3-remark}
In Proposition \ref{sss-thm3},
``='' holds if and only if $W(\cdot,m)=0$ in
$(\Phi_{\text{\rm up}}(\hat x),\Phi_{\text{\rm up}}(\check x))$.
\end{remark}

\begin{remark}
Denote
$$
{\mathscr F}(\gamma)={\mathscr H}(-\infty)=\int_{-\infty}^0 b^{1/2}(s)ds,\quad \gamma>1.
$$
Then
\begin{align*}
{\mathscr F}'(\gamma)<0\mbox{ for each }\gamma>1,\quad
\lim_{\gamma\to1^+}{\mathscr F}(\gamma)=+\infty,\quad
\lim_{\gamma\to+\infty}{\mathscr F}(\gamma)=0,\quad{\mathscr F}(5/3)=\pi/2.
\end{align*}
\end{remark}

Using Proposition \ref{sss-thm3}, one can get the following
sufficient conditions for formation of vacuum.

\begin{theorem}
\label{phys-thm62}
Assume that $f\in C^{2}([l_0,l_1))$ satisfies \eqref{sss-0},
$\Upsilon\in C^{2}([0,f(l_0)])$ satisfies \eqref{xc4} and
$q_0\in C^{1}([0,f(l_0)])$
satisfies \eqref{phy-8-q}--\eqref{phy-9}.
If
\begin{align}
\label{bbz2}
\arctan f'(l_1)-\arctan f'(l_0)>{\mathscr H}(-\infty)-{\mathscr H}(Q_0(f(l_0))),
\end{align}
then the smooth supersonic flow to the problem \eqref{phy1}--\eqref{phy5}
admits a vacuum.
Furthermore, the first vacuum point in the increasing $x$-direction must form at the upper wall before or at
$$
\hat x=\inf\Big\{x\in(l_0,l_1):\arctan f'(x)-\arctan f'(l_0)\ge
{\mathscr H}(-\infty)-{\mathscr H}(Q_0(f(l_0)))\Big\}.
$$
\end{theorem}

\begin{remark}
If $q_0$ satisfies \eqref{phy-8-positive} instead of \eqref{phy-8} in
Theorem \ref{phys-thm62},
then ``$>$'' can be relaxed by ``$\ge$'' in \eqref{bbz2}
and the first vacuum point in the increasing $x$-direction must form at the upper wall before $\hat x$.
\end{remark}

\begin{theorem}
Assume that $f\in C^{2}([l_0,l_1))$ satisfies \eqref{sss-0} but is not a linear function,
$\Upsilon\in C^{2}([0,f(l_0)])$ satisfies \eqref{xc4} and $q_0\in C^{1}([0,f(l_0)])$
satisfies \eqref{phy-8-q}--\eqref{phy-9}.
Then, the global smooth supersonic flow to the problem \eqref{phy1}--\eqref{phy5}
must admit a vacuum when $q_0(f(l_0))$ is near $c^*$.
Furthermore, the first vacuum point in the increasing $x$-direction tends to $(x_*,f(x_*))$ with
$x_*=\inf\{x\in[l_0,l_1):f''(x)>0\}$ as $q_0(f(l_0))\to c^*$.
\end{theorem}

\begin{theorem}
\label{sss-thm3-4}
Assume that $f\in C^{2}([l_0,l_1))$ satisfies \eqref{sss-0},
$\Upsilon\in C^{2}([0,f(l_0)])$ satisfies \eqref{xc4} and $q_0\in C^{1}([0,f(l_0)])$
satisfies \eqref{phy-8-q}, \eqref{phy-8-positive} and \eqref{phy-9}.
If
\begin{align}
\label{sss-thm3-4-0}
\left\{
\begin{aligned}
&{\lim_{x\to+\infty}}f''(x)x^{2\gamma/(\gamma+1)}=+\infty,
\quad&&\mbox{ when } l_1=+\infty \mbox{ and }\lim_{x\to+\infty}f'(x)<+\infty,
\\
&{\lim_{x\to l_1^-}}\frac{f''(x)}{(f'(x))^3}f^{2\gamma/(\gamma+1)}(x)=+\infty,
\quad&&\mbox{ when } l_1\le+\infty \mbox{ and }{\lim_{x\to l_1^-}}f'(x)=+\infty,
\end{aligned}
\right.
\end{align}
then the smooth supersonic flow to the problem \eqref{phy1}--\eqref{phy5}
admits a vacuum.
\end{theorem}

\Proof
We prove the theorem by a contradiction argument
and assume that the problem \eqref{phy1}--\eqref{phy5}
admits a global supersonic flow without vacuum,
or equivalently the problem \eqref{newsupp-eq}--\eqref{newsupp-q}
admits a global solution.
For any $l_0\le x\le \check x<X_{\text{\rm up}}(\zeta)$,
it follows from Proposition \ref{sss-thm3} and Theorem \ref{sss-thm1} that
\begin{align}
\label{fgg-3}
\arctan f'(\check x)-\arctan f'(x)\le&
{\mathscr H}(Q(\Phi_{\text{\rm up}}(\check x),m))-{\mathscr H}(Q(\Phi_{\text{\rm up}}(x),m))
\nonumber
\\
<&
{\mathscr H}(-\infty)-{\mathscr H}(Q(\Phi_{\text{\rm up}}(x),m))
\nonumber
\\
\le&
M_1\int^{+\infty}_{-Q(\Phi_{\text{\rm up}}(x),m)}s^{-(\gamma+1)/2}ds
\end{align}
with a positive constant $M_1$ independent of $x$ and $\check x$.
Proposition \ref{sss-thm1-lb} shows
\begin{align}
\label{fgg-4}
Q(\Phi_{\text{\rm up}}(x),m)\le-M_2 (1+\Phi_{\text{\rm up}}(x))^{2/(\gamma+1)},\quad
x\in[l_0,l_1)
\end{align}
with a positive constant $M_2$ independent of $x$.
The definition of $\Phi_{\text{\rm up}}$ yields
\begin{align}
\label{fgg-5}
\Phi_{\text{\rm up}}(x)\ge c_*\int_{l_0}^{x}
\sqrt{1+(f'(s))^2}ds\ge c_*(x-{l_0})+c_*(f(x)-f(l_0)),\quad x\in[l_0,l_1).
\end{align}

Assume that $l_1=+\infty$ and $f'(+\infty)=\lim_{x\to+\infty}f'(x)<+\infty$,
which imply that $f(x)=O(x)$ as $x\to+\infty$.
Letting $\check x\to+\infty$ in \eqref{fgg-3} gives
\begin{align*}
\arctan f'(+\infty)-\arctan f'(x)\le
M_1\int^{+\infty}_{-Q(\Phi_{\text{\rm up}}(x),m)}s^{-(\gamma+1)/2}ds,
\end{align*}
which, together with \eqref{fgg-4} and \eqref{fgg-5}, leads to
\begin{align}
\label{fgg-6}
\arctan f'(+\infty)-\arctan f'(x)=O\big(x^{2/(\gamma+1)-1}\big)\mbox{ as }x\to+\infty.
\end{align}
However, the first limit in \eqref{sss-thm3-4-0} shows
\begin{align*}
\lim_{x\to+\infty}(\arctan f'(+\infty)-\arctan f'(x))x^{1-2/(\gamma+1)}=+\infty,
\end{align*}
which contradicts \eqref{fgg-6}.

Turn to the case when $l_1\le+\infty$ and ${\lim_{x\to l_1^-}}f'(x)=+\infty$,
which, together with $\lim_{x\to l_1^-}(x+f(x))=+\infty$,
imply that $\lim_{x\to l_1^-}f(x)=+\infty$ and $x=o(f(x))$ as $x\to l_1^-$.
Letting $\check x\to l_1^-$ in \eqref{fgg-3} gives
\begin{align*}
\frac\pi2-\arctan f'(x)\le
M_1\int^{+\infty}_{-Q(\Phi_{\text{\rm up}}(x),m)}s^{-(\gamma+1)/2}ds,
\end{align*}
which, together with \eqref{fgg-4} and \eqref{fgg-5}, leads to
\begin{align}
\label{fgg-7}
\frac\pi2-\arctan f'(x)=O\big(f^{2/(\gamma+1)-1}(x)\big)\mbox{ as }x\to l_1^-.
\end{align}
However, the second limit in \eqref{sss-thm3-4-0} shows
\begin{align*}
\lim_{x\to l_1^-}\Big(\frac\pi2-\arctan f'(x)\Big)f^{1-2/(\gamma+1)}(x)=+\infty,
\end{align*}
which contradicts \eqref{fgg-7}.
$\hfill\Box$\vskip 4mm

\begin{remark}
Theorem \ref{sss-thm3-4} still holds if \eqref{sss-thm3-4-0} is relaxed by
$$
\left\{
\begin{aligned}
&\overline{\lim_{x\to+\infty}}\big(\arctan f'(+\infty)-\arctan f'(x)\big)x^{1-2/(\gamma+1)}=+\infty,
&&\mbox{ when } l_1=+\infty \mbox{ and }\lim_{x\to+\infty}f'(x)<+\infty,
\\
&\overline{\lim_{x\to+\infty}}\big(\pi/2-\arctan f'(x)\big)f^{1-2/(\gamma+1)}(x)=+\infty,
&&\mbox{ when } l_1\le+\infty \mbox{ and }{\lim_{x\to l_1^-}}f'(x)=+\infty.
\end{aligned}
\right.
$$
\end{remark}

\subsection{Lower bound estimates of the location of the first vacuum point}

\begin{proposition}
\label{sss-thm4}
Assume that $f\in C^{2}([l_0,l_1))$ satisfies \eqref{sss-0},
and $Q_0\in C^{1}([0,m])$ and $G_0\in C([0,m])$
satisfy \eqref{sss-1} and \eqref{sss-2}.
Let $Q\in C^{1}([0,\zeta)\times[0,m])$ with $0<\zeta<+\infty$
be the maximal smooth supersonic flow before vacuum formation
to the problem \eqref{newsupp-eq}--\eqref{newsupp-q}. Then,
\begin{align}
\label{sd20}
&\arctan f'(X_{\text{\rm up}}(\zeta))>{\mathscr H}(-\infty)-
{\mathscr H}\Big(\frac1m\Big(\int_0^m Q_0(\psi)d\psi+\zeta\int_0^m G_0(\psi)d\psi
\nonumber
\\
&\qquad-c^*\int_{l_0}^{X_{\text{\rm up}}(\zeta)}(\arctan f'(x)-\arctan f'(l_0))
{\sqrt{1+(f'(x))^2}}dx\Big)\Big).
\end{align}
\end{proposition}

\Proof
It follows from \eqref{newsupp-eq}--\eqref{newsupp-q} that
\begin{align*}
&\frac{d^2}{d\varphi^2}\int_0^m Q(\varphi,\psi)d\psi
=-\frac{f''(X_{\text{\rm up}}(\varphi))(1+(f'(X_{\text{\rm up}}(\varphi)))^2)^{-3/2}}
{A^{-1}_+(Q(\varphi,m))},\quad \varphi\in(0,\zeta),
\\
&\int_0^m Q(\varphi,\psi)d\psi\Big|_{\varphi=0}=\int_0^m Q_0(\psi)d\psi,
\quad
\frac{d}{d\varphi}\int_0^m Q(\varphi,\psi)d\psi\Big|_{\varphi=0}=\int_0^m G_0(\psi)d\psi.
\end{align*}
Therefore,
\begin{align}
\label{sd21}
\int_0^m Q(\zeta,\psi)d\psi
=&\int_0^m Q_0(\psi)d\psi+\zeta\int_0^m G_0(\psi)d\psi
\nonumber
\\
&\qquad-\int_{0}^\zeta\int_{0}^s
\frac{f''(X_{\text{\rm up}}(\varphi))(1+(f'(X_{\text{\rm up}}(\varphi)))^2)^{-3/2}}
{A^{-1}_+(Q(\varphi,m))}d\varphi ds.
\end{align}
Direct calculation gives
\begin{align}
\label{sd22}
&\int_{0}^\zeta\int_{0}^s
\frac{f''(X_{\text{\rm up}}(\varphi))(1+(f'(X_{\text{\rm up}}(\varphi)))^2)^{-3/2}}
{A^{-1}_+(Q(\varphi,m))}d\varphi ds
\nonumber
\\
=&\int_{l_0}^{X_{\text{\rm up}}(\zeta)}
{A^{-1}_+(Q(t,m))}\sqrt{1+(f'(t))^2}\int_{l_0}^t
\frac{f''(x)}{1+(f'(x))^2}dx dt
\nonumber
\\
<&c^*\int_{l_0}^{X_{\text{\rm up}}(\zeta)}
(\arctan f'(x)-\arctan f'(l_0))\sqrt{1+(f'(x))^2}dx.
\end{align}
By \eqref{sd21} and \eqref{sd22}, there exists $\psi_0\in(0,m)$ such that
\begin{align}
\label{sd23}
Q(\zeta,\psi_0)
>&\frac1m\Big(\int_0^m Q_0(\psi)d\psi+\zeta\int_0^m G_0(\psi)d\psi
\nonumber
\\
&\qquad-c^*\int_{l_0}^{X_{\text{\rm up}}(\zeta)}
(\arctan f'(x)-\arctan f'(l_0))\sqrt{1+(f'(x))^2}dx\Big).
\end{align}
It follows from Theorem \ref{sss-thm1} that
\begin{align*}
\pd{}{\psi}\big({\mathscr H}(Q(\varphi,\psi))-\theta(\varphi,\psi)\big)=
Q_\varphi(\varphi,\psi)-b^{1/2}(Q(\varphi,\psi))
Q_\psi(\varphi,\psi)\le0,\quad(\varphi,\psi)\in(0,\zeta)\times(0,m),
\end{align*}
which, together with Theorem \ref{sss-thm2}, leads to
\begin{align}
\label{sd24}
{\mathscr H}(-\infty)-\arctan f'(X_{\text{\rm up}}(\zeta))
\le {\mathscr H}(Q(\zeta,\psi_0))-\theta(\zeta,\psi_0)\le{\mathscr H}(Q(\zeta,\psi_0)).
\end{align}
Then, \eqref{sd20} follows from \eqref{sd23} and \eqref{sd24}.
$\hfill\Box$\vskip 4mm

\begin{theorem}
\label{sss-thm4-2}
Assume that $f\in C^{2}([l_0,l_1))$ satisfies \eqref{sss-0} but is not a linear function,
$\Upsilon\in C^{2}([0,f(l_0)])$ satisfies \eqref{xc4} and
$q_0\in C^{1}([0,f(l_0)])$
satisfies \eqref{phy-8-q}--\eqref{phy-9}.
Set
\begin{align*}
\check x=\inf\Big\{x\in(l_0,l_1):\arctan f'(x)\ge{\mathscr H}(-\infty)
-{\mathscr H}(h(x))\Big\},
\end{align*}
where for $x\in(l_0,l_1)$,
$$
h(x)=\frac1m\Big(\int_0^m Q_0(\psi)d\psi+\zeta\int_0^m G_0(\psi)d\psi
-c^*\int_{l_0}^{x}(\arctan f'(s)-\arctan f'(l_0))
{\sqrt{1+(f'(s))^2}}ds\Big).
$$
Then there is no vacuum on $\{(x,f(x)):l_0\le x\le \check x\}$
for the global smooth transonic flow to the problem \eqref{phy1}--\eqref{phy5}.
\end{theorem}

\subsection{Formation of a shock in a
nozzle whose upper wall is a non-convex perturbation of a straight line}

Below, we explain the condition that the upper wall of the nozzle is convex.
Indeed, if it is a non-convex perturbation of a straight line,
then there is a $q_0\in C^{1}([0,f(l_0)])$ satisfying \eqref{phy-8-q}--\eqref{phy-9}
such that a shock must form for the supersonic flow to the problem \eqref{phy1}--\eqref{phy5}.

\begin{theorem}
\label{ni2-thm-2}
Assume that
$f\in C^2([l_0,+\infty))$ satisfies
\begin{align*}
f(l_0)>0,\quad f'(l_0)>0,\quad f''(l_0)=0,\quad
\sup_{(l_0,+\infty)}|f''|\le\varepsilon,\quad \inf_{(l_0,+\infty)}f''<0,
\end{align*}
and
$\Upsilon\in C^{2}([0,f(l_0)])$ satisfies \eqref{xc4} and $\Upsilon''(0)=\Upsilon''(f(l_0))=0$,
where $\varepsilon\in(0,1)$ is a constant.
Then, there exist a function $q_0\in C^{1}([0,f(l_0)])$ satisfying \eqref{phy-8-q}--\eqref{phy-9}
and a positive constant $\varepsilon_0\in(0,1)$
such that if $0<\varepsilon<\varepsilon_0$,
then a shock must form for the global supersonic flow to the problem \eqref{phy1}--\eqref{phy5}.
\end{theorem}

\Proof
Set
$l_*=\inf\big\{x>l_0:f''(x)<0\big\}$ and $\zeta_*=c^*(1+f'(l_0)+l_*-l_0)(l_*-l_0)$.
Without loss of generality, it is assumed that $l_*>l_0$.
Choose $0<\delta_0<(c^*-c_*)/2$ such that
\begin{align}
\label{nonv-1}
\zeta_*b^{1/2}(A(c^*-2\delta_0))<m,
\quad
{\mathscr H}(A(c^*-2\delta_0))+\arctan f'(l_0)>{\mathscr H}(-\infty).
\end{align}
For $0<\delta<\delta_0$, consider the problem
\begin{gather}
\label{nonv-2}
q'_0(y)=\frac{-\Upsilon''(y)}{1+(\Upsilon'(y))^2}\,
\sqrt{\frac{-q_0^2(y)\rho(q_0^2(y))}{\rho(q_0^2(y))+2q_0^2(y)\rho'(q_0^2(y))}}\,,
\quad 0<y<f(l_0),
\qquad
q(l_0)=c^*-\delta.
\end{gather}
Since
$$
\lim_{s\to {c^*}^-}\sqrt{\frac{-s^2\rho(s^2)}{\rho(s^2)+2s^2\rho'(s^2)}}=0,
$$
there exists $\delta\in(0,\delta_0)$ such that the problem \eqref{nonv-2}
admits a solution $q_0\in C^{1}([0,f(l_0)])$ satisfying
$q_0(0)\ge c^*-\delta_0$.
Consider the problem \eqref{phy1}--\eqref{phy5} with this $q_0$,
which satisfies \eqref{phy-8-q}--\eqref{phy-9}.
Then, in the potential plane, the boundary data at the inlet satisfy
\begin{align}
\label{nonv-4}
Q_0(\psi)\le A(c^*-\delta_0)<A(c^*-2\delta_0),\quad
W_0(\psi)=0,\quad Z_0(\psi)\ge0,\quad 0<\psi<m.
\end{align}
Note that $f''\ge0$ on $[l_0,l_*]$.
According to Theorem \ref{sss-thm1},
the problem \eqref{newsupp-eq}--\eqref{newsupp-q} admits uniquely a solution
$Q\in C^{1}([0,\zeta)\times[0,m])$ for some $\zeta>0$,
where either
$\zeta=\Phi_{\text{\rm up}}(l_*)$
or
$X_{\text{\rm up}}(\zeta)<l_*$ with $\lim_{\varphi\to\zeta^-}\min_{[0,m]}Q(\varphi,\cdot)=-\infty$.
Furthermore, the solution satisfies
\begin{align}
\label{nonv-5}
-Q_\varphi(\varphi,\psi)\pm b^{1/2}(Q(\varphi,\psi))Q_\psi(\varphi,\psi)\ge0,\quad
Q(\varphi,\psi)\le A(c^*-\delta_0),\quad
(\varphi,\psi)\in(0,\zeta)\times(0,m).
\end{align}
It is clear that $\zeta<\zeta_*$.
By Theorem \ref{sss-thm2},
Proposition \ref{sss-thm3} and Remark \ref{sss-thm3-remark},
there exists a positive constant $\varepsilon_0\in(0,1)$
such that $\zeta=\Phi_{\text{\rm up}}(l_*)$ and $Q\in C^{1}([0,\zeta]\times[0,m])$
for each $0<\varepsilon<\varepsilon_0$.
Furthermore, \eqref{nonv-1}, \eqref{nonv-5} and Lemma \ref{lemma-sc4} imply that
the positive characteristic from $(\zeta,m)$ must touch $\{0\}\times(0,m)$,
while the negative characteristic from $(\zeta,m)$ never touches $(0,+\infty)\times\{0\}$.
Take $\tau>0$ sufficiently small such that $Q$ exists on $[\zeta,\zeta+\tau]\times[0,m]$,
the positive characteristic from $(\zeta+\tau,m)$ touches $\{0\}\times(0,m)$ and
\begin{gather}
\label{nonv-6}
Q(\varphi,m)\le A(c^*-2\delta_0),\quad
{\mathscr H}(A(c^*-2\delta_0))\ge\theta(\zeta,m)-\theta(\varphi,m),\quad
\varphi\in(\zeta,\zeta+\tau).
\end{gather}
It follows from \eqref{newsupp-neq1} and the second formula in \eqref{nonv-4} that
$W(\cdot,m)=0$ on $[0,\zeta+\tau]$.

\vskip10mm

\hskip30mm
\setlength{\unitlength}{0.6mm}
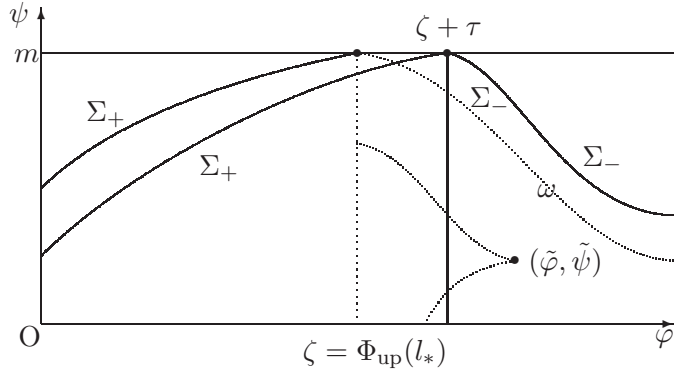
\begin{picture}(250,65)
\put(0,0){\vector(1,0){140}}
\put(0,0){\vector(0,1){70}}
\put(136,-4){$\varphi$} \put(-7,67){$\psi$}
\put(0,60){\qbezier(0,0)(55,0)(140,0)}

\put(0,0){\cbezier(0,30)(20,50)(60,58)(70,60)}

\put(70,60){\cbezier[100](0,0)(30,-4)(45,-45)(70,-46)}

\put(0,0){\cbezier(0,15)(20,35)(60,56)(90,60)}

\put(90,60){\cbezier(0,0)(15,-3)(25,-35)(50,-36)}

\put(70,60){\circle*{1.5}}
\put(90,60){\circle*{1.5}}

\put(120,35){$\Sigma_-$}

\put(95,48){$\Sigma_-$}

\put(10,45){$\Sigma_+$}

\put(35,33){$\Sigma_+$}

\put(90,0){\qbezier(0,0)(0,30)(0,60)}

\put(110,28){$\omega$}

\put(108,12){$(\tilde\varphi,\tilde\psi)$}
\put(105,14){\circle*{1.5}}
\put(105,14){\cbezier[60](0,0)(-15,3)(-20,23)(-35,26)}
\put(105,14){\cbezier[30](0,0)(-10,-2)(-15,-5)(-20,-14)}

\put(-6,58){$m$}

\put(-5,-5){O}

\put(70,0){\qbezier[40](0,0)(0,30)(0,60)}
\put(58,-8){$\zeta=\Phi_{\text{\rm up}}(l_*)$}

\put(83,64){$\zeta+\tau$}

\end{picture}

\vskip10mm

Let $\omega$ be the domain bounded by $\varphi=\zeta+\tau$, $\psi=0$,
and the negative characteristic from $(\zeta+\tau,m)$.
Lemma \ref{aaasss-thm2} gives
\begin{align}
\label{nonv-7}
-Q_\varphi(\varphi,\psi)\le M_1,\quad
(\varphi,\psi)\in\omega
\end{align}
with some positive constant $M_1$.
It follows from Lemma \ref{lemma-sc4} and \eqref{nonv-5} that
\begin{align}
\label{nonv-7-1}
Q(\varphi,0)\le A(c^*-\delta_0),\quad \varphi>\zeta.
\end{align}
Fix $(\tilde\varphi,\tilde\psi)\in\omega$.
Assume that
$$
\tilde\Sigma_{\pm}:
\tilde\Psi'_{\pm}(\varphi)=\pm b^{1/2}(Q(\varphi,\tilde\Psi_{\pm}(\varphi))),
\quad
\tilde\varphi_{\pm}<\varphi<\tilde\varphi,
\quad\tilde\Psi_{\pm}(\tilde\varphi_{\pm})=\tilde\psi_{\pm},\quad\tilde\Psi_{\pm}(\tilde\varphi)=\tilde\psi
$$
are the positive and negative characteristics approaching $(\tilde\varphi,\tilde\psi)$,
where
\begin{align}
\label{nonv-7-2}
(\tilde\varphi_+,\tilde\psi_+)\in\{\zeta\}\times(0,m)\cup [\zeta,+\infty)\times\{0\},
\quad
(\tilde\varphi_-,\tilde\psi_-)\in\{\zeta\}\times(0,m)\cup [\zeta,\zeta+\tau]\times\{m\}.
\end{align}
Lemma \ref{lemma-sc4} gives
\begin{align}
\label{nonv-7-3}
{\mathscr H}(Q(\tilde\varphi_\pm,\tilde\psi_\pm))\mp \theta(\tilde\varphi_\pm,\tilde\psi_\pm)
={\mathscr H}(Q(\tilde\varphi,\tilde\psi))\mp \theta(\tilde\varphi,\tilde\psi).
\end{align}
Then, \eqref{nonv-5}, \eqref{nonv-6}, \eqref{nonv-7-1}--\eqref{nonv-7-3} lead to
\begin{align*}
2{\mathscr H}(Q(\tilde\varphi,\tilde\psi))
=&{\mathscr H}(Q(\tilde\varphi_+,\tilde\psi_+))
+{\mathscr H}(Q(\tilde\varphi_-,\tilde\psi_-))-\theta(\tilde\varphi_+,\tilde\psi_+)+\theta(\tilde\varphi_-,\tilde\psi_-)
\\
\ge&{\mathscr H}(A(c^*-\delta_0))+{\mathscr H}(A(c^*-2\delta_0))-\theta(\tilde\varphi_+,\tilde\psi_+)+\theta(\tilde\varphi_-,\tilde\psi_-)
\ge{\mathscr H}(A(c^*-\delta_0)).
\end{align*}
Therefore, there exists a positive constant $M_2>0$ such that
\begin{align}
\label{nonv-7-4}
Q(\tilde\varphi,\tilde\psi)\le-M_2,\quad(\tilde\varphi,\tilde\psi)\in\omega,
\end{align}
which, together with \eqref{nonv-7}, shows that the flow admits no sonic and vacuum points
on $\overline\omega$.
According to the definition of $l_*$, there exists $0<\tilde\tau<\tau$ such that
$W(\zeta+\tilde\tau,m)+Z(\zeta+\tilde\tau,m)<0$, which, together with $W(\cdot,m)=0$ on $[0,\zeta+\tau]$,
leads to $Z(\zeta+\tilde\tau,m)<0$.
Then, Lemma \ref{ni}, \eqref{nonv-7} and \eqref{nonv-7-4} show that a shock must form in $\omega$.
$\hfill\Box$\vskip 4mm

\section{Global smooth supersonic flows}

In this section, we first investigate the properties of a vacuum boundary for a supersonic flow
and then consider the dynamics of the vacuum boundary separating the vacuum region from the non-vacuum  supersonic flow as a
free boundary problem of a smooth supersonic flow
from the potential level set where the first vacuum forms.
Subsequently, we solve the problem \eqref{phy1}--\eqref{phy5} globally.
Moreover, smooth supersonic flows with vacuum at the inlet are also considered.

\subsection{Properties of the vacuum boundary}

Generally, the vacuum boundary of a supersonic flow possesses the following properties.

\begin{theorem}
\label{asd}
Let $\Omega$ be a bounded domain and $\Gamma\subset\partial\Omega$ be a connected and open $C^1$ curve.
Assume that $\varphi\in C^1(\Omega\cup\Gamma)$ is a supersonic flow such that
it is vacuum on $\Gamma$ and there is no vacuum in $\Omega$.
Then at each point of $\Gamma$, the flow velocity is along $\Gamma$.
Furthermore,

{\rm (i)}
If $\varphi\in C^2(\Omega)$ additionally,
then $\Gamma$ is a straight line.

{\rm (ii)}
If $\varphi\in C^2(\Omega\cup\Gamma)$ additionally,
then $\pd{}{\vec\nu}|\nabla\varphi|=0$ on $\Gamma$,
where $\vec\nu$ is the normal to $\Gamma$.
\end{theorem}

\Proof
Assume that
there is a point $P_0\in\Gamma$ such that at $P_0$
the flow velocity is not along the tangent direction of $\Gamma$.
Fix a point $\tilde P_0$, which is located at the different side of
$\Gamma$ from $\Omega$, such that
the flow velocity at $P_0$ is along
the direct from $\tilde P_0$ to $P_0$ or its opposite direction.
For $r>0$, denote $\Omega_r$ to be the set of points
which belong to $\Omega$ and whose distances to $\tilde P_0$ are less than $r$.
Since $\varphi\in C^1(\Omega\cup\Gamma)$ and $\Gamma\in C^1$, there exists a number $r_0$,
which is greater than the distance from $\tilde P_0$ to $P_0$, such that
$\partial\Omega\cap\Omega_{r_0}\subset\Gamma$ and
for each point $P\in\Omega_{r_0}$, the angle between
the flow velocity at $P$ and the direction from $\tilde P_0$ to $P$
is always acute (if the flow velocity at $P_0$ is along the direct from $\tilde P_0$ to $P_0$)
or always obtuse
(if the flow velocity at $P_0$ is along the direct from $P_0$ to $\tilde P_0$).
Note that the flow is vacuum on $\Gamma$.
Hence
$$
\int_{{\partial{\Omega_{r}}}{\cap}{\Omega}}\rho(|\nabla\varphi(x,y)|^2)\nabla\varphi(x,y)\cdot
\vec\nu_r(x,y) d\sigma=0,\quad0<r\le r_0,
$$
where $\vec\nu_r$ is the normal to $\partial{\Omega_{r}}$.
The choice of $r_0$ implies that $\nabla\varphi\cdot\vec\nu_r$
is always positive or always negative on ${\partial{\Omega_{r}}}{\cap}{\Omega}$ for each $0<r\le r_0$.
Therefore, the flow is vacuum on $\overline\Omega_{r_0}$, which is a contradiction.

\vskip5mm

\hskip27mm
\setlength{\unitlength}{0.6mm}
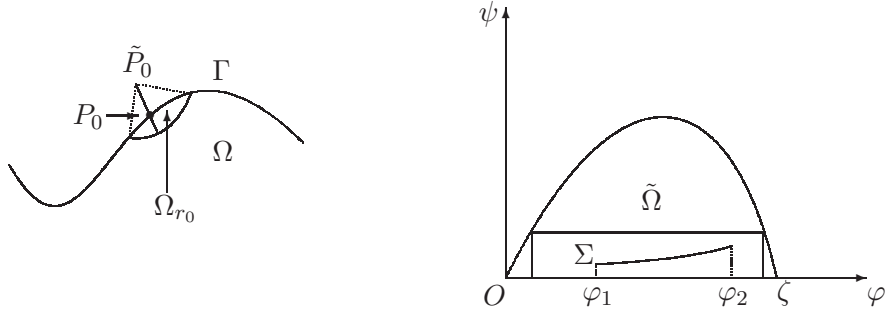
\begin{picture}(250,55)


\put(18,23){\qbezier(0,0)(2.5,-5.6)(5,-11.2)}
\put(18,23){\qbezier[16](0,0)(6.1,-1.1)(12.2,-2.2)}
\put(18,23){\qbezier[16](0,0)(-0.7,-6)(-1.4,-12)}

\put(18,23){\qbezier(12.2,-2.2)(7.8,-13)(-1.4,-12)}

\put(0,0){\cbezier(-10,5)(10,-30)(15,50)(55,10)}

\put(21.15,16.1){\circle*{1.5}}

\put(15,25.6){$\tilde P_0$}

\put(4,14){$P_0$}
\put(11.5,16){\vector(1,0){7}}

\put(35,24){$\Gamma$}
\put(35,6){$\Omega$}

\put(22,-6){$\Omega_{r_0}$}
\put(25,-1){\vector(0,1){18}}

\put(100,-20){\vector(1,0){80}}
\put(100,-20){\vector(0,1){60}}
\put(100,-20){\cbezier(0,0)(25,50)(50,45)(60,0)}
\put(95,-26){$O$}

\put(110,-10){\qbezier(-4.2,0)(15,0)(47,0)}
\put(105.8,-20){\qbezier(0,0)(0,5)(0,10)}
\put(157,-20){\qbezier(0,0)(0,5)(0,10)}

\put(120,-20){\qbezier(0,3)(20,4)(30,7)}

\put(120,-20){\qbezier[5](0,0)(0,1.5)(0,3)}
\put(150,-20){\qbezier[11](0,0)(0,3.5)(0,7)}

\put(115,-17){$\Sigma$}

\put(117,-25){$\varphi_1$}
\put(147,-25){$\varphi_2$}

\put(160,-25){$\zeta$}
\put(180,-25){$\varphi$}
\put(94,37){$\psi$}

\put(130,-5){$\tilde\Omega$}

\end{picture}
\vskip20mm

Transform the supersonic flow on $\Omega\cup\Gamma$ into the potential plane.
The transformation is one to one since
the flow velocity is along $\Gamma$ at each point of $\Gamma$, which also yields that
$\Gamma$ is a streamline.
Assume that $\Omega$ and $\Gamma$ are transformed into
$\tilde\Omega$ and $\tilde\Gamma$ in the coordinates transformation
from the physical plane to the potential plane.
Then, $\tilde\Gamma\subset\partial\tilde\Omega$.
Without loss of generality,
it is assumed that $\tilde\Gamma=(0,\zeta)\times\{0\}$
and $\{\zeta/2\}\times(0,\delta_0)\subset\tilde\Omega$, where
$\zeta$ and $\delta_0$ are positive constants.
For given $0<\varepsilon<\zeta/2$, there exists a positive constant $\delta<\delta_0$ such that
$[\varepsilon,\zeta-\varepsilon]\times(0,\delta)\subset\tilde\Omega$.
Assume that
$$
\Sigma:
\Psi'(\varphi)=b^{1/2}(A(q(\varphi,\Psi(\varphi)))),
\quad\varphi_1<\varphi<\varphi_2,
\quad\Psi(\varphi_1)=\psi_1
$$
is a positive characteristic in $[\varepsilon,\zeta-\varepsilon]\times(0,\delta)$.
Then, for each positive integer $n$,
$$
\Sigma_n:
\Psi'_n(\varphi)=b^{1/2}(A(q(\varphi,\Psi_n(\varphi)))),
\quad\varphi_1<\varphi<\varphi_2,
\quad\Psi_n(\varphi_1)={\psi_1}/{n}
$$
is also a positive characteristic in $[\varepsilon,\zeta-\varepsilon]\times(0,\delta)$,
and
\begin{align}
\label{asd1}
\lim_{n\to+\infty}\Psi_n(\varphi)=0,\quad
\varphi_1\le\varphi\le\varphi_2.
\end{align}
For each positive integer $n$,
$\varphi\in C^2(\Omega)$ and Lemma \ref{lemma-sc4} yield
\begin{align}
\label{asd2}
{\mathscr H}(A(q(\varphi,\Psi_{n}(\varphi))))-\theta(\varphi,\Psi_{n}(\varphi))=M_n,\quad
\varphi_1\le\varphi\le\varphi_2
\end{align}
with $M_n$ being a constant.
Note that
\begin{align}
\label{asd3}
\lim_{\psi\to0^+}q(\varphi,\psi)=c^*,\quad\varphi_1\le\varphi\le\varphi_2.
\end{align}
Then, \eqref{asd1}--\eqref{asd3} and $\varphi\in C^1(\Omega\cup\Gamma)$ lead to
that $\theta(\cdot,0)$ is invariant on $[\varphi_1,\varphi_2]$.
Since $\Sigma$ and $\varepsilon$ are arbitrary,
$\theta(\cdot,0)$ is equal to a constant identically in $(0,\zeta)$.
That is to say, the flow angle is invariant on $\Gamma$,
which, together with $\Gamma$ being a streamline, implies that
$\Gamma$ is a straight line.

Assume that $\varepsilon$ and $\delta$ are given as above.
For each $\xi\in C_0^\infty(\varepsilon,\zeta-\varepsilon)$, one gets
from $\varphi\in C^2(\Omega\cup\Gamma)$ and (i) that
\begin{align*}
&\int_\varepsilon^{\zeta-\varepsilon}\pd{B(q)}{\psi}(\varphi,0)\xi(\varphi)d\varphi
=\lim_{\psi\to 0^+}\int_\varepsilon^{\zeta-\varepsilon}\pd{B(q)}{\psi}(\varphi,\psi)\xi(\varphi)d\varphi
\\
=&\lim_{\psi\to 0^+}\int_\varepsilon^{\zeta-\varepsilon}
\pd{\theta}{\varphi}(\varphi,\psi)\xi(\varphi)d\varphi
=-\lim_{\psi\to 0^+}\int_\varepsilon^{\zeta-\varepsilon}\theta(\varphi,\psi)\xi'(\varphi)d\varphi
=0.
\end{align*}
Then, it follows from the arbitrariness of $\xi$ and $\varepsilon$ that
$$
\pd{B(q)}{\psi}(\varphi,0)=0,\quad0<\varphi<\zeta,
$$
which yields (ii) by the transformation from the potential plane to the physical plane.
$\hfill\Box$\vskip 4mm

\begin{remark}
Theorem \ref{asd} shows that
the vacuum of a $C^1$ supersonic potential flow
is never the so called physical vacuum in \cite{YL1,YL2}.
The reason lies in that there is not an external force on the flow.
\end{remark}

\subsection{Dynamics of the vacuum boundary and a free boundary problem
}

In this subsection, we study the dynamic of the vacuum boundary from its first formation
in the potential plane and treat it as a free boundary problem to extend globally a supersonic flow
from the potential level set where the first vacuum forms.
Assume that $\tilde Q_0\in C^{1}([0,m))$ and $\tilde G_0\in C([0,m))\cap L^1(0,m)$
satisfy \eqref{gg-1} and \eqref{gg-2}.
For $0<n<m$,
choose $\tilde Q_{0,n}\in C^{1}([0,m])$ and $\tilde G_{0,n}\in C([0,m])$ to satisfy
$\tilde Q'_{0,n}(m)=0$ and
\begin{gather*}
\tilde Q_{0,n}(\psi)=\tilde Q_{0}(\psi),\quad
\tilde G_{0,n}(\psi)=\tilde G_{0}(\psi),\quad\psi\in[0,n],
\\
\tilde Q_{0,n}(\psi)<0,\quad
-\tilde G_{0,n}(\psi)\pm b^{1/2}(\tilde Q_{0,n}(\psi))\tilde Q'_{0,n}(\psi)\ge0,\quad\psi\in[n,m].
\end{gather*}
Consider the problem
\begin{align}
\label{gg-3}
&\pd{^2Q_n}{\varphi^2}-\pd{}{\psi}\Big(b(Q_n)\pd{Q_n}{\psi}\Big)=0,
\quad&&(\varphi,\psi)\in(\zeta,+\infty)\times(0,m),
\\
\label{gg-4}
&{Q_n}(\zeta,\psi)=\tilde Q_{0,n}(\psi),\quad&&\psi\in(0,m),
\\
\label{gg-5}
&\pd{Q_n}{\varphi}(\zeta,\psi)=\tilde G_{0,n}(\psi),\quad&&\psi\in(0,m),
\\
\label{gg-8}
&\pd{Q_n}{\psi}(\varphi,0)=\pd{Q_n}{\psi}(\varphi,m)=0,
\quad&&\varphi\in(\zeta,+\infty).
\end{align}
It follows from Theorem \ref{phys-thm51} that
the problem \eqref{gg-3}--\eqref{gg-8} admits a solution
$Q_n\in C^{1}([\zeta,+\infty)\times[0,m])$ satisfying
\begin{align}
\label{gg-10}
-\pd{Q_n}\varphi(\varphi,\psi)
\pm b^{1/2}(Q_n(\varphi,\psi))\pd{Q_n}\psi(\varphi,\psi)\ge0,\quad
(\varphi,\psi)\in(\zeta,+\infty)\times(0,m).
\end{align}

\begin{lemma}
\label{gg-prop2}
Let $0<k\le n<m$.
Assume that
$$
\Sigma_n:
\left\{
\begin{aligned}
&\Psi'_n(\varphi)=-b^{1/2}(Q_{n}(\varphi,\Psi_n(\varphi))),
\quad
0<\Psi_n(\varphi)<k,\quad\zeta<\varphi<\varphi_n\le+\infty,
\\
&\Psi_n(\zeta)=k,\quad \Psi_n(\varphi_n)=0\mbox{ if }\varphi_n<+\infty
\end{aligned}
\right.
$$
is the negative characteristic from $(\zeta,k)$.
Then $\varphi_n\ge\zeta+b^{-1/2}(\tilde Q_0(k))k$ and
\begin{align*}
\Psi_n(\varphi)\ge k-b^{1/2}(\tilde Q_0(k))(\varphi-\zeta),
\quad
\zeta\le\varphi\le\zeta+b^{-1/2}(\tilde Q_0(k))k.
\end{align*}
\end{lemma}

\Proof
On $\Sigma_n$, $Q_{n}$ satisfies
\begin{align*}
\frac{d}{d\varphi}Q_{n}(\varphi,\Psi_n(\varphi))
=\pd{Q_{n}}\varphi(\varphi,\Psi_n(\varphi))
-b^{1/2}(Q_{n}(\varphi,\Psi_n(\varphi)))\pd{Q_{n}}\psi(\varphi,\Psi_n(\varphi))\le0,\quad
\zeta<\varphi<\varphi_n
\end{align*}
according to \eqref{gg-10}. Therefore,
\begin{align*}
Q_{n}(\varphi,\Psi_n(\varphi))\le Q_{n}(\zeta,\Psi_n(\zeta))=\tilde Q_{0,n}(k)=\tilde Q_0(k),\quad
\zeta<\varphi<\varphi_n,
\end{align*}
which implies $\Psi'_n\ge-b^{1/2}(\tilde Q_0(k))$ in $(\zeta,\varphi_n)$
and thus the lemma is proved.
$\hfill\Box$\vskip 4mm

For $0<k\le n<m$, Lemma \ref{gg-prop2} shows
\begin{align}
\label{newg-3}
Q_{n}(\varphi,\psi)=Q_{k}(\varphi,\psi),
\quad
\zeta\le\varphi\le\zeta+b^{-1/2}(\tilde Q_0(k))k,\,0\le\psi\le k-b^{1/2}(\tilde Q_0(k))(\varphi-\zeta).
\end{align}
By \eqref{newg-3} and \eqref{gg-2},
one can set
\begin{align}
\label{gg-11}
Q(\varphi,\psi)=\lim_{n\to m^-}{Q_n}(\varphi,\psi),\quad
(\varphi,\psi)\in[\zeta,+\infty)\times[0,m).
\end{align}

\begin{proposition}
\label{gg-thm1}
Assume that $\tilde Q_0\in C^{1}([0,m))$ and $\tilde G_0\in C([0,m))\cap L^1(0,m)$
satisfy \eqref{gg-1} and \eqref{gg-2}.
Then the function $Q$ defined by \eqref{gg-11} belongs to
$C^{1}([\zeta,+\infty)\times[0,m))$
and is the unique solution to the problem
\begin{align}
\label{gg-13}
&Q_{\varphi\varphi}-(b(Q)Q_{\psi})_\psi=0,
\quad&&(\varphi,\psi)\in(\zeta,+\infty)\times(0,m),
\\
\label{gg-14}
&{Q}(\zeta,\psi)=\tilde Q_0(\psi),\quad&&\psi\in(0,m),
\\
\label{gg-15}
&{Q}_{\varphi}(\zeta,\psi)=\tilde G_0(\psi),\quad&&\psi\in(0,m),
\\
\label{gg-16}
&Q_{\psi}(\varphi,0)=0,
\quad&&\varphi\in(\zeta,+\infty).
\end{align}
Furthermore,
\begin{gather}
\label{gg-17-0}
-Q_\varphi(\varphi,\psi)\pm b^{1/2}(Q(\varphi,\psi))Q_\psi(\varphi,\psi)\ge0,\quad(\varphi,\psi)\in(\zeta,+\infty)\times(0,m),
\end{gather}
and $q,\theta\in C([\zeta,+\infty)\times[0,m])$ satisfying
\begin{gather}
\label{gg-17-1}
q(\varphi,\psi)<c^*,\quad \theta(\varphi,\psi)<\theta(\zeta,m),
\quad(\varphi,\psi)\in[\zeta,+\infty)\times[0,m),
\\
\label{gg-17}
\lim_{\psi\to m^-}q(\varphi,\psi)=c^*\mbox{ and }
\lim_{\psi\to m^-}\theta(\varphi,\psi)=\theta(\zeta,m)
\mbox{ uniformly for }\varphi\ge\zeta.
\end{gather}
\end{proposition}

\Proof
It follows from Lemma \ref{gg-prop2} and \eqref{gg-2} that
$Q\in C^{1}([\zeta,+\infty)\times[0,m))$ is the unique solution to the problem \eqref{gg-13}--\eqref{gg-16},
and \eqref{gg-17-0} follows from \eqref{gg-10}.
The first estimate in \eqref{gg-17-1} is clear since $Q\in C^{1}([\zeta,+\infty)\times[0,m))$.
Moreover, \eqref{gg-2} and \eqref{gg-17-0} lead to the first limit in \eqref{gg-17}.

We now show
\begin{gather}
\label{gg-17-2}
\theta(\zeta,\psi)<\theta(\zeta,m),\quad0<\psi<m
\end{gather}
by a contradiction argument.
Otherwise, by \eqref{gg-1}, there exists $\psi_0\in(0,m)$ such that
$\theta(\zeta,\cdot)=\theta(\zeta,m)$ in $(\psi_0,m)$,
which yields $\tilde G_0=0$ in $(\psi_0,m)$.
Then, \eqref{gg-1} implies $\tilde Q'_0=0$ in $(\psi_0,m)$,
which contradicts the second formula in \eqref{gg-2} and thus \eqref{gg-17-2} is proved.

For $0<n<m$, set
$$
\Sigma_n:
\left\{
\begin{aligned}
&\Psi'_n(\varphi)=-b^{1/2}(Q(\varphi,\Psi_n(\varphi))),
\quad
0<\Psi_n(\varphi)<n,\quad\zeta<\varphi<\varphi_n\le+\infty,
\\
&\Psi_n(\zeta)=n,\quad \Psi_n(\varphi_n)=0\mbox{ if }\varphi_n<+\infty
\end{aligned}
\right.
$$
to be the negative characteristic from $(\zeta,n)$.
Lemma \ref{lemma-sc4} shows
\begin{align}
\label{fgg-1}
{\mathscr H}(\tilde Q_0(n))+\theta(\zeta,n)=
{\mathscr H}(Q(\varphi,\Psi_{n}(\varphi)))+\theta(\varphi,\Psi_{n}(\varphi)),
\quad \zeta<\varphi<\varphi_n.
\end{align}
The same proof as Lemma \ref{gg-prop2} gives
\begin{align}
\label{gg-17-4}
Q(\varphi,\Psi_n(\varphi))\le \tilde Q_0(n),\quad\zeta<\varphi<\varphi_n.
\end{align}
It follows from \eqref{fgg-1}, \eqref{gg-17-4} and \eqref{gg-17-2} that
\begin{align*}
\theta(\varphi,\Psi_{n}(\varphi))\le \theta(\zeta,n)<\theta(\zeta,m),
\quad \zeta<\varphi<\varphi_n,
\end{align*}
which, together with \eqref{gg-17-4} and \eqref{gg-2}, leads to
the second estimate in \eqref{gg-17-1}.
Furthermore, the second limit in \eqref{gg-17} follows from \eqref{fgg-1},
\eqref{gg-17-4} and \eqref{gg-2}.
Finally, \eqref{gg-17} shows $q,\theta\in C([\zeta,+\infty)\times[0,m])$.
$\hfill\Box$\vskip 4mm

In the physical plane, Proposition \ref{gg-thm1} is stated as follows.

\begin{proposition}
\label{phys-thm71}
Assume that $\tilde\Upsilon\in C^{2}([0,y_0))\cap C^{1}([0,y_0])$
and $\tilde q_0\in C^{1}([0,y_0))\cap C([0,y_0])$ satisfy
\eqref{gg-phy1} and \eqref{gg-phy2}.
Then the problem
\begin{align*}
&\mbox{div}(\rho(|\nabla\varphi|^2)\nabla\varphi)=0,\quad&&(x,y)\in\tilde\Omega,
\\
&\varphi(\tilde\Upsilon(y),y)=\zeta,\quad&&0<y<y_0,
\\
&|\nabla\varphi(\tilde\Upsilon(y),y)|=\tilde q_0(y),\quad&&0<y<y_0,
\\
&\pd{\varphi}y(x,0)=0,\quad&&x>\tilde\Upsilon(0)
\end{align*}
admits a unique solution $\varphi\in C^{2}(\overline{\tilde\Omega}
\setminus\tilde\Gamma_{\text{\rm up}})\cap C^{1}(\overline{\tilde\Omega})$,
where
$$
\tilde\Gamma_{\text{\rm up}}=\big\{(x,y):x\ge\tilde\Upsilon(y_0),
y=y_0-\tilde\Upsilon'(y_0)(x-\tilde\Upsilon(y_0))\big\}
$$
and $\tilde\Omega$ is the domain bounded by
$\{(\tilde\Upsilon(y),y):0\le y\le y_0\}$,
$\tilde\Gamma_{\text{\rm up}}$ and the $x$-axis.
Furthermore, $\varphi$ satisfies
\begin{gather*}
\inf_{(0,y_0)}\tilde q_0\le |\nabla\varphi(x,y)|<c^*,\quad
0\le \pd\varphi y(x,y)<-\tilde\Upsilon'(y_0)\pd\varphi x(x,y),
\quad (x,y)\in \overline{\tilde\Omega}\setminus\tilde\Gamma_{\text{\rm up}},
\\
|\nabla\varphi(x,y)|=c^*,\quad
\pd\varphi x(x,y)=\frac{c^*}{\sqrt{1+(\tilde\Upsilon'(y_0))^2}},\quad
\pd\varphi y(x,y)=\frac{-\tilde\Upsilon'(y_0)c^*}{\sqrt{1+(\tilde\Upsilon'(y_0))^2}},\quad
(x,y)\in\tilde\Gamma_{\text{\rm up}},
\\
\varphi(x,y)=\zeta+\frac{c^*}{\sqrt{1+(\tilde\Upsilon'(y_0))^2}}
\big((x-\tilde\Upsilon(y_0))-\tilde\Upsilon'(y_0)(y-y_0)\big),\quad
(x,y)\in\tilde\Gamma_{\text{\rm up}}.
\end{gather*}
\end{proposition}

Below we study the global behavior of positive characteristics
and the lower bound estimate of the flow speed to the problem
\eqref{gg-13}--\eqref{gg-16}.

\begin{proposition}
\label{gg-thm1-ch}
Assume that $Q\in C^{1}([\zeta,+\infty)\times[0,m))$ is the solution to the problem
\eqref{gg-13}--\eqref{gg-16}, where
$\tilde Q_0\in C^{1}([0,m))$ and $\tilde G_0\in C([0,m))\cap L^1(0,m)$
satisfying \eqref{gg-1} and \eqref{gg-2}. Then,
each positive characteristic never touches $(\zeta,+\infty)\times\{m\}$.
Furthermore, if, in addition, $\inf_{(0,m)}(-\tilde G_0-b^{1/2}(\tilde Q_0)|\tilde Q'_0|)>0$,
then
\begin{align}
\label{sd-thm3-re2}
Q(\varphi,\psi)\le-M (\varphi+1)^{2/(\gamma+1)},\quad
(\varphi,\psi)\in(\zeta,+\infty)\times(0,m),
\end{align}
where $M>0$ is a constant depending only on $\gamma$, $m$,
$\sup_{(0,m)}\tilde Q_0$
and $\inf_{(0,m)}(-\tilde G_0-b^{1/2}(\tilde Q_0)|\tilde Q'_0|)$.
\end{proposition}

\Proof
For $(\varphi,\psi)\in[\zeta,+\infty)\times[0,m)$, set
\begin{align*}
W(\varphi,\psi)=Q_\varphi(\varphi,\psi)-{b^{1/2}(Q(\varphi,\psi))}Q_\psi(\varphi,\psi),\quad
Z(\varphi,\psi)=-Q_\varphi(\varphi,\psi)-{b^{1/2}(Q(\varphi,\psi))}Q_\psi(\varphi,\psi).
\end{align*}
Then, $(W,Z)\in C([\zeta,+\infty)\times[0,m))\times C([\zeta,+\infty)\times[0,m))$
solves
\begin{align}
\label{lll-3}
&W_\varphi+b^{1/2}(Q)W_\psi
=\frac14b^{-1}(Q)p(Q)W(W+Z),
\quad&& (\varphi,\psi)\in(\zeta,+\infty)\times(0,m),
\\
\label{lll-4}
&Z_\varphi-b^{1/2}(Q)Z_\psi
=-\frac14b^{-1}(Q)p(Q)Z(W+Z),
\quad&&(\varphi,\psi)\in(\zeta,+\infty)\times(0,m).
\end{align}
It follows from \eqref{gg-17-0} that
\begin{align}
\label{lll-1}
W(\varphi,\psi)\le 0,\quad Z(\varphi,\psi)\ge 0,\quad(\varphi,\psi)\in[\zeta,+\infty)\times[0,m).
\end{align}
First, consider the case that $m\in\mbox{supp }W(\zeta,\cdot)$.
Note that
\begin{align*}
\frac{d}{d\psi}\big({\mathscr H}(Q(\zeta,\psi))-\theta(\zeta,\psi)\big)=
Q_\varphi(\zeta,\psi)-b^{1/2}(Q(\zeta,\psi))
Q_\psi(\zeta,\psi)=W(\zeta,\psi),\quad0<\psi<m,
\end{align*}
which, together with \eqref{lll-1}, $m\in\mbox{supp }W(\zeta,\cdot)$ and \eqref{gg-17}, leads to
\begin{align*}
{\mathscr H}(Q(\zeta,\psi))-\theta(\zeta,\psi)>{\mathscr H}(-\infty)-\theta(\zeta,m),\quad0<\psi<m.
\end{align*}
Then, Lemma \ref{lemma-sc4} and \eqref{gg-17} show that for each $\psi\in(0,m)$,
the positive characteristic from $(\zeta,\psi)$ never touches $(\zeta,+\infty)\times\{m\}$.

\vskip10mm

\hskip35mm
\setlength{\unitlength}{0.6mm}
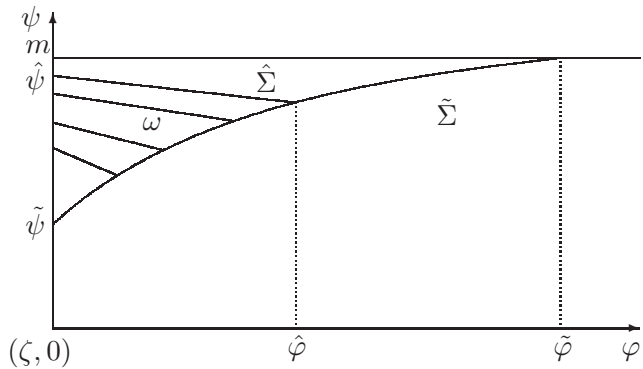
\begin{picture}(250,65)
\put(0,0){\vector(1,0){130}}
\put(0,0){\vector(0,1){70}}
\put(126,-6){$\varphi$} \put(-7,67){$\psi$}
\put(0,60){\qbezier(0,0)(65,0)(130,0)}

\put(0,23){\cbezier(0,0)(30,28)(80,33)(113,37)}

\put(0,40){\qbezier(0,0)(7,-3)(14,-6)}

\put(0,45.6){\qbezier(0,0)(12,-3)(24,-6)}

\put(0,52){\qbezier(0,0)(20,-3)(40,-6)}

\put(0,56){\qbezier(0,0)(27,-3)(54,-6)}

\put(54,0){\qbezier[50](0,0)(0,25)(0,50)}

\put(52,-6){$\hat\varphi$}

\put(20,44){$\omega$}

\put(45,52.5){$\hat\Sigma$}

\put(85,45){$\tilde\Sigma$}

\put(-6,61){$m$}

\put(-6,53){$\hat\psi$}

\put(-6,21){$\tilde\psi$}

\put(-10,-7){$(\zeta,0)$}

\put(112.5,0){\qbezier[60](0,0)(0,30)(0,60)}

\put(111,-6){$\tilde\varphi$}

\end{picture}

\vskip10mm

Next, consider the case that $m\not\in\mbox{supp }W(\zeta,\cdot)$,
which, together with \eqref{lll-1}, implies that
there exists $\psi_0\in(0,m)$ such that
\begin{align}
\label{lll-2}
W(\zeta,\psi)=0,\quad \psi_0<\psi<m.
\end{align}
It can be proved by a contradiction argument
that each positive characteristic never touches $(\zeta,+\infty)\times\{m\}$.
Otherwise, there exists $\tilde\psi\in(\psi_0,m)$ such that
the positive characteristic $\tilde\Sigma$ from $(\zeta,\tilde\psi)$
touches $(\zeta,+\infty)\times\{m\}$.
That is to say,
$$
\tilde\Sigma:
\tilde\Psi'(\varphi)=b^{1/2}(Q(\varphi,\tilde\Psi(\varphi))),
\quad
\tilde\psi<\tilde\Psi(\varphi)<m,\quad\zeta<\varphi<\tilde\varphi,
\quad\tilde\Psi(\zeta)=\tilde\psi,\quad\tilde\Psi(\tilde\varphi)=m.
$$
Let $\omega$ be the bounded domain bounded by $\varphi=\zeta$, $\psi=m$ and $\tilde\Sigma$.
Then, \eqref{lll-3} and \eqref{lll-2} yield that
\begin{align}
\label{lll-5}
W(\varphi,\psi)=0,\quad (\varphi,\psi)\in\omega.
\end{align}
For each $\hat\psi\in(\tilde\psi,m)$, let
$$
\hat\Sigma:
\hat\Psi'(\varphi)=-b^{1/2}(Q(\varphi,\hat\Psi(\varphi))),
\quad\zeta<\varphi<\hat\varphi,
\quad\hat\Psi(\zeta)=\hat\psi,\quad\hat\Psi(\hat\varphi)=\tilde\Psi(\hat\varphi)
$$
be the negative characteristic from $(\zeta,\hat\psi)$ to $\tilde\Sigma$.
It follows from \eqref{lll-5} that
$$
\frac{d}{d\varphi}Q(\varphi,\hat\Psi(\varphi))=W(\varphi,\hat\Psi(\varphi))=0,\quad\zeta<\varphi<\hat\varphi,
$$
which yields that $\hat\Sigma$ is a straight line and
\begin{align}
\label{lll-7}
Q(\varphi,\hat\Psi(\varphi)))=Q(\hat\varphi,\hat\Psi(\hat\varphi)),\quad\zeta<\varphi<\hat\varphi.
\end{align}
By \eqref{lll-4}, \eqref{lll-5}, \eqref{lll-7} and \eqref{a8-8},
there exists a positive constant $M$ depending only
on $\gamma$ such that
\begin{align*}
\frac{d}{d\varphi}Z(\varphi,\hat\Psi(\varphi))
=&-\frac14b^{-1}(Q(\varphi,\hat\Psi(\varphi)))p(Q(\varphi,\hat\Psi(\varphi)))Z^2(\varphi,\hat\Psi(\varphi))
\\
\le& M Q^{-1}(\varphi,\hat\Psi(\varphi)) Z^2(\varphi,\hat\Psi(\varphi))
= M Q^{-1}(\hat\varphi,\hat\Psi(\hat\varphi)) Z^2(\varphi,\hat\Psi(\varphi)),
\quad\zeta<\varphi<\hat\varphi,
\end{align*}
which, together with \eqref{lll-1} and $\hat\Psi(\hat\varphi)=\tilde\Psi(\hat\varphi)$, leads to
\begin{align*}
Z(\hat\varphi,\tilde\Psi(\hat\varphi))
\le-\frac{Q(\hat\varphi,\tilde\Psi(\hat\varphi))}{M(\hat\varphi-\zeta)}.
\end{align*}
By \eqref{gg-17},
$\hat\varphi$ will take all values over $(\zeta,\tilde\varphi)$
when $\hat\psi$ varies from $\tilde\psi$ to $m$.
Therefore,
\begin{align*}
\frac{d}{d\varphi}Q(\varphi,\tilde\Psi(\varphi))
=-Z(\varphi,\tilde\Psi(\varphi))
\ge \frac{2}{M(\tilde\varphi-\zeta)}
Q(\varphi,\tilde\Psi(\varphi)),\quad\frac{\zeta+\tilde\varphi}2<\varphi<\tilde\varphi,
\end{align*}
which contradicts $\tilde\Psi(\tilde\varphi)=m$ and \eqref{gg-17}.

Summing up, one gets that
each positive characteristic never touches $(\zeta,+\infty)\times\{m\}$.
Then, \eqref{sd-thm3-re2} can be proved
by a similar discussion as {Step III} in the proof of Proposition \ref{sss-thm1-lb}.
$\hfill\Box$\vskip 4mm

\subsection{Global smooth supersonic flows}

\vskip5mm

\hskip40mm
\setlength{\unitlength}{0.6mm}
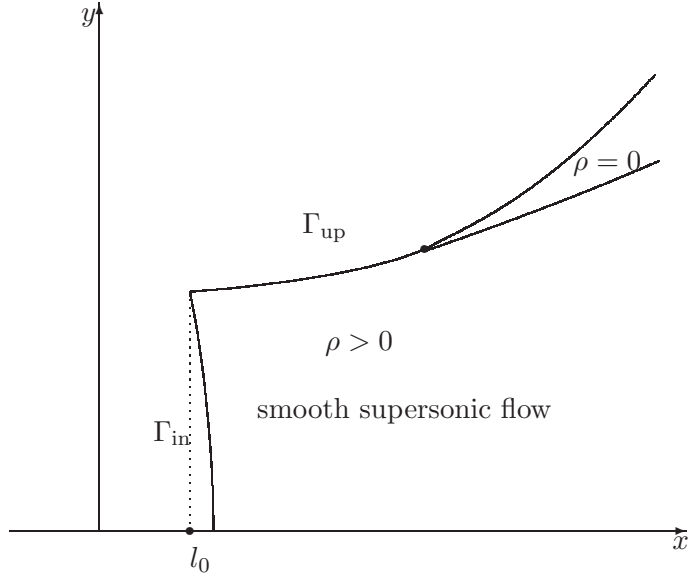
\begin{picture}(250,120)
\put(-20,0){\vector(1,0){150}}
\put(0,0){\vector(0,1){117}}
\put(127,-4){$x$} \put(-4,113){$y$}

\put(60,48){\cbezier(-40,5)(-20,6.5)(5,10)(15,16)}
\put(60,48){\cbezier(15,16)(25,21)(40,28)(63,53)}

\put(60,0){\cbezier(-40,53)(-37.5,42)(-34.8,20)(-34.7,0)}

\put(45,66){$\Gamma_{\text{\rm up}}$}
\put(12,20){$\Gamma_{\text{\rm in}}$}

\put(20,0){\qbezier[30](0,0)(0,26.5)(0,53)}
\put(20,0){\circle*{1.5}}
\put(20,-8){$l_0$}
\put(35,25){smooth supersonic flow}
\put(50,40){$\rho>0$}

\put(105,80){$\rho=0$}
\put(72,62.5){\circle*{1.5}}
\put(75,55){\qbezier(-3,7)(26,17)(49,27)}

\end{picture}
\vskip10mm

\begin{theorem}
\label{phys-thm81}
Assume that $f\in C^{2}([l_0,l_1))$ satisfies \eqref{sss-0},
$\Upsilon\in C^{2}([0,f(l_0)])$ satisfies \eqref{xc4} and $q_0\in C^{1}([0,f(l_0)])$
satisfies \eqref{phy-8-q}--\eqref{phy-9}.
Then the problem \eqref{phy1}--\eqref{phy5} admits uniquely a global smooth solution
$\varphi\in C^{1}(\overline\Omega)$ with $|\nabla\varphi|\in C^{0,1}(\overline\Omega)$.
Moreover,
the global smooth supersonic flow belongs to one and only one of the following two cases:
\\

{\rm\bf Case I} Global smooth supersonic flow without vacuum.
In this case, $\varphi\in C^{2}(\overline\Omega)$ and
$c_1\le |\nabla\varphi|<c^*$ on $\overline\Omega$.
If $f''(1+(f')^2)^{-3/2}\in L^\infty(l_0,l_1)$ and
$q_0$ satisfies \eqref{phy-8-positive} additionally, then $|\nabla\varphi|$ is globally Lipschitz continuous in $\Omega$.

{\rm\bf Case II} Global smooth supersonic flow with vacuum.
In this case, set
$$
x_0=\sup\big\{l\in(l_0,l_1):|\nabla\varphi(x,f(x))|<c^*\mbox{ for each }l_0<x<l\big\}
$$
and
$$
\Omega_{\text{\rm v}}=\big\{(x,y)\in\Omega:x>x_0,y>f(x_0)+f'(x_0)(x-x_0)\big\}.
$$
Then

{\rm (i)} $\overline\Omega_{\text{\rm v}}$ is the set of vacuum points.
Moreover, $\varphi\in C^{2}(\overline\Omega\setminus\overline\Omega_{\text{\rm v}})$ satisfying
\begin{gather*}
c_1\le |\nabla\varphi(x,y)|<c^*,\quad
0\le \pd\varphi y(x,y)<f'(x_0)\pd\varphi x(x,y),
\quad
(x,y)\in\overline\Omega\setminus\overline\Omega_{\text{\rm v}},
\\
\pd\varphi x(x,y)=\frac{c^*}{\sqrt{1+(f'(x_0))^2}},\quad
\pd\varphi y(x,y)=\frac{f'(x_0)c^*}{\sqrt{1+(f'(x_0))^2}},
\quad(x,y)\in\partial\Omega_{\text{\rm v}}\cap\Omega.
\end{gather*}

{\rm (ii)}
$|\nabla\varphi(x,f(x))|=c^*-O\big((x_0-x)^2\big)$ as $x\to x_0^-$.

{\rm (iii)} $|\nabla\varphi|$ is globally Lipschitz continuous in $\Omega$.

{\rm (iv)} $|\nabla\varphi|\in C^1(\overline\Omega\setminus{(x_0,f(x_0))})$
and
for each $(\tilde x,\tilde y)\in\partial\Omega_{\text{\rm v}}\cap\Omega$,
\begin{align*}
{\lim_{\stackrel{(x,y)\to(\tilde x,\tilde y)}
{(x,y)\in\Omega\setminus\overline\Omega_{\text{\rm v}}}}}
\nabla|\nabla\varphi(x,y)|=(0,0),
\quad
{\lim_{\stackrel{(x,y)\to(\tilde x,\tilde y)}
{(x,y)\in\Omega\setminus\overline\Omega_{\text{\rm v}}}}}
(c^*-|\nabla\varphi(x,y)|)^{1/2}\nabla\arctan\frac{\varphi_y(x,y)}{\varphi_x(x,y)}=(0,0).
\end{align*}
In particular, if $1<\gamma\le2$, then
$\rho(|\nabla\varphi|^2)\nabla\varphi\in C^1(\overline\Omega\setminus{(x_0,f(x_0))})$.

{\rm (v)} If $1<\gamma<3$ and $f''(x_0)>0$ additionally,
then for any $\varepsilon>0$,
\begin{align*}
{\lim_{\stackrel{(x,y)\to(\tilde x,\tilde y)}
{(x,y)\in\Omega\setminus\overline\Omega_{\text{\rm v}}}}}
(c^*-|\nabla\varphi(x,y)|)^{-(\gamma+1+\varepsilon)/(4\gamma-4)}
\nabla |\nabla\varphi(x,y)|\cdot(-1,f'(x_0))=+\infty,
\quad
(\tilde x,\tilde y)\in\partial\Omega_{\text{\rm v}}\cap\Omega.
\end{align*}
\end{theorem}

\Proof
If there is no vacuum for the flow, then the theorem is proved by
Theorems \ref{phys-thm31} and \ref{phys-thm53}.
Below we consider the case that there is a vacuum for the flow.
It follows from Theorem \ref{phys-thm31} and Proposition \ref{phys-thm71} that
$\overline\Omega_{\text{\rm v}}$ is the set of vacuum points
for a global smooth supersonic flow to the problem \eqref{phy1}--\eqref{phy5},
and the problem \eqref{phy1}--\eqref{phy5} in the gas region
admits uniquely a smooth solution
$\varphi\in C^{2}(\overline\Omega\setminus\overline\Omega_{\text{\rm v}})
\cap C^{1}(\overline{\Omega\setminus\Omega_{\text{\rm v}}})$.
By Theorem \ref{phys-thm31} and Proposition \ref{phys-thm71},
$\varphi$ satisfies
\begin{align*}
&\mbox{div}(\rho(|\nabla\varphi|^2)\nabla\varphi)=0,\quad&&(x,y)\in\hat\Omega,
\\
&\varphi(\Upsilon(y),y)=0,\quad&&0<y<f(l_0),
\\
&|\nabla\varphi(\Upsilon(y),y)|=q_0(y),\quad&&0<y<f(l_0),
\\
&\pd{\varphi}y(x,0)=0,\quad&&\Upsilon(0)<x<l_1,
\\
&\pd{\varphi}y(x,f(x))-f'(x)\pd{\varphi}x(x,f(x))=0,
\quad&&l_0<x<x_0,
\\
&|\nabla\varphi(x,\tilde f(x))|=c^*,
\quad&&x>x_0,
\end{align*}
where
\begin{gather*}
\tilde f(x)=f(x_0)+f'(x_0)(x-x_0),\quad x\ge x_0
\end{gather*}
and $\hat\Omega$ is the unbounded domain bounded by
$\Gamma_{\text{\rm in}}$, the $x$-axis, $y=f(x)\,(l_0\le x\le x_0)$
and $y=\tilde f(x)\,(x\ge x_0)$.
Furthermore,
\begin{align*}
\pd\varphi x(x,\tilde f(x))=\frac{c^*}{\sqrt{1+(f'(x_0))^2}},\quad
\pd\varphi y(x,\tilde f(x))=\frac{f'(x_0)c^*}{\sqrt{1+(f'(x_0))^2}},
\quad x\ge x_0
\end{align*}
and
\begin{align}
\label{boun}
\varphi(x,\tilde f(x))=\varphi(x_0,f(x_0))+\frac{c^*}{\sqrt{1+(f'(x_0))^2}}
\big((x-x_0+f'(x_0)(\tilde f(x)-f(x_0))\big),\quad x\ge x_0.
\end{align}
Extend $\varphi$ from $\overline{\hat\Omega}$ to $\overline\Omega$
by defining
\begin{align*}
\varphi(x,y)=\varphi(x_0,f(x_0))+\frac{c^*}{\sqrt{1+(f'(x_0))^2}}
\big((x-x_0+f'(x_0)(y-f(x_0))\big),\quad
(x,y)\in\overline\Omega\setminus\overline{\hat\Omega}.
\end{align*}
Then, $\varphi\in C^{2}(\overline\Omega\setminus\overline\Omega_{\text{\rm v}})\cap C^{1}(\overline\Omega)$
is just a solution to the problem \eqref{phy1}--\eqref{phy5} satisfying (i)
according to Theorem \ref{phys-thm31} and Proposition \ref{phys-thm71}.
Moreover, (ii) follows from Theorem \ref{sss-thm2}.

To verify the uniqueness of the global smooth solution to the problem \eqref{phy1}--\eqref{phy5},
it suffices to prove that if $\hat\varphi\in C^1(\overline\Omega_{\text{\rm v}})$ solves
\begin{align*}
&\mbox{div}(\rho(|\nabla\hat\varphi|^2)\nabla\hat\varphi)=0,\quad&&(x,y)\in\Omega_{\text{\rm v}},
\\
&\pd{\hat\varphi} x(x,\tilde f(x))=\frac{c^*}{\sqrt{1+(f'(x_0))^2}},\quad
\pd{\hat\varphi} y(x,\tilde f(x))=\frac{f'(x_0)c^*}{\sqrt{1+(f'(x_0))^2}},\quad&&x>x_0,
\\
&\rho(|\nabla\hat\varphi(x,f(x))|^2)\Big(\pd{\hat\varphi}y(x,f(x))-f'(x)\pd{\hat\varphi}x(x,f(x))\Big)=0,
\quad&&x_0<x<l_1,
\end{align*}
then $|\nabla\hat\varphi|=c^*$ in $\Omega_{\text{\rm v}}$.
Note that the flow angle on $\partial\Omega_{\text{\rm v}}\cap\Omega$ equals $\arctan f'(x_0)\in(0,\pi/2)$.
A similar argument as the beginning of the proof of Theorem \ref{asd} shows
\begin{align*}
|\nabla\hat\varphi(x,y)|=c^*,
\quad(x,y)\in\overline\Omega_{\text{\rm v}},\,x-x_0\le\delta
\end{align*}
with some positive constant $\delta$. Furthermore, one can prove that
\begin{align*}
\pd{\hat\varphi} x(x,y)=\frac{c^*}{\sqrt{1+(f'(x_0))^2}},\quad
\pd{\hat\varphi} y(x,y)=\frac{f'(x_0)c^*}{\sqrt{1+(f'(x_0))^2}},
\quad(x,y)\in\overline\Omega_{\text{\rm v}},\,x-x_0\le\delta.
\end{align*}
Repeating this discussion leads to
\begin{align*}
|\nabla\hat\varphi(x,y)|=c^*,
\quad\pd{\hat\varphi} x(x,y)=\frac{c^*}{\sqrt{1+(f'(x_0))^2}},\quad
\pd{\hat\varphi} y(x,y)=\frac{f'(x_0)c^*}{\sqrt{1+(f'(x_0))^2}},
\quad(x,y)\in\overline\Omega_{\text{\rm v}}.
\end{align*}

We now prove (iii)--(v).
Transform the supersonic flow on $\overline{\hat\Omega}$ into the potential plane.
This transformation is one to one according to \eqref{boun}.
It is assumed that $(\Upsilon(0),0)$ and $(x_0,f(x_0))$
are transformed into $(0,0)$ and $(\zeta,m)$, respectively,
in the coordinates transformation
from the physical plane to the potential plane.
For $(\varphi,\psi)\in[0,+\infty)\times[0,m]\setminus[\zeta,+\infty)\times\{m\}$,
set
\begin{gather*}
W(\varphi,\psi)=Q_\varphi(\varphi,\psi)-b^{1/2}(Q(\varphi,\psi))Q_\psi(\varphi,\psi),
\quad
Z(\varphi,\psi)=-Q_\varphi(\varphi,\psi)-b^{1/2}(Q(\varphi,\psi))Q_\psi(\varphi,\psi),
\\
{\mathscr U}(\varphi,\psi)
=b^{1/2}(Q(\varphi,\psi)W(\varphi,\psi),
\quad
{\mathscr V}(\varphi,\psi)
=b^{1/2}(Q(\varphi,\psi))Z(\varphi,\psi).
\end{gather*}
Then, $(W,Z,Q)$ satisfies
\begin{align}
\label{a8-6}
W(\varphi,\psi)\le0,\quad
Z(\varphi,\psi)\ge0,\quad
Q(\varphi,\psi)\le-\mu_1,\quad
(\varphi,\psi)\in(0,+\infty)\times(0,m),
\end{align}
$(W,Z)\in C([0,+\infty)\times[0,m]\setminus[\zeta,+\infty)\times\{m\})\times
C([0,+\infty)\times[0,m]\setminus[\zeta,+\infty)\times\{m\})$ solving
\begin{align}
\label{wz-1}
&W_\varphi+b^{1/2}(Q)W_\psi
=\frac14b^{-1}(Q)p(Q)W(W+Z),
\quad&&(\varphi,\psi)\in(0,+\infty)\times(0,m),
\\
\label{wz-2}
&Z_\varphi-b^{1/2}(Q)Z_\psi
=-\frac14b^{-1}(Q)p(Q)Z(W+Z),
\quad&&(\varphi,\psi)\in(0,+\infty)\times(0,m),
\end{align}
and $({\mathscr U},{\mathscr V})\in C([0,+\infty)\times[0,m]\setminus[\zeta,+\infty)\times\{m\})\times
C([0,+\infty)\times[0,m]\setminus[\zeta,+\infty)\times\{m\})$ solving
\begin{align}
\label{ccv1}
&{\mathscr U}_\varphi+b^{1/2}(Q){\mathscr U}_\psi
=\frac14b^{-3/2}(Q)p(Q){\mathscr U}({\mathscr U}-{\mathscr V})\ge0,
\quad&&(\varphi,\psi)\in(0,+\infty)\times(0,m),
\\
\label{ccv2}
&{\mathscr V}_\varphi-b^{1/2}(Q){\mathscr V}_\psi
=\frac14b^{-3/2}(Q)p(Q){\mathscr V}({\mathscr U}-{\mathscr V})\le0,
\quad&&(\varphi,\psi)\in(0,+\infty)\times(0,m),
\\
&{\mathscr U}(0,\psi)={\mathscr U}_0(\psi)=
b^{1/2}(Q_0(\psi))G_0(\psi)-b(Q_0(\psi))Q'_0(\psi),\quad&&\psi\in(0,m),
\nonumber
\\
&{\mathscr V}(0,\psi)={\mathscr V}_0(\psi)=
-b^{1/2}(Q_0(\psi))G_0(\psi)-b(Q_0(\psi))Q'_0(\psi),\quad&&\psi\in(0,m),
\nonumber
\\
&{\mathscr U}(\varphi,0)+{\mathscr V}(\varphi,0)=0,\quad&&\varphi\in(0,+\infty),
\nonumber
\\
&{\mathscr U}(\varphi,m)+{\mathscr V}(\varphi,m)=
\frac{2f''(X_{\text{\rm up}}(\varphi))(1+(f'(X_{\text{\rm up}}(\varphi)))^2)^{-3/2}}
{A^{-1}_+(Q(\varphi,m))},
\quad&&\varphi\in(0,\zeta),
\nonumber
\\
&\lim_{\psi\to m^-}Q(\varphi,\psi)=-\infty,\quad&&\varphi\in(\zeta,+\infty).
\nonumber
\end{align}
By the method of characteristics and Proposition \ref{gg-thm1-ch}, one can prove
\begin{align*}
-M_1\le{\mathscr U}(\varphi,\psi)\le 0,\quad
0\le {\mathscr V}(\varphi,\psi)\le M_1,\quad
(\varphi,\psi)\in(0,+\infty)\times(0,m)
\end{align*}
with $M_1=M_1(\gamma,m,\mu_1,\mu_3,\zeta,f)$ being a positive constant.
Therefore,
\begin{align*}
-M_1\le b^{1/2}(Q(\varphi,\psi))Q_\varphi(\varphi,\psi)\le 0,\quad
|b(Q(\varphi,\psi))Q_\psi(\varphi,\psi)|\le M_1,\quad
(\varphi,\psi)\in(0,+\infty)\times(0,m),
\end{align*}
which yields that $|\nabla\varphi|$ is globally Lipschitz continuous in
$\Omega\setminus\overline\Omega_{\text{\rm v}}$.
Thus (iii) is proved since $|\nabla\varphi|=c^*$ on $\overline\Omega_{\text{\rm v}}$.

\vskip10mm

\hskip30mm
\setlength{\unitlength}{0.6mm}
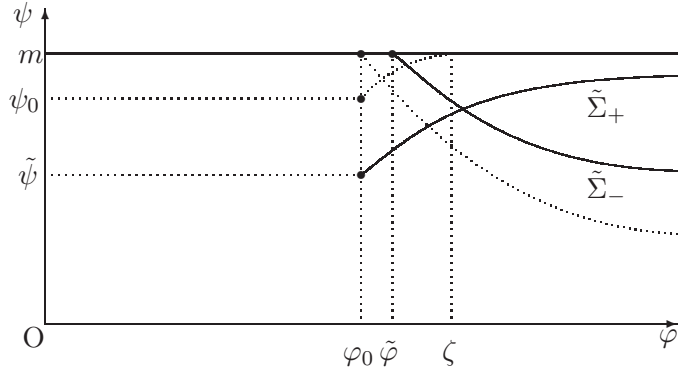
\begin{picture}(250,65)
\put(0,0){\vector(1,0){140}}
\put(0,0){\vector(0,1){70}}
\put(136,-4){$\varphi$} \put(-7,67){$\psi$}
\put(0,60){\qbezier(0,0)(55,0)(140,0)}

\put(90,0){\qbezier[40](0,0)(0,30)(0,60)}

\put(-6,58){$m$}

\put(-5,-5){O}

\put(70,0){\qbezier[40](0,0)(0,30)(0,60)}
\put(66,-8){$\varphi_0$}

\put(88,-8){$\zeta$}

\put(77,60){\cbezier(0,0)(10,-8)(25,-25)(63,-26)}

\put(77,60){\circle*{1.5}}

\put(120,28){$\tilde\Sigma_-$}

\put(77,0){\qbezier[40](0,0)(0,30)(0,60)}
\put(74,-8){$\tilde\varphi$}

\put(0,33){\qbezier[45](0,0)(35,0)(70,0)}

\put(-6,31){$\tilde\psi$}

\put(70,33){\cbezier(0,0)(15,12)(25,21)(70,22)}

\put(70,33){\circle*{1.5}}

\put(120,46){$\tilde\Sigma_+$}

\put(70,60){\circle*{1.5}}

\put(70,60){\qbezier[60](0,0)(30,-38)(70,-40)}

\put(0,50){\qbezier[45](0,0)(35,0)(70,0)}

\put(-8,48){$\psi_0$}

\put(70,50){\circle*{1.5}}

\put(70,50){\qbezier[20](0,0)(8,9)(20,9.8)}

\end{picture}

\vskip10mm

By Lemma \ref{lemma-sc4} and
Proposition \ref{gg-thm1-ch},
there exists $\varphi_0\in(0,\zeta)$ such that
both the negative characteristic from $(\varphi_0,m)$
and the positive characteristic from $(\varphi_0,m/2)$ are unbounded.
Denote
\begin{align*}
\psi_0=\sup\big\{\psi\in(0,m):\mbox{ the positive characteristic from }(\varphi_0,\psi)
\mbox{ is unbounded}\big\}.
\end{align*}
Fix $\tilde\varphi\in(\varphi_0,\zeta)$ and $\tilde\psi\in(0,\psi_0)$.
Let $\tilde\Sigma_+$ be the positive characteristic from $(\varphi_0,\tilde\psi)$
and $\tilde\Sigma_-$ be the negative characteristic from $(\tilde\varphi,m)$, i.e.,
\begin{align*}
&\tilde\Sigma_+:
\tilde\Psi'_+(\varphi)=b^{1/2}(Q(\varphi,\tilde\Psi_+(\varphi))),
\quad\tilde\psi<\tilde\Psi_+(\varphi)<m,
\quad\varphi>\varphi_0,
\quad\tilde\Psi_+(\varphi_0)=\tilde\psi,
\\
&\tilde\Sigma_-:
\tilde\Psi'_-(\varphi)=-b^{1/2}(Q(\varphi,\tilde\Psi_-(\varphi))),
\quad0<\tilde\Psi_-(\varphi)<m,
\quad\varphi>\tilde\varphi,
\quad\tilde\Psi_-(\tilde\varphi)=m.
\end{align*}
It follows from \eqref{a8-6}, \eqref{ccv1} and \eqref{ccv2} that
\begin{align*}
&\frac{d}{d\varphi}{\mathscr U}(\varphi,\tilde\Psi_+(\varphi))
\ge M_2(-Q(\varphi,\tilde\Psi_+(\varphi)))^{(\gamma-1)/2}{\mathscr U}^2(\varphi,\tilde\Psi_+(\varphi)),
\quad&&\varphi>\varphi_0,
\\
&\frac{d}{d\varphi}{\mathscr V}(\varphi,\tilde\Psi_-(\varphi))
\le-M_2(-Q(\varphi,\tilde\Psi_-(\varphi)))^{(\gamma-1)/2}{\mathscr V}^2(\varphi,\tilde\Psi_-(\varphi)),
\quad&&\varphi>\tilde\varphi
\end{align*}
with $M_2=M_2(\gamma,\mu_1)$ being a positive constant.
Therefore,
\begin{align}
\label{ccv3}
&{\mathscr U}(\varphi,\tilde\Psi_+(\varphi))\ge
-\frac1{M_2}\Big(\int_{\varphi_0}^\varphi
(-Q(s,\tilde\Psi_+(s)))^{(\gamma-1)/2}ds\Big)^{-1},
\quad&&\varphi>\varphi_0,
\\
\label{ccv4}
&{\mathscr V}(\varphi,\tilde\Psi_-(\varphi))\le
\frac1{M_2}\Big(\int_{\tilde\varphi}^\varphi
(-Q(s,\tilde\Psi_-(s)))^{(\gamma-1)/2}ds\Big)^{-1},
\quad&&\varphi>\tilde\varphi.
\end{align}
Letting $\tilde\psi\to\psi_0^-$ in \eqref{ccv3}
and $\tilde\varphi\to\zeta^-$ in \eqref{ccv4},
and using (i) and \eqref{a8-6}, one can get
\begin{align*}
\lim_{{\psi\to m^-}}\min_{\zeta_1\le\varphi\le\zeta_2}{\mathscr U}(\varphi,\psi)=0,
\quad
\lim_{{\psi\to m^-}}\max_{\zeta_1\le\varphi\le\zeta_2}{\mathscr V}(\varphi,\psi)=0,
\quad
\zeta<\zeta_1<\zeta_2.
\end{align*}
Therefore,
\begin{align*}
\lim_{{\psi\to m^-}}\min_{\zeta_1\le\varphi\le\zeta_2}
b^{1/2}(Q(\varphi,\psi))Q_\varphi(\varphi,\psi)=0,
\quad
\lim_{{\psi\to m^-}}\max_{\zeta_1\le\varphi\le\zeta_2}|b(Q(\varphi,\psi))Q_\psi(\varphi,\psi)|=0,
\quad
\zeta<\zeta_1<\zeta_2,
\end{align*}
which yield that for each $(\tilde x,\tilde y)\in\partial\Omega_{\text{\rm v}}\cap\Omega$,
\begin{align*}
{\lim_{\stackrel{(x,y)\to(\tilde x,\tilde y)}
{(x,y)\in\Omega\setminus\overline\Omega_{\text{\rm v}}}}}
\nabla|\nabla\varphi(x,y)|=(0,0),
\quad
{\lim_{\stackrel{(x,y)\to(\tilde x,\tilde y)}
{(x,y)\in\Omega\setminus\overline\Omega_{\text{\rm v}}}}}
(c^*-|\nabla\varphi(x,y)|)^{1/2}\nabla\arctan\frac{\varphi_y(x,y)}{\varphi_x(x,y)}=(0,0).
\end{align*}
Thus (iv) is proved since $|\nabla\varphi|=c^*$ on $\overline\Omega_{\text{\rm v}}$.

Finally, we prove (v).
It follows from $f''(x_0)>0$ and \eqref{a8-6} that
\begin{align}
\label{a8-6-1}
\lim_{\overline{\varphi\to\zeta^-}}b^{1/2}(Q(\varphi,m))Z(\varphi,m)
\ge&\lim_{{\varphi\to\zeta^-}}b^{1/2}(Q(\varphi,m))(W(\varphi,m)+Z(\varphi,m))
\nonumber
\\
=&\frac{2f''(x_0)}{c^*(1+(f'(x_0))^2)^{3/2}}>0,
\end{align}
which yields that $Z(\cdot,m)$ is positive in the left neighborhood of $\zeta$.
Without loss of generality, it is assumed that
$Z(\cdot,m)>0$ in $(\varphi_0,\zeta)$,
where $\varphi_0$ is given above.
Let $\tilde\Sigma_+$ and $\tilde\Sigma_-$ be the above characteristics.
Denote
\begin{gather*}
\tilde{\mathscr P}_+(\varphi)=Q(\varphi,\tilde\Psi_+(\varphi)),\quad
\tilde{\mathscr W}(\varphi)=W(\varphi,\tilde\Psi_+(\varphi)),\quad
\varphi\ge\varphi_0,
\\
\tilde{\mathscr P}_-(\varphi)=Q(\varphi,\tilde\Psi_-(\varphi)),\quad
\tilde{\mathscr Z}(\varphi)=Z(\varphi,\tilde\Psi_-(\varphi)),\quad
\varphi\ge\tilde\varphi.
\end{gather*}
Then, \eqref{wz-1} and \eqref{wz-2} yield
\begin{align}
\label{a8-7-1}
(h(\tilde{\mathscr P}_+(\varphi))\tilde{\mathscr W}(\varphi))'=
\frac14b^{-1}(\tilde{\mathscr P}_+(\varphi))p(\tilde{\mathscr P}_+(\varphi))
h(\tilde{\mathscr P}_+(\varphi))\tilde{\mathscr W}^2(\varphi),\quad
\varphi>\varphi_0
\end{align}
and
\begin{gather}
\label{a8-7-0}
\tilde{\mathscr P}'_-(\varphi)=W(\varphi,\tilde\Psi_-(\varphi)),\quad
\varphi>\tilde\varphi,
\\
\label{a8-7}
(h(\tilde{\mathscr P}_-(\varphi))\tilde{\mathscr Z}(\varphi))'=
-\frac14b^{-1}(\tilde{\mathscr P}_-(\varphi))p(\tilde{\mathscr P}_-(\varphi))h(\tilde{\mathscr P}_-(\varphi))
\tilde{\mathscr Z}^2(\varphi),\quad
\varphi>\tilde\varphi,\quad\tilde{\mathscr Z}(\tilde\varphi)>0,
\end{gather}
where
$$
h(s)=\mbox{exp}\Big\{\frac14\int_{-1}^sb^{-1}(t)p(t)dt\Big\},\quad s<0.
$$
First estimate $\tilde{\mathscr W}$.
It follows from \eqref{a8-7-1} that
\begin{align}
\label{a8-7-2}
h(\tilde{\mathscr P}_+(\varphi))\tilde{\mathscr W}(\varphi)\ge
h(\tilde{\mathscr P}_+(\varphi_0))\tilde{\mathscr W}(\varphi_0),\quad
\varphi>\varphi_0.
\end{align}
Note
$\inf_{0<\tilde\psi<\psi_0}h(\tilde{\mathscr P}_+(\varphi_0))\tilde{\mathscr W}(\varphi_0)>-\infty$.
Then, \eqref{a8-7-2} and (i) yield
\begin{align}
\label{a8-10-1}
\lim_{\overline{\psi\to m^-}}
\min_{\zeta_1\le\varphi\le\zeta_2}h(Q(\varphi,\psi))W(\varphi,\psi)>-\infty,
\quad
\zeta<\zeta_1<\zeta_2.
\end{align}
Turn to $\tilde{\mathscr Z}$.
By \eqref{a8-6}, the second limit in \eqref{a8-8} and $1<\gamma<3$, there exists
a positive constant $M_3=M_3(\gamma,\mu_1)$ such that
\begin{align}
\label{a8-9}
\frac14b^{-1}(Q(\varphi,\psi))p(Q(\varphi,\psi))h^{-1}(Q(\varphi,\psi))\le
M_3(-Q(\varphi,\psi))^{(\gamma-3)/8},
\quad&&(\varphi,\psi)\in(0,+\infty)\times(0,m).
\end{align}
It follows from \eqref{a8-6}, \eqref{a8-7-0}, \eqref{a8-7} and \eqref{a8-9} that
\begin{gather*}
\tilde{\mathscr P}_-(\varphi)\le\tilde{\mathscr P}_-(\tilde\varphi),\quad
\tilde{\mathscr Z}(\varphi)>0,\quad\varphi>\tilde\varphi,
\\
(h(\tilde{\mathscr P}_-(\varphi))\tilde{\mathscr Z}(\varphi))'
\ge-M_3(-\tilde{\mathscr P}_-(\tilde\varphi))^{(\gamma-3)/8}
(h(\tilde{\mathscr P}_-(\varphi))\tilde{\mathscr Z}(\varphi))^2,\quad
\varphi>\tilde\varphi
\end{gather*}
and thus
\begin{align*}
h(\tilde{\mathscr P}_-(\varphi))\tilde{\mathscr Z}(\varphi)
\ge\frac{1}{(h(\tilde{\mathscr P}_-(\tilde\varphi))\tilde{\mathscr Z}(\tilde\varphi))^{-1}
+M_3(-\tilde{\mathscr P}_-(\tilde\varphi))^{(\gamma-3)/8}(\varphi-\tilde\varphi)},
\quad
\varphi\ge\tilde\varphi.
\end{align*}
Note that \eqref{a8-6-1}, the second limit in \eqref{a8-8} and (i) give
${\lim_{\tilde\varphi\to\zeta^-}}h(\tilde{\mathscr P}_-(\tilde\varphi))\tilde{\mathscr Z}(\tilde\varphi)
=+\infty$
and $\lim_{\tilde\varphi\to\zeta^-}\tilde{\mathscr P}_-(\tilde\varphi)=-\infty$.
Therefore,
\begin{align}
\label{a8-10}
{\lim_{\psi\to m^-}}\min_{\zeta_1\le\varphi\le\zeta_2}h(Q(\varphi,\psi))Z(\varphi,\psi)=+\infty,
\quad
\zeta<\zeta_1<\zeta_2.
\end{align}
Combine \eqref{a8-10-1} and \eqref{a8-10} to get
\begin{align*}
\lim_{{\tilde\varphi\to\zeta^-}}\max_{\zeta_1\le\varphi\le\zeta_2}h(Q(\varphi,\psi))
b^{1/2}(Q(\varphi,\psi))Q_\psi(\varphi,\psi)=-\infty,
\quad\zeta<\zeta_1<\zeta_2,
\end{align*}
which, together with the second limit in \eqref{a8-8}, yields
\begin{align}
\label{a8-11}
\lim_{{\tilde\varphi\to\zeta^-}}\min_{\zeta_1\le\varphi\le\zeta_2}
(c^*-q(\varphi,\psi))^{(3-\gamma-\varepsilon)/(4\gamma-4)}
q_\psi(\varphi,\psi)=+\infty,
\quad\zeta<\zeta_1<\zeta_2
\end{align}
for each $\varepsilon>0$.
Then, (v) is proved by transforming \eqref{a8-11} into the physical plane.
$\hfill\Box$\vskip 4mm

\subsection{Smooth supersonic flows with vacuum at the inlet}

For an incoming flow with vacuum, the smooth supersonic flow problem can be formulated as follows
\begin{align}
\label{aint-phy1}
&\mbox{div}(\rho(|\nabla\varphi|^2)\nabla\varphi)=0,\quad&&(x,y)\in\Omega,
\\
\label{aint-phy2}
&\rho(q_0^2(y))\varphi(\Upsilon(y),y)=0,\quad&&0<y<f(l_0),
\\
\label{aint-phy3}
&|\nabla\varphi(\Upsilon(y),y)|=q_0(y),\quad&&0<y<f(l_0).
\\
\label{aint-phy4}
&\rho(|\nabla\varphi(x,0)|^2)\pd{\varphi}y(x,0)=0,\quad&&\Upsilon(0)<x<l_1,
\\
\label{aint-phy5}
&\rho(|\nabla\varphi(x,f(x))|^2)\Big(\pd{\varphi}y(x,f(x))-f'(x)\pd{\varphi}x(x,f(x))\Big)=0,
\quad&&l_0<x<l_1.
\end{align}

Assume that $q_0\in C([0,f(l_0)])\cap
C^{1}([0,y_1])\cap C^{1}([y_2,f(l_0)])$ with $0<y_1<y_2<f(l_0)$
satisfies
\begin{align}
\label{aint-1}
q_0(y)\left\{
\begin{aligned}
&\in(c_*,c^*),\quad&&y\in [0,y_1)\cup (y_2,f(l_0)],
\\
&=c^*,\quad&&y\in [y_1,y_2]
\end{aligned}
\right.
\end{align}
and
\begin{align}
\label{aint-2}
|q'_0(y)|\le\frac{-\Upsilon''(y)}{1+(\Upsilon'(y))^2}\,
\sqrt{\frac{-q_0^2(y)\rho(q_0^2(y))}{\rho(q_0^2(y))+2q_0^2(y)\rho'(q_0^2(y))}}\,,
\quad
y\in [0,y_1)\cup (y_2,f(l_0)].
\end{align}
Similar to Proposition \ref{phys-thm71},
one can solve the smooth supersonic flow problems with free boundary on $[0,y_1]$ and $[y_2,f(l_0)]$
individually.
The two free boundaries are
$$
\tilde\Gamma_k: y=y_k-\Upsilon'(y_k)(x-\Upsilon(y_k)),\quad x\ge\Upsilon(y_k),
$$
for $k=1,2$.
Then, one can get a global smooth supersonic flow
to the problem \eqref{aint-phy1}--\eqref{aint-phy5} by
connecting the two smooth supersonic flows with free boundary by the vacuum.

\vskip5mm

\hskip40mm
\setlength{\unitlength}{0.6mm}
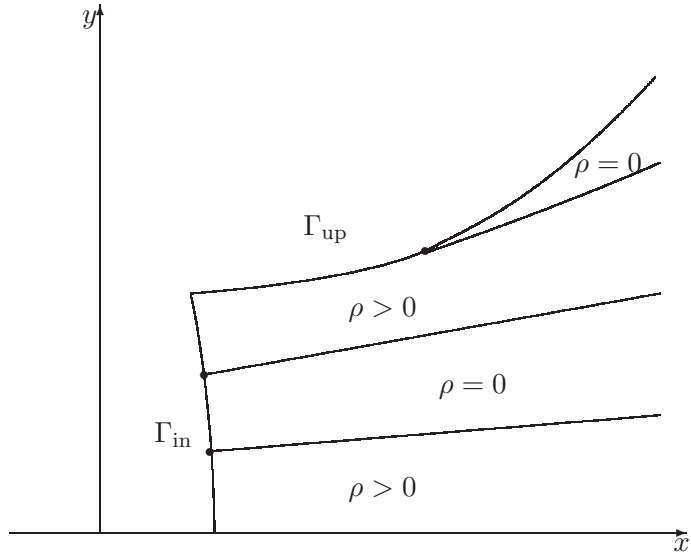
\begin{picture}(250,120)
\put(-20,0){\vector(1,0){150}}
\put(0,0){\vector(0,1){117}}
\put(127,-4){$x$} \put(-4,113){$y$}

\put(60,48){\cbezier(-40,5)(-20,6.5)(5,10)(15,16)}
\put(60,48){\cbezier(15,16)(25,21)(40,28)(63,53)}

\put(60,0){\cbezier(-40,53)(-37.5,42)(-34.8,20)(-34.7,0)}

\put(45,66){$\Gamma_{\text{\rm up}}$}
\put(12,20){$\Gamma_{\text{\rm in}}$}

\put(55,48){$\rho>0$}

\put(105,80){$\rho=0$}
\put(72,62.5){\circle*{1.5}}
\put(75,55){\qbezier(-3,7)(26,17)(49,27)}

\put(23,35){\circle*{1.5}}

\put(24.3,18){\circle*{1.5}}

\put(0,0){\qbezier(23,35)(73.5,44)(124,53)}

\put(0,0){\qbezier(24.3,18)(74.15,22)(124,26)}

\put(75,31){$\rho=0$}

\put(55,8){$\rho>0$}

\end{picture}
\vskip10mm

\begin{theorem}
Assume that $f\in C^{2}([l_0,l_1))$ satisfies \eqref{sss-0},
$\Upsilon\in C^{2}([0,f(l_0)])$ satisfies \eqref{xc4} and
$q_0\in C([0,f(l_0)])\cap
C^{1}([0,y_1])\cap C^{1}([y_2,f(l_0)])$
satisfies \eqref{aint-1} and \eqref{aint-2}.
Then the problem \eqref{aint-phy1}--\eqref{aint-phy5} admits uniquely a global smooth solution
$\varphi\in C^{1}(\overline\Omega)$ with $|\nabla\varphi|\in C^{0,1}(\overline\Omega)$.
Besides the possible vacuum enclosure attached to the upper wall,
there is another inner vacuum region which is
the closed domain bounded by $\Gamma_{\text{\rm in}}$, $\tilde\Gamma_1$ and $\tilde\Gamma_2$.
Moreover, $|\nabla\varphi|\in C^{1}(\Omega)$,
$\pd{|\nabla\varphi|}{x}\Big|_{{\tilde\Gamma_k}{\cap}{\Omega}}=\pd{|\nabla\varphi|}{y}\Big|_{{\tilde\Gamma_k}{\cap}{\Omega}}=0$
for $k=1,2$.
\end{theorem}

\begin{remark}
For the problem \eqref{aint-phy1}--\eqref{aint-phy5},

{\rm (i)} If $0<y_1=y_2<f(l_0)$, then $\tilde\Gamma_1$ and $\tilde\Gamma_2$ coincide,
and thus the inner vacuum region of the global smooth supersonic flow
is degenerate to a half-line.

{\rm (ii)} If $y_1=0$, then $\tilde\Gamma_1$ is the lower wall of the nozzle
and the global smooth supersonic flow is vacuum on the lower wall.

{\rm (iii)} If $y_2=f(l_0)$, then the vacuum enclosure attached to the upper wall
and the inner vacuum region are connected, and
thus the whole vacuum region is
the closed domain bounded by $\Gamma_{\text{\rm in}}$, $\Gamma_{\text{\rm up}}$ and $\tilde\Gamma_1$.
\end{remark}

\end{document}